\documentclass[11pt]{article}
\usepackage{graphicx}%
\usepackage{multirow}%
\usepackage{amsmath,amssymb,amsfonts}%
\usepackage{amsthm}%
\usepackage{mathrsfs}%
\usepackage[title]{appendix}%
\usepackage{xcolor}%
\usepackage{textcomp}%
\usepackage{manyfoot}%
\usepackage{booktabs}%
\usepackage{algorithm}%
\usepackage{algorithmicx}%
\usepackage{algpseudocode}%
\usepackage{listings}%
\usepackage{amsmath,amsfonts}
\usepackage{txfonts,wasysym,lscape,amssymb,mathrsfs,extarrows}
\usepackage[all,pdf]{xy}
\usepackage{graphicx}
\usepackage{subfigure}
 \usepackage{multirow}
\usepackage{booktabs}
\usepackage{color}

\newtheorem{theorem}{Theorem}[section]
\newtheorem{proposition}{Proposition}[section]
\newtheorem{lemma}{Lemma}[section]
\newtheorem{corollary}{Corollary}[section]

\newtheorem{remark}{Remark}[section]

\renewcommand{\theequation}{\thesection.\arabic{equation}}

\title{Stochastic Approximation Proximal Subgradient Method for
Stochastic Convex-Concave Minimax Optimization
}

\author{Yu-Hong Dai\footnote{LSEC, ICMSEC, AMSS, Chinese Academy of Sciences, Beijing 100190, China. {\sl  Email}: dyh@lsec.cc.ac.cn.
This author was supported by the Natural Science Foundation of China (Nos. 11991021, 12021001, 11991020 and 11971372) and the Strategic Priority Research Program of Chinese Academy of Sciences (No. XDA27000000).}\quad Jiani Wang\footnote{School of Science, Beijing University of Posts and Telecommunications, Beijing 100876,  China. {\sl  Email}: wjiani@lsec.ac.cc.cn }\quad and \quad Liwei Zhang \footnote{National Frontiers Science Center for Industrial Intelligence and Systems Optimization, Northeastern University, Shenyang 110819, China; Key Laboratory of Data Analytics and Optimization for Smart Industry  (Northeastern University),  Ministry of Education, Shenyang 110819, China. {\sl Email}: Zhanglw@ mail.neu.edu.cn. This author was supported by the Major Program of National Natural Science Foundation of China (Nos. 72192830 and 72192831), National Natural Science Foundation of China (No.12371298) and the 111 Project (B16009)}}

\date{}
\begin{document}

\maketitle

\begin{abstract}
This paper  presents a stochastic approximation proximal subgradient (SAPS) method for stochastic convex-concave minimax optimization. By accessing unbiased and variance bounded approximate subgradients, we show that this algorithm exhibits  ${\rm O}(N^{-1/2})$  expected convergence rate of the minimax optimality measure if the parameters in the algorithm are properly chosen,  where $N$ denotes the number of iterations.  Moreover, we show that the algorithm has  ${\rm O}(\log(N)N^{-1/2})$  minimax optimality measure bound with high probability. Further we study a  specific stochastic convex-concave minimax optimization problems arising from stochastic convex conic optimization problems, which the the bounded subgradient condition is fail. To overcome the lack of the bounded subgradient conditions in convex-concave minimax problems, we propose a linearized stochastic approximation augmented Lagrange (LSAAL) method and prove that this algorithm exhibits  ${\rm O}(N^{-1/2})$ expected convergence rate for the minimax optimality measure and ${\rm O}(\log^2(N)N^{-1/2})$  minimax optimality measure bound with high probability as well.  Preliminary numerical results demonstrate the effect of the SAPS and LSAAL methods.
\vskip 6 true pt \noindent \textbf{Key words}: stochastic  convex-concave minimax optimization, stochastic convex conic optimization, stochastic approximation, proximal  point method, linearized stochastic approximation augmented Lagrange  method, expected convergence rate, high probability bound.
\vskip 12 true pt \noindent \textbf{AMS subject classification}: 90C30
\end{abstract}
\bigskip\noindent

\section{Introduction}
Consider the following stochastic minimax optimization
\begin{equation}\label{minimax}
\min_{x \in \Re^n} \max_{y \in \Re^m}\Big\{\phi(x,y)=\vartheta (x)+\mathbb E[F(x,y,\xi)]-\omega (y)\Big\},
\end{equation}
where $\vartheta:\Re^n \rightarrow \overline \Re$ and  $\omega:\Re^m\rightarrow \overline \Re$ are proper lower semicontinuous convex functions,  $\xi$ is a random vector whose probability distribution   is supported on  $\Xi \subseteq \Re^q$ and $F: \Re^n \times \Re^m \times \Xi \rightarrow \Re$ is a real-valued function.
Let us denote the expected value function
$$
f(x,y):=\mathbb{E}[F(x,y,\xi)]= \int_{\Xi} F(x,y,\xi)dP(\xi).
$$
Assume that this function $f(\cdot,\cdot)$ is well defined and finite valued 
for every $(x,y)\in {\rm dom}\,\vartheta \times {\rm dom}\,\omega$.  Meanwhile, assume that  $f(\cdot,\cdot)$ is continuous and convex-concave on ${\rm dom}\,\vartheta \times {\rm dom}\,\omega$, but  may not be differentiable (this clearly holds if for every
$\xi \in \Xi$, the function $F(\cdot,\cdot,\xi)$ is convex-concave on ${\rm dom}\,\vartheta \times
{\rm dom}\,\omega$).
 Under the above assumptions, the problem (\ref{minimax}) becomes a convex-concave minimax problem.

\subsection{Special Cases}
As a stochastic minimax problem with a general convex-concave structure,  the problem \eqref{minimax} has a wide range of applications such as stochastic convex programming \cite{Duchi2015}, linear regression \cite{Du2017}, robust optimization \cite{Ben2009} and adversarial generative networks \cite{Good14}. Below we list three special cases of the problem \eqref{minimax},  where $\vartheta$ and $\omega$ are defined as different forms.

\begin{description}
  \item[Minimax problems over convex sets.] If $\vartheta(x)=\delta_X(x)$ and $\omega (y)=\delta_Y(y)$, where $X \subset \Re^n,\ Y \subset \Re^m$ are two convex sets and $\delta_X,\ \delta_Y$ are the indicator functions, the problem (\ref{minimax}) is reduced to
 \begin{equation}\label{minimaxXY}
\min_{x \in X} \max_{y \in Y}\Big\{\phi(x,y)=\mathbb E[F(x,y,\xi)]\Big\}.
\end{equation}
 Many machine learning problems such as reinforcement learning (\cite{Wai18, Wai19}),
black-box adversarial attack (\cite{Liu19b, Wang2020, Xu2020}) and adversarial training (\cite{Good14, Liu19a}) can be expressed as the minimax problem \eqref{minimaxXY}.
  \item[Minimax problems with regularized functions.] Another case is the minimax problem with the functions $\vartheta(x)=r_x(x)$ and $\omega (y)=r_y(y)$, where $r_x$ and $r_y$ are regularized functions; namely,
      \begin{equation}\label{minimaxregularizedf}
\min_{x \in \Re^n} \max_{y \in \Re^m}\Big\{\phi(x,y)=r_x(x)+\mathbb E[F(x,y,\xi)]-r_y(y)\Big\}.
\end{equation}
To solve the overfitting problem while making the model sparse or low-rank, some regularized function is usually added to the objective function. However, many regularized functions are non-differentiable, such as $\ell_0$-norm $\|\cdot\|_0$, $\ell_1$-norm $\|\cdot\|_1$, $\ell_2$-norm $\|\cdot\|_2$, which are common in regression analysis in statistics \cite{Domingos2000}, model training in machine learning \cite{Dietterich1995,  Kolluri2020, Ying2019} and sparse two-matrix game problems \cite{Carmon2020}.
  \item[Stochastic convex conic optimization.] Consider the stochastic convex conic optimization problem
 \begin{equation}\label{CP}
 \begin{array}{ll}
 \displaystyle\min_{x \in X} & f(x)=\mathbb E [F(x,\xi)]\\[8pt]
 {\rm s.\ t.} & g(x)=\mathbb E [G(x,\xi)] \in {\cal K}.
  \end{array}
 \end{equation}
In the above, $\xi$ is a random vector whose probability distribution   is supported on  $\Xi \subseteq \Re^q$, $F: {\cal O} \times \Xi \rightarrow \Re$ is a  real-valued function, $G:{\cal O} \times \Xi \rightarrow {\cal Y}$ is a mapping for an open convex set ${\cal O} \supset X$, $X \subset \Re^n$ is a nonempty convex compact set, ${\cal K} \subset {\cal Y}$ is a closed convex cone and ${\cal Y}$ is a finite-dimensional Hilbert space. The conjugate dual of the stochastic convex programming problem (\ref{CP}) is defined as
\begin{equation}\label{CD}
         \max_{y}\inf_{x\in X}\big\{ l(x,y)-\delta^*_{{\cal K}}(y)\big\},
        \end{equation}
where the Lagrangian function for (\ref{CP}) is defined by $l(x,y)=f(x)+\langle y, g(x)\rangle$ for any $(x, y) \in {\cal O}\times {\cal Y}$ and $\delta^*_{{\cal K}}$ is the conjugate function of the indicator function $\delta_{{\cal K}}$. Under some regularity conditions (see Section 3 for more details), the conjugate dual  (\ref{CD}) can be expressed as a convex-concave minimax problem of the form \eqref{minimax} with functions $\vartheta$ and $\omega$ being indicator functions related to the constraints and it has the same optimal value as the stochastic convex programming problem (\ref{CP}).
\end{description}

\subsection{Motivation and Contributions}
One difficulty in solving the problem (\ref{minimax}) is that  the probability distribution function may not be available. Even  if the distribution function is easy to obtain, the expectation with respect to $\xi$ may be difficult to calculate within a high accuracy in the large scale case. In order to overcome this difficulty, a popular approach is to utilize the stochastic approximation (SA) technique, where  approximations can be accessed via calls to stochastic oracles considered as the noisy computable version of
the ``real'' function.  To this aim, the following general assumptions are used throughout this paper.
\begin{itemize}
\item[({\bf A1})] The samples $\xi_1,\xi_2,\ldots$ of realizations of random vector $\xi$ are
 generated by independent identical distribution (i.i.d.) ;
\item[({\bf A2})] First-order oracles are unbiased estimators of the subgradient of $F(x,y,\xi)$; namely, for any point $(x,y,\xi)\in {\rm dom}\,\vartheta \times {\rm dom}\,\omega \times \Xi$, return a stochastic subgradient
    $$
    G(x,y,\xi)=\left(
    \begin{array}{l}
    G_x(x,y,\xi)\\[3pt]
    -G_y(x,y,\xi)
    \end{array}
    \right
    )
    $$
     such that
          $$
    g(x,y)=\left(
    \begin{array}{l}
    g_x(x,y)\\[3pt]
    -g_y(x,y)
    \end{array}
    \right
    )=\left(
    \begin{array}{l}
   \mathbb E[G_x(x,y,\xi)]\\[3pt]
    -\mathbb E[G_y(x,y,\xi)]
    \end{array}
    \right
    )
    $$
   is well-defined, where $g_x(x,y)\in \partial_x f(x,y)$ and $-g_y(x,y) \in \partial_y[-f(x,y)]$.
\end{itemize}

For every $\xi \in \Xi$, if the function $F(\cdot,\cdot,\xi)$ is convex-concave and its
respective subdifferential and integral operators are interchangeable, we can ensure ({\bf A2})
by setting
$$
G(x,y,\xi)=\left(
    \begin{array}{l}
   G_x(x,y,\xi)\\[3pt]
    -G_y(x,y,\xi)
    \end{array}
    \right
    )\in \left(
    \begin{array}{l}
   \partial_x F(x,y,\xi)\\[3pt]
    - \partial_y F(x,y,\xi)
    \end{array}
    \right
    ).$$

In a pioneering and profound work \cite{Lan2009}, Nemirovski {\it et al.} proposed the robust stochastic approximation (RSA) method for stochastic convex optimization. This method stimulates the development of stochastic approximation algorithms in machine learning  \cite{Bottou2018, Davis2019, Gower19, Levy20, Li19, Pu21, Schmidt2017} and now is widely used in distributed deep learning \cite{Ben-Nun2019}, multiple access channel \cite{Sery2020} and  low-rank matrix factorization \cite{chi2019}.
However, to the best of our knowledge, there is no algroithm for solving stochastic general convex-concave non-differentiable minimax optimization (\ref{minimax}). Therefore, we shall propose  the stochastic approximation proximal subgradient method  for the problem (\ref{minimax}), which is a fundamental idea that extends the algorithm from stochastic convex optimization to stochastic general convex-concave minimax problems.

The main contributions of this paper are as follows.
\begin{itemize}
  \item For the stochastic minimax problem (\ref{minimax}) with lower semicontinuous convex items $\vartheta$ and $\omega$, we design a stochastic approximation proximal subgradient (SAPS) method and verify the sublinear convergence rates of the method for general convex-concave case, where the estimators are generated by unbiased first-order oracles with bounded variance.  Moreover, SAPS can obtain ${\rm O}(N^{-1/2}\log(N))$ minimax optimality measure bound with high probability under the bounded subgradient condition, where $N$ is the number of iterations. To the  best of our knowledge, \emph{SAPS is the first stochastic algorithm  for  solving the stochastic general convex-concave non-differentiable minimax problem (\ref{minimax}).}

  \item  A linearized stochastic approximation augmented Lagrange (LSAAL) method is developed for the stochastic convex-concave minimax problem \eqref{CD} arising from stochastic convex conic constrained programming \eqref{CP}, in which case the bounded subgradient condition dose not satisfied.  Under
mild conditions,  the optimal value of the stochastic convex cone constrained
optimization is equivalent to its dual problem. The sublinear convergence of the LSAAL method is proved with respect to the expected minimax optimal measure of the augmented Lagrangian function of (1.3). Meanwhile, we show that LSAAL exhibits ${\rm O}(N^{-1/2}\log^2(N))$ minimax optimality measure bound with high probability as well. It is worth noting that \emph{LSAAL can be widely used to solve stochastic convex cone constrained programming including nonlinear programming, second-order cone optimization and semidefinite programming.}

      \item We  test the effect of the SAPS method for  strongly convex-concave case and general convex-concave case and conduct numerical experiments of the LSAAL method for the multi-class
Neyman-Pearson classification on four
real data sets. Numerical results reveal promising performances of the SAPS and LSAAL methods.
\end{itemize}

\subsection{Related Work}
\begin{description}

  \item[Stochastic minimax problems] Many methods for stochastic minimization problems are proposed to solve the stochastic minimax problem  (\ref{minimaxXY}), such as zero-order (namely, gradient-free) algorithms \cite{Akhavan2021, Dvinskikh2022, Xu2020, Xu2021, Xu2023} and first-order algorithms \cite{HuangZ2022, LinT20, Luo21, Lan2009, Xian21}. First-order stochastic algorithms have been extended to the case where function $\mathbb E[F(x,y,\xi)]$ is continuous differentiable. Zhang \emph{et al.} \cite{ZAG21} proposed a  stochastic accelerated primal-dual algorithm for strongly-convex-strongly-concave saddle point problems and established the  linearly convergence to a neighborhood of the unique saddle point.
      Yang \emph{et al.}   \cite{Yang2020} proposed the stochastic alternating gradient descent ascent (SAGDA) algorithm for a subclass of the nonconvex-nonconcave minimax problem, and showed a sublinear rate of SAGDA under the two-sided Polyak-Lojasiewicz condition.
       In \cite{Q20}, Tran-Dinh \emph{et al.} considered the following stochastic nonconvex-concave minimax problem
       \begin{equation}\label{minimaxTD}
\min_{x \in \Re^n} \max_{y \in \Re^m}\Big\{\phi(x,y)=\vartheta (x)+\langle Ky,\ \mathbb E[F(x,\xi)]\rangle-\omega (y)\Big\},
\end{equation}
where the  randomness of the objective function is only on the variable $x$ and $\mathbb E[F(x,\xi)]$ is continuous differentiable. They developed a single-loop variance-reduced algorithm and  achieved the  convergence rate of order $O(N^{-2/3})$.

      Many  stochastic zero-order and subgradient algorithms are used to solve the case where function $\mathbb E[F(x,y,\xi)]$ in the problem (\ref{minimaxXY}) is non-differentiable. Under conditions {\bf (A1)}-{\bf (A2)},  Nemirovski \emph{et al.} \cite{Lan2009} proposed a robust stochastic approximation approach for solving the problem (\ref{minimaxXY}) when $X$ and $Y$ are nonempty convex compact sets. They proved that the expected
convergence in terms of an minimax optimality measure is of the order ${\rm O}(N^{1/2})$. Moreover, they established a high probability guarantee of the algorithm. The RSA approach in  \cite{Lan2009} is in fact a projection gradient method for solving the problem (\ref{minimaxXY})
combining the averaging technique  developed by \cite{Polyak1990} and \cite{PJ1992}. The epoch-wise stochastic gradient descent ascent method is an extension of the epoch gradient descent method for solving strongly-convex-strongly-concave minimax problems in \cite{YY20},
and achieved the optimal rate of $O(1/T)$ for the duality gap. The ZO-Min-Max framework was proposed by \cite{Liu19b}, where
   Liu \emph{{et al.}} integrated the zeroth-order gradient estimator with an alternating projected stochastic gradient descent-ascent method, and obtained the sublinear convergence rate and scales with problem size.  In \cite{Xu2020}, for the nonconvex-strongly-concave minimax stochastic optimization,  a zeroth-order variance reduced gradient descent ascent (ZO-VRGDA) algorithm was proposed, and ZO-VRGDA achieved the iteration  complexity of order $O(\epsilon^{-3})$ for finding an $\epsilon$-stationary point. A class of accelerated zeroth-order and first-order momentum
methods  were analyzed in \cite{Huang2022} for both nonconvex minimization and minimax  optimization of the form (\ref{minimaxXY}). Nevertheless, to the best of our knowledge, there is no algorithm for solving the stochastic general convex-concave minimax problem (\ref{minimax}), which is a more general form including the problems (\ref{minimaxXY}) and (\ref{minimaxTD}). This promotes us to explore the properties of the problem (\ref{minimax}) and develop a fast stochastic algorithm for solving (\ref{minimax}).

  \item[Stochastic convex programming problems] Currently, there have been some methods to solve the stochastic nonlinear programming, such as \cite{Duchi2015, Mahdavi(2011), Mehrdad(2013),  ZZXW22}.  Mahdavi \emph{et al.} \cite{Mahdavi(2011)} studied an online gradient descent method for the online convex optimization and obtained ${\rm O}(N^{-1/2})$ objective regret and ${\rm O}(N^{-1/4})$ constraint violation by analyzing the convergence of the optimal value depending on the sufficiently large sample and the boundedness of the gradient of the constraint functions. Moreover, they proposed a stochastic primal-dual algorithm for stochastic nonlinear programming  in \cite{Mehrdad(2013)} and attained the optimal convergence rate of ${\rm O}(N^{-1/2})$ for the general Lipschitz continuous objective functions with high probability. Zhang \emph{et al.} in \cite{ZZXW22} designed a proximal
point method for solving convex stochastic nonlinear programming, which had no more than ${\rm O}(N^{-1/4})$ objective regret and no more than ${\rm O}(N^{-1/8})$ constraint violation with high probability. In this paper, stochastic convex conic optimization problems, which  contain the convex nonlinear programming, are solved by reformulating as convex-concave minimax problems of the form (\ref{minimax}).  However, for these stochastic convex-concave minimax problems, $Y$ is an unbounded set and hence the bounded gradient condition, which is normally required for the convergence analysis of RSA, does not hold. This means that we can not use RSA to solve stochastic convex programming problems. Such an observation motivates us to put forward the linearized stochastic approximation augmented Lagrange method, which is based on the augmented Lagrange duality of the convex programming.
\end{description}

%


\subsection{Organization and notations}
{\bf Organization.} The remainder  of this paper is organized as follow. In Section
 \ref{Sec:3}, we propose the stochastic approximation proximal subgradient method  for solving the problem (\ref{minimax}), and prove the convergence rates for general convex-concave minimax problems. In Section \ref{Sec:4}, we  put forward the linearized stochastic approximation augmented Lagrange method for solving stochastic convex conic optimization problems. In  Section \ref{Sec:5}, we  report our numerical results of the SAPS and LSAAL methods. Finally, we draw some discussion in Section \ref{Sec:6}.\\[2mm]
{\bf Notations.} Throughout the paper, we use the following notations. By $\|x\|$, we denote  the Euclidean norm of vector
$x \in \Re^n$; namely, $\|x\|=\sqrt{x^Tx}$. Let ${\cal Y}$ be a finite-dimensional Hilbert space. The inner product in ${\cal Y}$ is defined as $\langle \cdot,\cdot \rangle$ and the induced norm of vector $y \in {\cal Y}$ is $\|y\|=\sqrt{\langle y, y\rangle}$.
 For a lower semicontinuous convex function $\theta:{\cal Y} \rightarrow \overline \Re$, by ${\rm Prox}_{\gamma \theta}$ for $\lambda>0$, we denote the proximal mapping of $\theta$; namely,
$
{\rm Prox}_{\gamma \theta}(y)={\rm argmin}\left\{ \theta (w)+\displaystyle \frac{1}{2\gamma}\|w-y\|^2\right\}.
$
It follows from Page 878 of  \cite{Rockafellar76a} that the proximal mapping  ${\rm Prox}_{\gamma \theta}$ is a nonexpanding
operator; namely,
$
\|{\rm Prox}_{\gamma \theta}(y')-{\rm Prox}_{\gamma \theta}(y)\|\leq \|y'-y\|$ for all $y',y \in {\cal Y}.
$
The conjugate function of $\theta$, denoted by $\theta^*:{\cal Y} \rightarrow \overline \Re$, is defined by
$
\theta^*(v)=\sup_{y \in {\cal Y}}\Big\{\langle v, y\rangle -\theta (y) \Big\}.
$
For a nonempty close convex set $D \subseteq {\cal Y}$, we use $\delta_D$ to denote the indicator function
of $D$,
$$
\delta_D(y)=\left
\{
\begin{array}{ll}
0, & y \in D;\\[4pt]
+\infty, & y \notin D.
\end{array}
\right.
$$
The conjugate of $\delta_D$ is the supporting function of $D$; namely,
$
\delta^*_D(v)=\sup_{y \in D} \langle v, y\rangle.
$
The proximal mapping of $\delta_D$ is the metric projection operator onto
the set $D$; namely,  ${\rm Prox}_{\gamma \delta_D}(y)=\Pi_D(y)$, where
$
\Pi_D(y) ={\rm argmin}\{\|y'-y\|: y' \in D\}.
$
 Obviously we have that  $\Pi_D$ is a nonexpanding
operator; namely,
$
\|\Pi_D(y')-\Pi_D(y)\|\leq \|y'-y\|$ for all $y',y \in {\cal Y}.
$
For a closed convex cone ${\cal K}\subset {\cal Y}$, any $y \in {\cal Y}$ has the following decomposition
$
y=\Pi_{K}(y)+\Pi_{{\cal K}^{\circ}}(y)
$
with $\langle \Pi_{K}(y),\Pi_{{\cal K}^{\circ}}(y)\rangle=0$, where ${\cal K}^{\circ}$ is the polar cone of ${\cal K}$.
 \section{Stochastic Approximation Proximal Subgradient Method }\label{Sec:3}
 \setcounter{equation}{0}
 In this section,  we propose a stochastic approximation proximal subgradient method  for solving the convex-concave minimax problem (\ref{minimax}) and discuss the convergence results of the iterations with general convex-concave structures. Particularly, the sublinear convergence rate with expectation is established for the general convex-concave case if the approximate subgradients of the function $\phi$ are unbiased provided with bounded variance. Ulteriorly, the convergence rate of SAPS with high probability can be obtained when the approximate subgradients are bounded.

Define
$
 z=(x,y) \ \mbox{for any}\ z\in \mathbb Z={\rm dom}\,\vartheta \times {\rm dom}\,\omega.
$
 We give the existence assumption of the minimax point of the problem (\ref{minimax}).
\begin{itemize}
\item[({\bf A3})]There exists a point $z^*=(x^*,y^*) \in {\rm dom}\,\vartheta \times {\rm dom}\,\omega$ such that
\begin{equation}\label{eq:SaddleP}
\phi (x^*,y) \leq \phi (x^*,y^*) \leq \phi (x,y^*), \quad \forall (x,y) \in {\rm dom}\,\vartheta \times {\rm dom}\,\omega.
\end{equation}
\end{itemize}

For the minimax point $z^*=(x^*,y^*)$ defined as \eqref{eq:SaddleP},  the optimality conditions of the problem (\ref{minimax}) at $z^*=(x^*,y^*)$ are expressed as
 \begin{equation}\label{eq:opt-minmax}
 \left\{
 \begin{array}{ll}
 0\in g_x (z^*)+\partial \vartheta (x^*), & 0\in -g_y (z^*)+\partial \omega (y^*);\\[6pt]
  g_x(z^*)\in \partial_x f(z^*),& -g_y(z^*)\in \partial_y [-f](z^*),
  \end{array}
  \right.
 \end{equation}
 which are equivalent to the following equalities
  \begin{equation}\label{eq:opt-minmaxE}
 x^*={\rm Prox}_{\gamma \vartheta}(x^*-\gamma g_x (z^*)), \ \  y^*={\rm Prox}_{\gamma \omega}(y^*+\gamma g_y (z^*)).
 \end{equation}
We introduce the function $\psi: \Re^n \times \Re^m \rightarrow \overline \Re$ by
$
\psi (z)=\vartheta (x)+\omega(y).
$
Then for $z =(x,y)$, we have
 $$
 {\rm Prox}_{\gamma \psi}(z)= ({\rm Prox}_{\gamma \psi}(x),{\rm Prox}_{\gamma \psi}(y)).
 $$
 With this notation,  the optimal conditions (\ref{eq:opt-minmaxE}) can be expressed as
 \begin{equation}\label{eq:compact-zop}
z^*= {\rm Prox}_{\gamma \psi}(z^*-\gamma g(z^*)),
 \end{equation}
where $g(z^*)=(g_x(z^*);-g_y(z^*))$.

Suppose that the first-order oracles ${G}$ are $Carath\acute{e}odory$ functions\footnote{It is said that the function ${G}: {\rm dom}\,\vartheta \times {\rm dom}\,\omega \times \Xi\rightarrow\mathbb{R}^{n+m}$ is a $Carath\acute{e}odory$ function \cite{Shapiro2009} if ${G}$ is continuous w.r.t. $(x,y)\in{{\rm dom}\,\vartheta \times {\rm dom}\,\omega }$ for any $\xi$ in ${\Xi}$ and  measurable w.r.t. $\xi\in{\Xi}$ for any $(x,y)\in{{\rm dom}\,\vartheta \times {\rm dom}\,\omega }$.}, which are satisfied conditions {\bf (A1)}-{\bf (A2)}. By the optimality condition \eqref{eq:compact-zop}, a natural idea for finding the minimax point $z^*$ is that at  the $(k+1)$-th iteration, we compute the following subproblem
\begin{equation}\label{eq:compact-z}
{z^{k+1}= {\rm Prox}_{\gamma_k \psi}(z^k-\gamma_k G(z^k,\xi_k))},
 \end{equation}
where $G(z^k,\xi_k)=(G_x(x^k,y^k,\xi_k);-G_y(x^k,y^k,\xi_k))$ .
Specifically, in each iteration, the algorithm alternately  solves proximal subproblems with respect to  $x$ and $y$, which are generated by the first-order oracles ${G}$. A detailed algorithmic framework is presented as follows.
\begin{algorithm}
\caption{Stochastic approximation proximal subgradient (SAPS) method}\label{ASPGalg}
\begin{algorithmic}[1]
\State Given $(x^1,y^1) \in \Re^n\times \Re^m$. Generate i.i.d. samples $\xi_1,\xi_2, \ldots$ of random vector $\xi$. Generate a stochastic subgradient $G(x^1,y^1,\xi^1)=(G_x(x^1,y^1,\xi^1);G_y(x^1,y^1,\xi^1))$. Set $k=1$.
\State Choose a step size $\gamma_k >0$ and compute
\begin{equation}\label{Itformulas}
\begin{array}{ll}
x^{k+1}={\rm argmin}\,\left\{ \vartheta (x)+F(x^k,y^k,\xi_k)+\langle G_x(x^k,y^k,\xi_k),x-x^k \rangle +\displaystyle
\frac{\|x-x^k\|^2}{2\gamma_k}\right\},\\[12pt]
y^{k+1}={\rm argmax}\,\left\{-\omega (y)+F(x^k,y^k,\xi_k)+\langle G_y(x^k,y^k,\xi_k),y-y^k \rangle -\displaystyle
\frac{\|y-y^k\|^2}{2\gamma_k}\right\}.
\end{array}
\end{equation}
\State  Generate a stochastic subgradient $$ G(x^{k+1},y^{k+1},\xi^{k+1})=(G_x(x^{k+1},y^{k+1},\xi^{k+1});G_y(x^{k+1},y^{k+1},\xi^{k+1})).$$
      Update $k+1$ to $k$, and go to {\bf Step 2}.
\end{algorithmic}
\end{algorithm}

 Obviously, we can express $(x^{k+1},y^{k+1})$ as
 \begin{equation}\label{eq:Pr}
 x^{k+1}={\rm Prox}_{\gamma_k\vartheta}(x^k-\gamma_k G_x(x^k,y^k,\xi_k)),\,\,{y^{k+1}={\rm Prox}_{\gamma_k\omega}(y^k+\gamma_k G_y(x^k,y^k,\xi_k))}.
 \end{equation}
If $F(\cdot,\cdot,\xi)$ is continuously differentiable on $\Re^n\times \Re^m$ and $\vartheta(x)=\omega (y)\equiv0$, the first-order oracles ${G}$ become the gradients of $F$ and Algorithm \ref{ASPGalg}  is reduced to stochastic  alternating gradient descent
ascent (SAGDA) algorithm as \cite{Yang2020}. In addition, if functions $\vartheta(x)$ and $\omega (y)$ are the indicator functions shown in \eqref{minimaxXY}. Algorithm \ref{ASPGalg}  is reduced to the stochastic gradient descent (SGD) method  in \cite{Farnia2020,Lei21}.

  The averaging technique, developed by \cite{Polyak1990} and \cite{PJ1992}, considers  an averaged iterated point defined as follows. Let
 $$
 \lambda^k_t=\left(\displaystyle \sum_{j=1}^k \gamma_j\right)^{-1}\gamma_t,\ \, t=1,2,\ldots, k.
 $$
 The averaged iterated point is defined as
 \begin{equation}\label{s2.16}
 \widetilde z^k=\displaystyle \sum_{j=1}^k \lambda^k_j z^j.
 \end{equation}
 Let us analyze convergence properties of the update (\ref{s2.16}), where $z^t,t=1,2,\ldots$ are generated by Algorithm \ref{ASPGalg}. Under the condition {\bf (A3)}, for a saddle point $z^*=(x^*,y^*)$ of $\phi$, from (\ref{eq:SaddleP}), it is reasonable to measure the quality of
an approximate solution $z=(x,y) \in {\cal Z}$ by the error
\begin{equation}\label{errorM}
\epsilon_{\phi}(z)=[\phi (x,y^*)-\phi(x^*,y^*)]+[\phi (x^*,y^*)-\phi(x^*,y)]=\phi (x,y^*)-\phi(x^*,y),
\end{equation}
 which is called the minimax optimality measure at $z$. Obviously, $\epsilon_{\phi}(z)\geq 0$ for any $z \in \mathbb Z$ and $z$ is a saddle point of $\phi$ if and only if $\epsilon_{\phi}(z)=0$.
 \subsection{Expected convergence rate analysis}
Now we consider the convergence properties of the averaged iterated point $\widetilde z^N$ defined as \eqref{s2.16},  where $\{z^k:k=1,\ldots, N\}$  is generated by Algorithm \ref{ASPGalg} for solving a general stochastic convex-concave minimax problem (\ref{minimax}). By assuming first-order oracles with bounded variance, we discuss the convergence rate of the function value of the iterations with expectation.

 The following proposition gives a key property for the sequence generated by Algorithm \ref{ASPGalg}.
 \begin{proposition}\label{impInProp}
 Let $\{z^k\}$ be generated by Algorithm \ref{ASPGalg}, for integer $k>1$, and $\widetilde z^N$ be computed by formula (\ref{s2.16}). Assume that conditions {\bf (A1)} and {\bf (A2)} are satisfied. Denote $\Delta_k=G(z^k,\xi_k)-g(z^k)$. Then for integer $N>0$, for any $z=(x,y)\in \mathbb Z$, the following inequality
 \begin{equation}\label{eq:ineq1}
 \begin{array}{l}
  \phi (\widetilde x^N,y)- \phi (x,\widetilde y^N)\\[1pt]
   \leq
 \displaystyle \left(2\displaystyle\sum_{k=1}^N\gamma_k\right)^{-1}\left[\|z-z^1\|^2-\|z^{N+1}-z\|^2+\displaystyle \sum_{k=1}^N \left[ \gamma_k^2\|v(z^k)+G(z^k,\xi_k)\|^2+2\gamma_k\Delta_k(z-z^k)\right]\right]
 \end{array}
 \end{equation}
 holds for any $v(z^k)=(v_x(x^k),v_y(y^k))$ with $v_x(x^k) \in \partial \vartheta (x^k)$ and $v_y(y^k) \in \partial \omega(y^k)$.

  \end{proposition}
  {\bf Proof}.
  Define
  $
  \widehat x^k=x^k-\gamma_kG_x(z^k,\xi_k),\,\, \widehat y^k=y^k-\gamma_kG_y(z^k,\xi_k).
  $
  Then the two problems for generating $(x^{k+1},y^{k+1})$ in (\ref{Itformulas})  can be expressed as
  \begin{equation}\label{problem1}
  x^{k+1}={\rm argmin}\,\left\{ \vartheta (x) +\displaystyle
\frac{\|x-\widehat x^k\|^2}{2\gamma_k}\right\}\quad \mbox{and}\quad y^{k+1}={\rm argmin}\,\left\{\omega (y)+\displaystyle
\frac{\|y-\widehat y^k\|^2}{2\gamma_k}\right\},
  \end{equation}
   respectively. It follows from Lemma \ref{prox-ineq} (see the Appendix) that  for any $(x,y) \in \mathbb Z$,
   $$
   \begin{array}{l}
   \vartheta (x)+\displaystyle \frac{1}{2\gamma_k}\|x-\widehat x^k\|^2-\displaystyle \frac{1}{2\gamma_k}\|x-x^{k+1}\|^2
   \geq \vartheta (x^{k+1})+\displaystyle \frac{1}{2\gamma_k}\|x^{k+1}-\widehat x^k\|^2,\\[12pt]
   \omega(y)+\displaystyle \frac{1}{2\gamma_k}\|y-\widehat y^k\|^2-\displaystyle \frac{1}{2\gamma_k}\|y-y^{k+1}\|^2
   \geq \omega (y^{k+1})+\displaystyle \frac{1}{2\gamma_k}\|y^{k+1}-\widehat y^k\|^2.
   \end{array}
   $$
   Summing the above two inequalities yields
   \begin{equation}\label{hp1}
\vartheta(x)+\omega(y)\geq  \vartheta(x^{k+1})+\omega(y^{k+1})+\displaystyle \frac{1}{2\gamma_k}\|z^{k+1}-\widehat z^k\|^2-
   \displaystyle \frac{1}{2\gamma_k}\|z-\widehat z^k\|^2+\displaystyle \frac{1}{2\gamma_k}\|z-z^{k+1}\|^2.
   \end{equation}
      Since $\vartheta$ and $\omega$ are convex functions, for any $v_x(x^k) \in \partial \vartheta (x^k)$ and $v_y(y^k) \in \partial \omega(y^k)$, we have
   \begin{equation}\label{hp2}
   \begin{array}{l}
   \vartheta (x^{k+1})\geq \vartheta (x^k)+\langle v_x(x^k), x^{k+1}-x^k \rangle,\\[6pt]
   \omega (y^{k+1})\geq \omega (y^k)+\langle v_y(y^k), y^{k+1}-y^k \rangle.
   \end{array}
   \end{equation}
   Combining (\ref{hp1}) and (\ref{hp2}), we obtain the following inequality
   \begin{equation}\label{hp3}
\vartheta(x)+\omega(y)\geq  \vartheta(x^{k})+\omega(y^{k})+\langle v(z^k), z^{k+1}-z^k \rangle+\displaystyle \frac{1}{2\gamma_k}\left[\|z^{k+1}-\widehat z^k\|^2-
   \|z-\widehat z^k\|^2+\|z-z^{k+1}\|^2\right].
   \end{equation}
   Noting that
   $$
   \begin{array}{rl}
   &\|z^{k+1}-\widehat z^k\|^2-
   \|z-\widehat z^k\|^2\\[4pt]
   &\ \ =\|z^{k+1}-z^k+\gamma_kG(z^k,\xi_k)\|^2-
   \|z-z^k+\gamma_kG(z^k,\xi_k)\|^2\\[4pt]
   &\ \ =\|z^{k+1}-z^k\|^2-\|z-z^k\|^2+2\gamma_k\langle z^{k+1}-z^k, G(z^k,\xi_k)\rangle
   -2\gamma_k \langle z-z^k, G(z^k,\xi_k)\rangle,
   \end{array}
   $$
   we have from (\ref{hp3}) that
   \begin{equation}\label{hp4}
   \begin{array}{rcl}
\vartheta(x)+\omega(y)-  [\vartheta(x^{k})+\omega(y^{k})]
&\geq &\langle v(z^k)+G(z^k,\xi_k), z^{k+1}-z^k \rangle+\displaystyle \frac{1}{2\gamma_k}\|z^{k+1}-z^k\|^2\\[6pt]
&&+   \displaystyle \frac{1}{2\gamma_k}\big[\|z-z^{k+1}\|^2-\|z-z^k\|^2\big]-\langle z-z^k,G(z^k,\xi_k)\rangle\\[10pt]
&\geq& -\displaystyle \frac{\gamma_k}{2}\|v(z^k)+G(z^k,\xi_k)\|^2\\[6pt]
&& +   \displaystyle \frac{1}{2\gamma_k}\big[\|z-z^{k+1}\|^2-\|z-z^k\|^2\big]-\langle z-z^k,G(z^k,\xi_k)\rangle.\\[10pt]
\end{array}
   \end{equation}
Since $g_x(z^k) \in \partial_x f(z^k)$ and $-g_y(z^k)\in \partial_y [-f](z^k)$, we have the following relations
  \begin{equation}\label{hp5}
   \begin{array}{rcl}
   -\langle z-z^k,g(z^k)\rangle &=&-\langle x-x^k,g_x(z^k)\rangle+\langle y-y^k,g_y(z^k)\rangle\\[6pt]
   & \geq& f(z^k)-f(x,y^k)+f(x^k,y)-f(z^k)\\[6pt]
   &=&f(x^k,y)-f(x,y^k).
\end{array}
   \end{equation}
   From (\ref{hp4}) and (\ref{hp5}), we get
    $$
   \begin{array}{l}
[\vartheta(x)+f(x,y^k)-\omega(y^k)]-  [\vartheta(x^{k})+f(x^k,y)-\omega(y)]\\[6pt]
\ \ \geq -\displaystyle \frac{\gamma_k}{2}\|v(z^k)+G(z^k,\xi_k)\|^2+\displaystyle \frac{1}{2\gamma_k}\big[\|z-z^{k+1}\|^2-\|z-z^k\|^2\big]+\langle z-z^k,g(z^k)-G(z^k,\xi_k)\rangle.
\end{array}
   $$
   Therefore, noting $\phi (x,y)=\vartheta(x)+f(x,y)-\omega(y)$ and $\Delta_t=G(x^t,\xi_t)-g(z^t)$ for $t=1,\ldots, N$, we obtain
   \begin{equation}\label{hp6}
   \begin{array}{ll}
&\phi (x^k,y)-\phi (x,y^k)\\[2pt]&\ \ \leq \displaystyle \frac{\gamma_k}{2}\|v(z^k)+G(z^k,\xi_k)\|^2+\displaystyle \frac{1}{2\gamma_k}\big[\|z-z^{k}\|^2-\|z-z^{k+1}\|^2\big]+\langle \Delta_k, z-z^k\rangle.
\end{array}
   \end{equation}
Multiplying $\gamma_k$ on both sides of \eqref{hp6} and summing up over $k=1,\ldots, N$, we obtain
   \begin{equation}\label{hp7}
  \begin{array}{ll}
&\displaystyle \sum_{k=1}^N\gamma_k[\phi (x^k,y)-\phi (x,y^k)]\\[6pt]
&\ \ \leq \displaystyle \sum_{k=1}^N\displaystyle \frac{\gamma_k^2}{2}\|v(z^k)+G(z^k,\xi_k)\|^2+\displaystyle \frac{1}{2}\big[\|z-z^1\|^2-\|z-z^{N+1}\|^2\big] +\displaystyle \sum_{k=1}^N\gamma_k\langle \Delta_k, z-z^k\rangle.
\end{array}
  \end{equation}
   From the convexity of $\phi (\cdot,y)$ and $-\phi(x, \cdot)$ for $(x,y) \in \mathbb Z$, using the definition of $\widetilde z^N$, we obtain
   \begin{equation}\label{hp8}
   \begin{array}{rcl}
   \left[\displaystyle \sum_{k=1}^N \gamma_k\right]\left(\phi (\widetilde x^N, y)-\phi(x, \widetilde y^N)\right)& \leq & \left[\displaystyle \sum_{k=1}^N \gamma_k\right]\displaystyle \sum_{k=1}^N \lambda^N_k \left(\phi (x^k, y)-\phi(x,y^k)\right)\\[12pt]
   &\leq&   \displaystyle \sum_{k=1}^N \gamma_k \left(\phi (x^k, y)-\phi(x,y^k)\right).
   \end{array}
  \end{equation}
  Thus, by combining (\ref{hp7}) and (\ref{hp8}), we obtain  (\ref{eq:ineq1}).\hfill $\Box$\\

  To develop the expected convergence rate, we need modify the condition (\ref{s2.5}) as follows.
 \begin{itemize}
 \item[({\bf A4})]
 There is a positive constant $M_*$ such that for any $z=(x,y)\in  \mathbb Z$, there exists $v(z)=(v_x(x),v_y(y))$, where
 $v(z) \in \partial \vartheta (x)$ and $v_y(y) \in \partial \omega (y)$, such that
 \begin{equation}\label{s2.5a}
 \mathbb E\big[\|v(z)+G(z,\xi)\|^2\big] \leq M^2_*.
 \end{equation}

 \end{itemize}

One sufficient condition for {\bf (A4)}  is that stochastic subgradients $G$ has bounded variance and the subgradients of $\vartheta $ and $\omega$ are bounded. The latter is true for the special cases \eqref{minimaxXY} and \eqref{minimaxregularizedf} and hence they satisfy {\bf (A4)}.

 \begin{theorem}\label{ConvergenceTh}
 Let $\{z^k\}$ be generated by Algorithm \ref{ASPGalg}, for integer $k>1$, and $\widetilde x^k$ be computed by formula (\ref{s2.16}). Assume that conditions {\bf (A1)-(A4)} are satisfied. Then for integer $N>0$, the following properties hold.
 \begin{description}
 \item[(a)] For $\gamma_k=1/\sqrt{N}$, one has
 $
 \mathbb E[\epsilon_{\phi}(\widetilde z^N)] \leq (\|z^1-z^*\|^2+M_*^2)/(2\sqrt{N}).
$
 \item[(b)] For $\gamma_k=\|z^1-z^*\|/(M_*\sqrt{N})$, one has
$
 \mathbb E[\epsilon_{\phi}(\widetilde z^N)] \leq \|z^1-z^*\|M_*/\sqrt{N}.
$
 \item[(c)] For $\gamma_k=\theta\|z^1-z^*\|/(M_*\sqrt{N})$ for $\theta>0$, one has
$
 \mathbb E[\epsilon_{\phi}(\widetilde z^N)] \leq \max\{\theta,\theta^{-1}\}\|z^1-z^*\|M_*/\sqrt{N}.
$
 \end{description}
  \end{theorem}
  {\bf Proof}. It is obvious from (\ref{eq:ineq1}) that
    \begin{equation}\label{eq:ineq2}
   \left(2\displaystyle\sum_{k=1}^N\gamma_k\right)[\phi (\widetilde x^N,y)- \phi (x,\widetilde y^N)]
   \leq
 \displaystyle \|z-z^1\|^2+\displaystyle \sum_{k=1}^N \left[ \gamma_k^2\|v(z^k)+G(z^k,\xi_k)\|^2+2\gamma_k\Delta_k(z-z^k)\right].
  \end{equation}
Notice that $\mathbb E[z^k|\,\xi_{[k-1]}]=z^k$  and $\mathbb E[\Delta_k|\,\xi_{[k-1]}]=0$, where $\xi_{[k-1]}=(\xi_1,\ldots, \xi_{k-1})$.  Taking expectations of both sides of  (\ref{eq:ineq2}),  we obtain from {\bf (A4)} and the definition of $\epsilon_{\phi}$ that
\begin{equation}\label{eq:ineq4}
 \left(2\displaystyle\sum_{k=1}^N\gamma_k\right)  \mathbb E\big[\epsilon_{\phi}(\widetilde z^N)\big]
   \leq
 \displaystyle \|z^1-z^*\|^2+M_*^2\displaystyle \sum_{k=1}^N  \gamma_k^2.
  \end{equation}
The results in \textbf{(a)},  \textbf{(b)} and  \textbf{(c)} are easily obtained by the selection of the step size $\gamma_k$. \hfill $\Box$
   \begin{remark}\label{remark-2.2}
    Let $\mathbb Z={\rm dom}\, \vartheta\times {\rm dom}\,\omega$ be bounded. Then there exists a positive number $D_0$ such that
   $
   \|z-z^1\|\leq D_{\mathbb Z} \, \mbox{ for any } z \in \mathbb Z.
   $
   We can easily obtain
   \begin{description}
 \item[(1)] For $\gamma_t=1/\sqrt{N}$, one has
 $
 \mathbb E[\epsilon_{\phi}(\widetilde z^N)] \leq \displaystyle (D_{\mathbb Z}^2+M_*^2)/(2\sqrt{N}).
 $
 \item[(2)] For $\gamma_t=D_{\mathbb Z}/(M_*\sqrt{N})$, one has
 $
 \mathbb E[\epsilon_{\phi}(\widetilde z^N)] \leq \displaystyle D_{\mathbb Z}M_*/\sqrt{N}.
$
 \item[(3)] For $\gamma_t=\theta D_{\mathbb Z}/(M_*\sqrt{N})$ for $\theta>0$, one has
 $
 \mathbb E[\epsilon_{\phi}(\widetilde z^N)] \leq \max\{\theta,\theta^{-1}\}D_{\mathbb Z}M_*/\sqrt{N}.
$
 \end{description}
   \end{remark}
   \begin{remark}For strongly-convex-strongly-concave minimax problems,  the convergence results of SAPS can be strengthened to the convergence rate of the iterations by \cite{Lan2009}. Select step sizes $\gamma_k=\theta/k$ with the positive constant $\theta>\displaystyle 1/(2\mu)$, where $f$ is $\mu$-strongly convex-concave function on ${\rm dom}\,\vartheta \times {\rm dom}\,\omega$. Assume that conditions {\bf (A1)} and {\bf (A2)} are satisfied,  and  there is a positive constant $M$ such that
$
 \mathbb E\big[\|G(z,\xi)\|^2\big] \leq M^2, \ \forall z \in \mathbb Z,
$ then
one has
 \begin{equation*}\label{s2.9}
 \mathbb E [\|z^k-z^*\|] \leq \displaystyle \max\{\|z^1-z^*\|,\theta M (2\mu \theta-1)^{-1/2}\}/\sqrt{k}.
 \end{equation*}
   \end{remark}

Notice that the convergence results of the stochastic algorithm SAPS is an extension of the nonsmooth convex minimization problem to the convex-concave minimax optimization. In the next subsection, we shall give the related high probability convergence results for general convex-cancave minimax optimization.
  \subsection{ High probability performance analysis}\label{Sec3}
In this subsection, we focus on the high probability convergence of Algorithm \ref{ASPGalg}, which is a stronger conclusion than Theorem \ref{ConvergenceTh}.  To develop the high probability guarantee of the update (\ref{s2.16}) with $\{z^k\}$ being generated by Algorithm \ref{ASPGalg}, we need the following conditions.
   \begin{itemize}
   \item[({\bf B1})] There is a positive constant $M_*$ such that for any $z=(x,y)\in  \mathbb Z$, there exists $v(z)=(v_x(x),v_y(y))$, where
 $v_x(x) \in \partial \vartheta (x)$ and $v_y(y) \in \partial \omega (y)$, such that
 \begin{equation}\label{s2.5a}
 \mathbb E\left[\exp\left\{\displaystyle \frac{\|v(z)+G(z,\xi)\|^2}{M^2_*}\right\}\right] \leq \exp\big\{1\big\}.
 \end{equation}
 \item[({\bf B2})]
There is a positive constant $\kappa_0$ such that for any $z=(x,y)\in  \mathbb Z$ and $\xi \in \Xi$,
 \begin{equation}\label{s2.5aa}
 \|G(z,\xi))\| \leq \kappa_0.
 \end{equation}
 \item[({\bf B3})] Let $\{z^k\}$ be generated by Algorithm \ref{ASPGalg} with the step size $\gamma_k$  defined by
  $
  \gamma_t=\theta\|z^1-z^*\|/(M_*\sqrt{N})
  $
 for $\theta>0$, and for $k\in \textbf{N}, k>1$. Then one has
 \begin{equation}\label{dis-k}
 \|z^k-z^*\|\leq (1+q(\theta)) \|z^1-z^*\|,
 \end{equation}
 where $q(\theta)>0$ is some positive constant depending on $\theta$.
 \end{itemize}
 \begin{remark}\label{remark3.1}
 The condition {\bf (B1)} implies the condition {\bf (A4)}. Indeed, it follows from the Jensen inequality that for the convex function $s \rightarrow \exp \{s\}$,
 $$
  \exp \left[\mathbb E\left\{\displaystyle \frac{\|v(z)+G(z,\xi)\|^2}{M^2_*}\right\}\right]
  \leq \mathbb E\left[\exp\left\{\displaystyle \frac{\|v(z)+G(z,\xi)\|^2}{M^2_*}\right\}\right] \leq \exp\big\{1\big\},
 $$
which induces the condition {\bf (A4)}.
 \end{remark}

 \begin{remark}\label{rem:B3}
 There are two natural cases that the
 assumption {\bf (B3)} holds.  One is the case when ${\rm dom}\, \vartheta \times
 {\rm dom}\,\omega$ is a bounded set. The other is when $(x,y)\rightarrow \phi (x,y^*)-\phi (x^*,y)$ is level-bounded.
 \end{remark}

The follow theorem derives the convergence result of the  minimax optimality measure with probability.
\begin{theorem}\label{high-prob}
 Let {\bf (A1)}--{\bf (A3)}, {\bf (B1)}--{\bf (B3)} be satisfied. Let $\{z^k\}$ be generated by Algorithm \ref{ASPGalg} with the step sizes $\gamma_k$ being defined by
  $
  \gamma_k=\theta\|z^1-z^*\|/(M_*\sqrt{N})
  $
 for $\theta>0$. Then, one has for any $\Omega >4$,
 \begin{equation}\label{high-prob-f}
 {\rm Prob}\, \left\{ \epsilon_{\phi}(\widetilde z^N)\geq \displaystyle \frac{\|z^1-z^*\|M_*\max\{\theta,\theta^{-1}\}}{2\sqrt{N}}\left\{ \left(3+2q(\theta)(1+q(\theta))\right)\Omega\right\} \right\}\leq 2\exp \{-\Omega\}.
 \end{equation}
\end{theorem}
{\bf Proof}. Define $\Gamma_N=\gamma_1+\cdots+\gamma_N$.  It is obvious from (\ref{eq:ineq1}) that
    \begin{equation}\label{eq:ineq2a}
   \Gamma_N\big[\phi (\widetilde x^N,y^*)- \phi (x^*,\widetilde y^N)\big]
   \leq
 \displaystyle \frac{1}{2}\|z-z^1\|^2+\displaystyle \frac{1}{2}\displaystyle \sum_{k=1}^N \gamma_k^2\|v(z^k)+G(z^k,\xi_k)\|^2+\displaystyle \sum_{k=1}^N \gamma_k\Delta_k(z^*-z^k).
  \end{equation}
  Define
  $
  \alpha_N=\displaystyle \frac{1}{2}\displaystyle \sum_{k=1}^N \gamma_k^2\|v(z^k)+G(z^k,\xi_k)\|^2\ \mbox{ and \ }\ \beta_N=\displaystyle \sum_{k=1}^N \gamma_k\Delta_t(z^*-z^k).
  $
  Firstly, we estimate the term $\alpha_N$. From the condition {\bf (B1)}, we have for any $k \in \textbf{N}$,
  $$
   \mathbb E\left[\exp\left\{\displaystyle \frac{\frac{1}{2}\gamma_k^2\|v(z^k)+G(z^k,\xi_k)\|^2}{\frac{1}{2}\gamma_k^2M^2_*}\right\}\right] \leq \exp\big\{1\big\}.
  $$
  Observe that if $r_1,\ldots,r_i$ are nonnegative random variables satisfying $\mathbb E [\exp\{r_k/\sigma_k\}] \leq \exp\{1\}$ for some deterministic $\sigma_k>0$. It follows from the convexity of the function $s\rightarrow \exp (s)$ that
  \begin{equation}\label{eq6.13}
  \mathbb E \left[ \displaystyle \frac{\sum_{k=1}^i r_k}{\sum_{k=1}^i \sigma_k}\right]\leq
   \mathbb E \left[ \displaystyle \sum_{k=1}^i\displaystyle \frac{\sigma_k}{\sum_{l=1}^i \sigma_l}\exp \left\{\displaystyle
   \frac{r_k}{\sigma_k}\right\}\right]\leq \exp \big\{1\big\}.
  \end{equation}
  Substituting $r_k=1/2\gamma_k^2\|v(z^k)+G(z^k,\xi_k)\|^2$ and $\sigma_k=1/2M_*^2\gamma_k^2$ into (\ref{eq6.13}), we obtain for $\gamma_{\alpha_N}=1/2M_*^2 \sum_{k=1}^N\gamma_k^2$ that
  \begin{equation}\label{eq6.11}
  \mathbb E \left[ \exp \big\{ \alpha_N/\gamma_{\alpha_N} \big\}\right]\leq  \exp \big\{1\big\}.
  \end{equation}
  Thus, by the Markov inequality,  we have for any $\Omega>0$,
  $
  {\rm Prob}\,\big\{\alpha_N\geq (1+\Omega) \gamma_{\alpha_N}\big\} \leq \exp\{-\Omega\}.
  $
  Namely, we obtain the following estimate
  \begin{equation}\label{eq6.12}
  {\rm Prob}\,\left\{\alpha_N\geq \frac{1}{2}(1+\Omega) \theta^2\|z^1-z^*\|^2\right\} \leq \exp\{-\Omega\}.
  \end{equation}

  Secondly, we calculate the boundedness of the term $\beta_N$. Let $\varsigma_k=\gamma_k \langle \Delta_k, z^*-z^k\rangle$ for $k \in \textbf{N}$.  Observing that $z^*-z^k$ is a deterministic function
of $\xi_{[k-1]}$ while $\mathbb E\left[ \Delta_k\,|\,\xi_{[k-1]}\right]=0$, we know that the sequence $\{\varsigma_k\}_{k=1}^N$ of random real
variables forms a martingale difference. It follows from the condition {\bf (B2)} that $\|\Delta_k\|=\|G(z^k,\xi_k)-g(z^k)\|\leq 2\kappa_0$. We have from the
 assumption {\bf (B3)}  that $|\varsigma_k|\leq (1+q(\theta))\|z^1-z^*\|\gamma_k \|\Delta_k\|\leq  2\kappa_0 (1+q(\theta))\|z^1-z^*\|\gamma_k$. Therefore
 $$
   \mathbb E\left[\exp\left\{\displaystyle \frac{\varsigma_k^2}{4\gamma_k^2\|z^1-z^*\|^2(1+q(\theta))^2\kappa_0^2}\right\}\,\Big|\, \xi_{[k-1]}\right] \leq \exp\big\{1\big\}.
  $$

 Define $\eta_k=2\kappa_0(1+q(\theta))\|z^1-z^*\|\gamma_k$. Since $\varsigma_k$ is a deterministic function of $\xi_{[k]}$ with $\mathbb E[\varsigma_k\,|\,\xi_{[k-1]}]=0$ and $\mathbb E\left [\exp\{\varsigma_k^2/\eta_k^2\}\,|\, \xi_{[k-1]}\right] \leq \exp\{1\}$,  we have that for any $\tau >0$ and $0 < \tau \eta_k <1$,
 $$
 \begin{array}{rcl}
 \mathbb E \left[ \exp \{\tau \varsigma_k\}\,\Big|\, \xi_{[k-1]}\right] \leq \mathbb E \left[ \exp \{\tau^2 \varsigma_k^2\}\, |\, \xi_{[k-1]}\right] =\mathbb E \left[ \left(\exp \{\varsigma_k^2/\eta_k^2\}\right)^{\tau^2\eta_k^2}\, \Big|\, \xi_{[k-1]}\right]\leq \exp\{\tau^2\eta_k^2\}.
 \end{array}
 $$
where the first inequality is  derived from $ \exp (x) \leq x +\exp (x^2)$.  When $\tau \eta_k \geq 1$, we have
  $$
  \begin{array}{rcl}
 \mathbb E \left[ \exp \{\tau \varsigma_k\}\,\Big|\, \xi_{[k-1]}\right] \leq\mathbb E \left[ \exp \{\frac{1}{2}\tau^2\eta_k^2+\frac{1}{2} \varsigma_k^2/\eta_k^2\}\, |\, \xi_{[k-1]}\right]
 \leq \exp \left\{\frac{1}{2}\tau^2\eta_k^2+\frac{1}{2}\right\}\leq \exp\{\tau^2\eta_k^2\}.
 \end{array}
  $$
  Thus, in both cases, we have
  $
  \mathbb E \left[ \exp \{\tau \varsigma_t\}\,\Big|\, \xi_{[t-1]}\right] \leq  \exp\{\tau^2\eta_t^2\}.
  $
  Therefore, we have
  $$
  \begin{array}{rcl}
  \mathbb E \left[\exp\left\{\tau \beta_i\right\}\right] = \mathbb E \left[\exp\left\{\tau \beta_{i-1} \right\}
  \mathbb E \left[\exp\left\{\tau \varsigma_i\right\}\,|\, \xi_{[i-1]}\right]\right]
  \leq \exp\left\{\tau^2\eta_i^2\right\}\mathbb E \left[\exp\left\{\tau \beta_{i-1} \right\}\right],
    \end{array}
  $$
 which implies
 $$
  \mathbb E \left[\exp\left\{\tau \beta_N\right\}\right]\leq \exp \left\{\tau^2\displaystyle \sum_{t=1}^N \eta_t^2  \right\}.
 $$
 By the Markov inequality for $\Theta>0$, it holds
 $$
 {\rm Prob}\, \left\{\beta_N > \Theta \sqrt{\displaystyle \sum_{t=1}^N \eta_t^2}  \right\}
 \leq \exp \left\{\tau^2 \displaystyle \sum_{t=1}^N \eta_t^2 \right\}\exp \left\{-\tau\Theta \sqrt{\displaystyle \sum_{t=1}^N \eta_t^2} \right\}.
 $$
 When choosing $\tau=\frac{1}{2}\Theta\left(\displaystyle \sum_{t=1}^N \eta_t^2\right)^{-1/2}$, we get the following estimate
 \begin{equation}\label{eq6.14}
 \begin{array}{rl}
  &{\rm Prob}\, \left\{\beta_N > 2\kappa_0(1+q(\theta))\theta M_*^{-1}\|z^1-z^*\|^2\Theta   \right\}\\[6pt]
  &\ \ ={\rm Prob}\, \left\{\beta_N > q(\theta)(1+q(\theta))\|z^1-z^*\|^2\Theta   \right\}
 \leq \exp \left\{-\displaystyle \frac{\Theta^2}{4}\right\}.
 \end{array}
  \end{equation}

Finally, combining (\ref{eq:ineq2a}), (\ref{eq6.12}) and (\ref{eq6.14}), we get the following for any positive $\Omega$ and $\Theta$,
  \begin{equation}\label{eq6.14a}
  \begin{array}{rl}
  &{\rm Prob}\,\left\{ \Gamma_N\epsilon_{\phi}(\widetilde z^N)
   >
 \displaystyle \frac{\|z-z^1\|^2}{2}\left\{1+(1+\Omega)\theta^2+2q(\theta)(1+q(\theta))\Theta\right\}\right\}
 \\[12pt]
 &\ \  \leq\exp \{-\Omega\}+\exp \left\{-\displaystyle \frac{\Theta^2}{4}\right\}.
 \end{array}
 \end{equation}
Let $\Theta=2\sqrt{\Omega}$. Noting $\Omega>4$, we have that $\Theta < \Omega$. Hence, (\ref{eq6.14a}) implies
\begin{equation}\label{eq6.14b}
  \begin{array}{l}
  {\rm Prob}\,\left\{ \Gamma_N\epsilon_{\phi}(\widetilde z^N)
   >
 \displaystyle \frac{\|z^1-z^*\|^2}{2}\left\{1+\theta^2+[\theta^2+2q(\theta)(1+q(\theta))]\Omega\right\}\right\}
 \leq 2\exp \{-\Omega\}.
 \end{array}
 \end{equation}
 Noting $\Gamma_N=\theta M_*^{-1}\|z^1-z^*\|\sqrt{N}$, we have $$
 \begin{array}{rl}
  \displaystyle \frac{\|z^1-z^*\|}{\Gamma_N}\left\{1+\theta^2+[\theta^2+2q(\theta)(1+q(\theta))]\Omega\right\}
  < \displaystyle \displaystyle \frac{M_*\max\{\theta,\theta^{-1}\}}{\sqrt{N}}\left\{
  \left( 3+2q(\theta)(1+q(\theta))\right)\Omega\right\}.
 \end{array}
 $$

 Therefore, we obtain (\ref{high-prob}) from (\ref{eq6.14b}). \hfill $\Box$\\

 Taking $\Omega=\log N$ in Theorem \ref{high-prob}, we obtain the following result.
 \begin{corollary}\label{cor-hp}
 Let the conditions of Theorem \ref{high-prob} be satisfied. Then for $N> \exp\{4\}$,
 \begin{equation}\label{high-prob-cor}
 {\rm Prob}\, \left\{ \epsilon_{\phi}(\widetilde z^N)\geq \displaystyle \|z^1-z^*\|M_*\bar{q}(\theta) N^{-1/2}\log N\right\}\leq \displaystyle \frac{2}{N},
 \end{equation}
 where $\bar{q}(\theta)=\max\{\theta,\theta^{-1}\}(3+2q(\theta)(1+q(\theta)))$.
 \end{corollary}

By choosing $\Omega$ appropriately, Corollary \ref{cor-hp} shows that the minimax optimality measure $\epsilon_{\phi}(\widetilde z^N)$ converges to zero with probability 1, as $N$ tends to infinity.
 \begin{remark}\label{remark-3.2}
 Noting that the condition (\ref{s2.5aa}) implies
 $
 \mathbb E \|G(z,\xi)\|^2 \leq \kappa_0^2
 $
 and the condition (\ref{s2.5a}) implies
 $$
 \mathbb E \|v(z)\|^2+\mathbb E \|G(z,\xi)\|^2=\mathbb E \|(v(z)+G(z,\xi))\|^2\leq M^2_*,
$$
we think that the condition $M_*\geq \kappa_0$ in \eqref{high-prob-cor} is reasonable under {\bf (B1)} and {\bf (B2)}.
 \end{remark}
 \section{Linearized Stochastic Approximation Augmented Lagrange  Method }\label{Sec:4}
 \setcounter{equation}{0}
 In this section, we consider the special stochastic convex-concave minimax problem coming from the stochastic convex conic optimization problem \eqref{CP}. However, this is a typical problem which does not satisfy conditions {\bf (A4)} or {\bf (B1)}. This motivates us to consider an improved version of Algorithm \ref{ASPGalg} for solving the minimax problem in terms of  the augmented Lagrangian function.

We assume that $f$ is a convex function and the set-valued mapping $x \rightarrow g(x)-{\cal K}$ is graph convex in the problem (\ref{CP}). Under these conditions,  (\ref{CP}) becomes a convex programming problem, see Definition 2.163 of \cite{BS2000}.  Certainly, if we assume that   for every
$\xi \in \Xi$, $F(\cdot,\xi)$ is convex  and the set-valued mapping $x \rightarrow G(x,\xi)-{\cal K}$ is graph-convex over ${\cal O}$, then $f$ is a convex function and $x \rightarrow g(x)-{\cal K}$ is graph-convex over ${\cal O}$.  Moreover, we assume that for every
$\xi \in \Xi$, the function $F(\cdot,\xi)$  and the mapping $G(\cdot, \xi)$ are smooth over  ${\cal O}$ so that $f(x)=\mathbb E(F(x,\xi))$ and $g(x)=\mathbb E (G(x,\xi))$ are smooth over ${\cal O}$, too. Under the above conditions, the problem (\ref{CP}) is a smooth convex conic optimization problem.

It can be easily check that the mapping $x \rightarrow g(x)-{\cal K}$ is graph-convex over ${\cal O}$ if and only if
\begin{equation}\label{eq:c1}
 g((1-\lambda)x^1+\lambda x^2) -[(1-\lambda) g(x^1)+\lambda g(x^2)] \in {\cal K},\,\, \forall x^1,x^2 \in {\cal O}, \forall \lambda \in [0,1].
\end{equation}
The following lemma lists some important properties of the graph-convexity of the set-valued mapping $x \rightarrow g(x)-{\cal K}$.
\begin{lemma}\label{lemm-4.0}
Let the mapping $x \rightarrow g(x)-{\cal K}$ be graph-convex over an open convex set ${\cal O}$ and $g$ be differentiable at $x \in {\cal O}$. Then
\begin{equation}\label{eq:4.0}
g(x)+ {\rm D} g(x)(z-x) -g(z)\in {\cal K}, \quad \forall\ x, z \in {\cal O}
\end{equation}
and
\begin{equation}\label{eq:4.01}
\|\Pi_{{\cal K}^{\circ}}[g(x)+ {\rm D} g(x)(z-x)]\| \leq \|\Pi_{{\cal K}^{\circ}}[g(z)]\|.
\end{equation}
\end{lemma}
{\bf Proof}. For  $x, z \in {\cal O}$ and $t \in (0,1)$, we have $x+t(z-x) \in {\cal O}$. From (\ref{eq:c1}), we have
$$
g(x+t(z-x))-(1-t)g(x)-tg(z)=g((1-t)x+tz)-(1-t)g(x)-tg(z)\in {\cal K},
$$
which implies
$
t^{-1}[g(x+t(z-x))-g(x)]+g(x)-g(z)\in {\cal K}, \forall t \in (0,1).
$
Taking the limit when $t \searrow 0$ yields
$
g(x)+{\rm D}g(x)(z-x)-g(z) \in {\cal K}.
$

Now we prove (\ref{eq:4.01}). Let $\eta_1=\|\Pi_{{\cal K}^{\circ}}(g(z))\|^2$ and $\eta_2=\|\Pi_{{\cal K}^{\circ}}(g(x)+{\rm D}g(x)(z-x))\|^2$. Then
\begin{equation}\label{eq:aa2}
\eta_1=\min\{ \|u'\|^2: u'+g(z) \in {\cal K}\}\ \mbox{ and }\ \eta_2=\min\{\|u'\|^2: u'+g(x)+{\rm D}g(x)(z-x) \in {\cal K}\}.
\end{equation}
In view of the graph-convexity of the set-valued mapping  $x \rightarrow g(x)-{\cal K}$, we have for $x\in X$,
$
g(x)+{\rm D}g(x)(z-x)-g(z) \in {\cal K}.
$
Thus, if $u'+g(z)\in {\cal K}$, then
$$
u'+g(x)+{\rm D}g(x)(z-x) =u'+g(z)+[g(x)+{\rm D}g(x)(z-x) -g(z)] \in {\cal K},
$$
which implies that the feasible region for defining $\eta_1$ is a subset of that for defining $\eta_2$. Therefore we have $\eta_1 \geq \eta_1$. The proof is completed.\hfill $\Box$\\

The general assumptions by using the stochastic approximation (SA) technique  include the independent identically distributed (i.i.d.) samples $\xi_1,\xi_2,\ldots,$ and the unbiased gradients $\nabla_x F(x,\xi)$ and ${\rm D}_xG(x,\xi)$.  For the stochastic convex problem \eqref{CP}, assume that the condition \textbf{(A1)} is satisfied and replace \textbf{(A2)} with the following assumption.
 \begin{itemize}
 \item[{\bf (C1)}] For any point $(x,\xi)\in {\cal O}\times \Xi$, stochastic gradients $\nabla_x F(x,\xi)$ and ${\rm D}_xG(x,\xi)$ satisfy $\nabla f(x)=\mathbb{E}[\nabla_x F(x,\xi)]$
       and ${\rm D}g(x)=\mathbb{E}[{\rm D}_x G(x,\xi)]$.
\end{itemize}

Furthermore, we assume the following assumptions about the problem (\ref{CP}).
\begin{itemize}
\item[{\bf (D1)}] Let  $R>0$  be a positive parameter such that for any $x',x'' \in X$,
        $
        \|x'-x''\|\leq R.
        $
\item[{\bf (D2)}]  There exists a constant $\nu_g>0$  such that for any $\xi_k$ and  $x \in {\cal O}$,
        $
        \|G(x,\xi_k)\| \leq \nu_g.
        $
\item[{\bf (D3)}] There exist constants $\kappa_f>0$ and $\kappa_g>0$ such that for each $\xi_k$ and $x\in {\cal O}$,
                $
        \|\nabla_x F(x,\xi_k)\| \leq \kappa_f, \,\, \|{\rm D}_xG(x,\xi_k)\| \leq \kappa_g.
        $
\item[{\bf (D4)}]
     Assume that the Slater condition holds; namely, there exists
     $
      \hat x \in {\rm int\,X}  \mbox{ such that } g(\hat x) \in {\rm int}\, {\cal K}.
     $
     \end{itemize}

                If the optimal solution set of the problem (\ref{CP}) is nonempty, it follows from  Theorem 2.165 of \cite{BS2000} that, under the above Slater condition \textbf{(D4)}, the conjugate dual \eqref{CD} of the problem (\ref{CP}) has a compact optimal solution set and the optimal value of the problem (\ref{CP}) is equal to the optimal value of the problem (\ref{CD}).  Since ${\cal K}$ is a closed convex cone, we have $\delta^*_{{\cal K}}=\delta_{{\cal K}^{\circ}}$. In this case, the problem (\ref{CP}) can be expressed as the following convex-concave minmax problem
        \begin{equation}\label{minimax-cp}
        \min_x \max_y \phi (x,y)=\delta_X(x)+l(x,y)   -\delta_{{\cal K}^{\circ}}(y)
        \end{equation}
        with
        \begin{equation}\label{eq:Lxi}
        l(x,y)=\mathbb E \big[L(x,y,\xi)\big] \,\, \mbox{ and }\,\,
        L(x,y,\xi)= F(x,\xi)+\langle y, G(x,\xi)\rangle.
        \end{equation}

       Let $x^*$ and $y^*$ be optimal solutions to the problem (\ref{CP}) and the problem (\ref{CD}), respectively. Then
       $$
       \phi (x^*,y) \leq \phi(x^*,y^*)\leq \phi(x,y^*), \mbox{ for } (x,y) \in \Re^n  \times {\cal Y}.
       $$
Obviously we have the following expression
$$
\nabla_x L(x,y,\xi)=\nabla_x F(x,\xi)+{\rm D} G(x,\xi)^*y, \  \  {\rm D}_y L(x,y,\xi)=G(x,\xi).
$$
Noting that ${\cal K}^{\circ}$ is an unbounded set, it is not possible to guarantee the boundedness of $\nabla_x L(x,y,\xi)$, which implies that neither {\bf (A4)} nor {\bf  (B1)} is satisfied. This observation motives us to consider a variant of Algorithm \ref{ASPGalg}, which is able to handle this difficulty. For this purpose, instead of the ordinary Lagrange dual (\ref{CD}),
we use the following augmented Lagrange dual of the problem (\ref{CP}):
\begin{equation}\label{ACD}
         \max_{y\in {\cal K}^{\circ}}\min_{x\in X} l_{\sigma}(x,y), \,\,\mbox{ where }\,\, l_{\sigma}(x,y)
         =\mathbb E[F(x,\xi)]+\displaystyle \frac{1}{2\sigma}\left[ \|\Pi_{{\cal K}^{\circ}}(y+\sigma \mathbb E[ G(x,\xi)]) \|^2-\|y\|^2\right].
        \end{equation}

Due to the convexity of the problem (\ref{CP}),  the optimal value of \eqref{ACD}  is equal to that of the following minimax problem
\begin{equation}\label{minmaxAL}
         \min_{x\in X}\max_{y\in{\cal K}^{\circ}} l_{\sigma}(x,y).
        \end{equation}
        Now we focus on solving the stochastic convex-concave minimax problem (\ref{minmaxAL}). Define the linearized approximations of $F(\cdot,\xi_k)$ at $x^k$ and $G(\cdot,\xi_k)$ at $x^k$ as
        \begin{equation}\label{eq:lflg}
        l_f^k(x):=F(x^k,\xi_k)+\langle \nabla_xF(x^k,\xi_k),x-x^k \rangle,\quad l_g^k(x):=G(x^k,\xi_k)+{\rm D} G(x^k,\xi_k)(x-x^k)),
        \end{equation}
        respectively. The corresponding augmented Lagrangian function is expressed  as
      \begin{equation}\label{eq:lsigma}
        l^k_{\sigma}(x,y)=l_f^k(x)+
\displaystyle \frac{1}{2\sigma}\left[ \left\|\Pi_{{\cal K}^{\circ}}(y+\sigma l_g^k(x)) \right\|^2-\|y\|^2\right].
        \end{equation}
        We propose the following linearized stochastic approximation method in terms of the augmented Lagrangian function for solving the minimax problem (\ref{minmaxAL}).
        \begin{algorithm}
\caption{Linearized stochastic approximation augmented Lagrange (LSAAL)  method}\label{ASalm}
\begin{algorithmic}[1]
\State
 Choose an initial point $x^1 \in X$, $y^1=0 \in {\cal Y}$ and select parameters
$\sigma >0$. Generate i.i.d. sample $\xi_1,\xi_2, \ldots$ of random vector $\xi$.  Compute the stochastic gradient $\nabla_x F(x^1,\xi^1)$ and the stochastic derivative ${\rm D}_x G(x^1,\xi^1)$. Set $k=1$.
\State Compute
\begin{equation}\label{Itformula}
\begin{array}{ll}
x^{k+1}&=\displaystyle{\rm argmin}_{x \in X}\,\left\{l^k_{\sigma}(x,y^k)  +\displaystyle
\frac{1}{2\sigma}\|x-x^k\|^2\right\},\\[12pt]
y^{k+1}&=\displaystyle{\rm argmax}_{y}\left\{l^k_{\sigma}(x^{k+1},y^k)+{\rm D}_yl^k_{\sigma}(x^{k+1},y^k)(y-y^k)-  \displaystyle
\frac{1}{2\sigma}\|y-y^k\|^2\right\}\\[12pt]
&=\Pi_{{\cal K}^{\circ}}(y^k+\sigma l_g^k(x^{k+1}))=\Pi_{{\cal K}^{\circ}}(y^k+\sigma (G(x^k,\xi_k)+{\rm D} G(x^k,\xi_k)(x^{k+1}-x^k))).
\end{array}
\end{equation}
\State Compute the stochastic gradient $\nabla_x F(x^{k+1},\xi^{k+1})$ and the stochastic derivative ${\rm D}_x G(x^{k+1},\xi^{k+1})$.
      Update $k+1$ to $k$, and go to {\bf Step 2}.
\end{algorithmic}
\end{algorithm}

The following auxiliary lemma will be used for several times in the sequel.
\begin{lemma}\label{lem:opt-x}
 For any $x\in X$, we have
\begin{equation}\label{eq:opt-x-1}
\begin{array}{ll}
\displaystyle\langle \nabla_x F(x^k,\xi_k),x^{k+1}-x^k \rangle + \frac{1}{2\sigma}\|y^{k+1}\|^2
+ \frac{1}{2\sigma} \|x^{k+1}-x^k\|^2  \\[10pt]
\ \ \leq \displaystyle\langle \nabla_x F(x^k,\xi_k),x-x^k \rangle + \frac{1}{2\sigma}\left[ \|\Pi_{{\cal K}^{\circ}}\left(y^k+\sigma
(G(x^k,\xi_k)+{\rm D}_x(x^k,\xi_k)(x-x^k))\right)\|^2\right]\\[15pt]
\ \ \quad+ \displaystyle\frac{1}{2\sigma}(\|x-x^k\|^2-\|x-x^{k+1}\|^2).
\end{array}
\end{equation}
In particular, if we take $x=x^k$, it yields
\begin{equation}\label{eq:opt-x-2}
\begin{array}{ll}
\displaystyle\langle \nabla_x F(x^k,\xi_k),x^{k+1}-x^k \rangle + \frac{1}{2\sigma}\|y^{k+1}\|^2
+ \frac{1}{\sigma} \|x^{k+1}-x^k\|^2
\leq \displaystyle \frac{1}{2\sigma}\left[ \|\Pi_{{\cal K}^{\circ}}\left(y^k+\sigma
G(x^k,\xi_k)\right)\|^2\right].
\end{array}
\end{equation}
\end{lemma}
{\bf Proof}.
From the definition of $x^{k+1}$ in  (\ref{Itformula}) and its optimality conditions, we have that $x^{k+1}$ is also the optimal solution to the following problem
\[
\begin{array}{ll}
\displaystyle\min_{x}\ \displaystyle\langle \nabla_x F(x^k,\xi_k),x-x^k \rangle + \frac{1}{2\sigma}\left[ \|\Pi_{{\cal K}^{\circ}}\left(y^k+\sigma
(G(x^k,\xi_k)+{\rm D}_x(x^k,\xi_k)(x-x^k))\right)\|^2\right]\\[10pt]
\quad\quad+ \displaystyle\frac{1}{2\sigma}(\|x-x^k\|^2-\|x-x^{k+1}\|^2).
\end{array}
\]
Then,
the claimed results are obvious.
\hfill $\Box$
\begin{corollary}\label{coro-xd}
Let $\{(x^k,y^k)\}$ be generated by Algorithm \ref{ASalm}. Then for $k=1,2,\ldots$,
\begin{equation}\label{eq:xdd}
\|x^{k+1}-x^k\|\leq \sigma \Big\{[\kappa_f+\nu_g\kappa_g\sigma]+\kappa_g\|y^k\|\Big\}.
\end{equation}
\end{corollary}
{\bf Proof}. From (\ref{eq:opt-x-2}), noting the convexity and differentiability of mapping $y \rightarrow \|\Pi_{{\cal K}^{\circ}}(y)\|^2$,  we have that
$$
\begin{array}{l}
-\kappa_f\|x^{k+1}-x^k \|
+ \frac{1}{\sigma} \|x^{k+1}-x^k\|^2 \\[10pt]
\ \ \leq -\|\nabla_x F(x^k,\xi_k)\|\|x^{k+1}-x^k \|
+ \frac{1}{\sigma} \|x^{k+1}-x^k\|^2 \\[10pt]
\ \ \leq \displaystyle \frac{1}{2\sigma}\left[ \|\Pi_{{\cal K}^{\circ}}\left(y^k+\sigma
G(x^k,\xi_k)\right)\|^2-\|\Pi_{{\cal K}^{\circ}}\left(y^k+\sigma
[G(x^k,\xi_k)+{\rm D}G(x^k,\xi_k)(x^{k+1}-x^k)]\right)\|^2\right]\\[10pt]
\ \ \leq \langle \Pi_{{\cal K}^{\circ}}\left(y^k+\sigma
G(x^k,\xi_k)\right), {\rm D}G(x^k,\xi_k)(x^k-x^{k+1})\rangle\\[10pt]
\ \ \leq \|\Pi_{{\cal K}^{\circ}}\left(y^k+\sigma
G(x^k,\xi_k)\right)\|\| {\rm D}G(x^k,\xi_k)(x^k-x^{k+1})\|\\[10pt]
\ \ \leq [\|y^k\|+\sigma \|G(x^k,\xi_k)\|]\| {\rm D}G(x^k,\xi_k)\|\|x^{k+1}-x^k\|\\[10pt]
\ \ \leq  [\|y^k\|+\sigma \nu_g]\kappa_g\|x^{k+1}-x^k\|,
\end{array}
$$
which implies (\ref{eq:xdd}). The proof is completed. \hfill $\Box$\\

Under the Slater condition \textbf{(D4)}, we can prove  the following  conditional expected estimate related to the multipliers.
\begin{lemma}\label{lem:aux4}
Let assumptions {\bf (A1), (C1), (D4)}  be satisfied and $\{(x^k,y^k)\}$ be generated by Algorithm \ref{ASalm}. Then, there exists a positive number $\varepsilon_0 >0$ such that for any positive integers $\tilde{k} \leq k-1$,
\[
 \mathbb{E}\left[\langle y^{k}, G(\widehat x, \xi_{k}) \rangle \,|\, {\cal F}_{\tilde{k}}\right]
\leq -\varepsilon_0 \mathbb{E}\left[\|y^{k}\| \,|\, {\cal F}_{\tilde{k}}\right].
\]
\end{lemma}
{\bf Proof}.
It follows from {\bf (C4)} that there exists an positive number $\varepsilon_0$ such that
$
g(\hat{x} ) +\varepsilon_0\textbf{B}\subset {\cal K}.
$
This implies that for any nonzero $y \in {\cal K}^{\circ}$,
$
\delta^*_{ g(\hat{x})+\varepsilon_0 \textbf{B}}(y)\leq 0,
$
which yields
$$
\langle y, g(\hat{x})\rangle\leq -\delta^*_{\varepsilon_0 \textbf{B}}(y)=-\varepsilon_0\|y\|.
$$
Noticing that $y^{k} \in {\cal F}_{k}$  and ${\cal F}_{\tilde{k}}\subseteq {\cal F}_{k1}$ for $\tilde{k} \leq k-1$, we have
from $y^{k} \in {\cal K}^{\circ}$ that
\[
\begin{array}{rcl}
\mathbb{E}\left[\langle y^{k}, G(\hat{x}, \xi_{k})\rangle \,|\, {\cal F}_{\tilde{k}}\right] =\mathbb{E} \left[\mathbb{E}\left[\langle y^{k}, G(\hat{x}, \xi_{k})\rangle \,|\, {\cal F}_{k}\right]\,|\,{\cal F}_{\tilde{k}}\right]=\mathbb{E}\left[\langle y^{k},g(\hat{x})\rangle\,|\,{\cal F}_{\tilde{k}}\right]
\leq-\varepsilon_0 \mathbb{E}\left[\|y^{t_1}\| \,|\, {\cal F}_{t_2}\right].
\end{array}
\]
The proof is completed.
\hfill $\Box$\\

The next result shows a self-adjusting property  of $\|y^k\|$, which is essential to establish the expected convergence rate of the minmax optimality measure.
\begin{lemma}\label{lem:aux5}
Let assumptions {\bf (A1), (C1), (D1)-(D4)} be satisfied. Let $\{(x^k,y^k)\}$ be generated by Algorithm \ref{ASalm} and $s > 0$ be an arbitrary integer.  Define
 \begin{equation}\label{eq:theta9}
 \vartheta (\sigma,s):= \frac{\varepsilon_0\sigma s}{2}+\sigma \beta_0(s-1)+ \frac{ R^2}{\varepsilon_0\sigma s}+ \frac{(\kappa_f+2 \nu_g^2+2\kappa_g^2R^2)\sigma}{\varepsilon_0},
 \end{equation}
where $\varepsilon_0 >0$ is defined in Lemma \ref{lem:aux4}. Then, we have
\begin{equation}\label{eq:6}
|\|y^{k+1}\|-\|y^k\||\leq \sigma \beta_0
\end{equation}
and
\begin{equation}\label{eq:7}
\mathbb{E}\left [ \|y^{k+s}\|-\|y^k\| \,|\, {\cal F}_k\right]
\leq \left
\{
\begin{array}{ll}
s \sigma \beta_0, & \mbox{if } \|y^k\| < \vartheta (\sigma, s);\\[6pt]
-s \displaystyle \frac{\sigma \varepsilon_0}{2}, & \mbox{if } \|y^k\| \geq  \vartheta (\sigma,s),
\end{array}
\right.
\end{equation}
where $\beta_0=\nu_g+\kappa_gR$.
\end{lemma}
{\bf Proof}.
From conditions {\bf (D1)-(D3)}, it follows from   (\ref{Itformula}) and the nonexpansion property of the projection $\Pi_{{\cal K}^{\circ}}(\cdot)$ that
\begin{equation}\label{eq:c1}
\begin{array}{rcl}
|\|y^{k+1}\|-\|y^k\||
&\leq&\|y^{k+1}-y^k\|  =\|\Pi_{{\cal K}^{\circ}}[y^k+\sigma (G(x^k,\xi_k)+{\rm D}_xG(x^k,\xi_k)(x^{k+1}-x^k))] -\Pi_{{\cal K}^{\circ}}[y^k]\|\\[6pt]
& \leq&
\sigma \|G(x^k,\xi_k)+{\rm D}_xG(x^k,\xi_k)(x^{k+1}-x^k)\|\\[6pt]
& \leq& \sigma [\nu_g+\kappa_g R],
\end{array}
\end{equation}
which implies (\ref{eq:6}).  This also gives that $\|y^{k+s}\|-\|y^k\|\leq s \sigma \beta_0$. Hence, we only need to demonstrate the second part  in (\ref{eq:7}) under the case $\|y^k\| \geq  \vartheta (\sigma,s)$.

For a given positive integer $s$, suppose that $\|y^k\| \geq \vartheta (\sigma,s)$. For any $t \in \{k,k+1,\ldots,k+s-1\}$, from (\ref{eq:opt-x-1}) and  Lemma \ref{lemm-4.0} about the convexity of $G(\cdot,\xi_t)$ with respect to ${\cal K}$, one has
$$
\begin{array}{rcl}
 \|y^{t+1}\|^2&\leq & \left[\| \Pi_{{\cal K}^{\circ}}(y^t+\sigma (G(x^t,\xi_t)+{\rm D} G(x^t,\xi_t)(\hat{x}-x^t)))\|^2\right]\\[6pt]
&\leq &\left[\| (y^t+\sigma (G(x^t,\xi_t)+{\rm D} G(x^t,\xi_t)(\hat{x}-x^t))\|^2\right]\\[6pt]
&\leq& \|y^t\|^2+2\sigma
 \langle y^t,(G(x^t,\xi_t)+{\rm D} G(x^t,\xi_t)(\hat{x}-x^t))\rangle+\sigma^2\|(G(x^t,\xi_t)+{\rm D} G(x^t,\xi_t)(\hat{x}-x^t))\|^2 \\[6pt]
&\leq& \|y^t\|^2+2\sigma
 \langle y^t,G(\hat{x},\xi_t)\rangle+\sigma^2\|(G(x^t,\xi_t)+{\rm D} G(x^t,\xi_t)(\hat{x}-x^t))\|^2,
\end{array}
$$
which implies
$$
\begin{array}{l}
\langle \nabla_x F(x^t,\xi_t),x^{t+1}-x^t\rangle + \frac{1}{2\sigma}\|y^{t+1}\|^2+ \frac{1}{2\sigma}\|x^{t+1}-x^t\|^2\\[10pt]
\ \ \leq \langle \nabla_x F(x^t,\xi_t),\hat{x}-x^t\rangle + \frac{1}{2\sigma}\|y^t\|^2+
 \langle y^t,G(\hat{x},\xi_t)\rangle\\[10pt]
\ \ \quad+\frac{\sigma}{2}\|(G(x^t,\xi_t)+{\rm D} G(x^t,\xi_t)(\hat{x}-x^t))\|^2 + \frac{1}{2\sigma}(\|\hat{x}-x^t\|^2-\|\hat{x}-x^{t+1}\|^2).
\end{array}
$$
Noting that
$$
\langle \nabla_x F(x^t,\xi_t),x^{t+1}-x^t\rangle+\frac{1}{2\sigma}\|x^{t+1}-x^t\|^2\geq -\frac{\sigma}{2}\|\nabla_x F(x^t,\xi_t)\|^2\geq -\frac{\sigma}{2}\kappa_f
$$
and
$$
\frac{\sigma}{2}\|(G(x^t,\xi_t)+{\rm D} G(x^t,\xi_t)(\hat{x}-x^t))\|^2\leq \sigma [\nu_g^2+\kappa_g^2R^2],
$$
we obtain
\[
\begin{array}{ll}
 \frac{1}{2\sigma} \left[\|y^{t+1}\|^2-\|y^t\|^2\right]
\leq \sigma [\kappa_f/2+\nu_g^2+\kappa_g^2R^2]+ \langle y^t,G(\hat{x},\xi_t)\rangle
+ \frac{1}{2\sigma}(\|\hat{x}-x^t\|^2-\|\hat{x}-x^{t+1}\|^2).
\end{array}
\]
Summing  over $\{k,k+1,\ldots,k+s-1\}$ and taking the conditional expectation, we obtain from
Lemma \ref{lem:aux4} and $\|y^{k+s}\|-\|y^k\|\leq s \sigma \beta_0$ that
\begin{align*}
\frac{1}{2\sigma} \mathbb{E}\left[\|y^{k+s}\|^2-\|y^k\|^2\,|\, {\cal F}_k\right]
&\ \ \leq (\frac{\kappa_f}{2}+\nu_g^2+\kappa_g^2R^2)\sigma s+ \sum_{t=k}^{k+s-1} \mathbb{E}\left[\langle y^t,G(\widehat x,\xi_t)\rangle\,|\, {\cal F}_k\right]+\frac{1}{2\sigma}\|\widehat x-x^k\|^2\\
&\ \ \leq (\frac{\kappa_f}{2}+\nu_g^2+\kappa_g^2R^2  )\sigma s -\varepsilon_0 \sum_{t=0}^{s-1} \mathbb{E}\left[\|y^{k+t}\|\,|\, {\cal F}_k\right]+\frac{1}{2\sigma}R^2
\\
&\ \ \leq (\frac{\kappa_f}{2}+\nu_g^2+\kappa_g^2R^2  )\sigma s -\varepsilon_0 \sum_{t=0}^{s-1} \mathbb{E}\left[\|y^{k}\|-\sigma\beta_0t \,|\, {\cal F}_k\right]+\frac{1}{2\sigma}R^2\\
&\ \ \leq (\frac{\kappa_f}{2}+\nu_g^2+\kappa_g^2R^2  )\sigma s+ \varepsilon_0\sigma\beta_0 \frac{s(s-1)}{2}
-\varepsilon_0 s \|y^{k}\| +\frac{1}{2\sigma}R^2.
\end{align*}
Further, we get from   (\ref{eq:theta9}) that
\[
\begin{array}{l}
\mathbb{E}\left[\|y^{k+s}\|^2\,|\, {\cal F}_k\right]\\[10pt]
\ \ \leq
\|y^k\|^2+2(\kappa_f/2+\nu_g^2+\kappa_g^2R^2  )\sigma^2 s
+\varepsilon_0\sigma^2\beta_0s(s-1)-2\varepsilon_0\sigma s \|y^{k}\| + R^2\\[10pt]
\ \ \leq(\|y^k\|- \frac{\varepsilon_0\sigma}{2}s)^2
+\varepsilon_0\sigma^2 \beta_0s(s-1)+ 2(\kappa_f/2+\nu_g^2+\kappa_g^2R^2  )\sigma^2 s+ R^2-\varepsilon_0\sigma s \|y^{k}\| \\[10pt]
\ \ \leq(\|y^k\|- \frac{\varepsilon_0\sigma}{2}s)^2
+\varepsilon_0\sigma s\left[\sigma \beta_0(s-1)+ \displaystyle\frac{(\kappa_f+2\nu_g^2+2\kappa_g^2R^2)\sigma}{\varepsilon_0}+\frac{ R^2}{\varepsilon_0\sigma s}- \vartheta (\sigma,s)\right]\\[10pt]
\ \ \leq (\|y^k\|- \frac{\varepsilon_0\sigma}{2}s)^2.
\end{array}
\]
Together with the Jensen inequality and the fact that $\|y^k\|\geq  \vartheta (\sigma,s)\geq \frac{\varepsilon_0\sigma}{2}s$, we have
$
\mathbb{E}\left[\|y^{k+s}\|\,|\, {\cal F}_k\right]\leq
\|y^k\|- \varepsilon_0\sigma s/2.
$
The proof is completed.
\hfill $\Box$\\

If we take
$
\theta=\vartheta (\sigma,s),\ \delta_{\max}=\sigma \beta_0,\ \zeta = \frac{\sigma}{2}\varepsilon_0,\ t_0=s,
$ it follows from Lemma \ref{lem:aux5} that the conditions of Lemma \ref{lem:Yu1}  are satisfied with respect to $\|y^k\|$.
For simplicity,  we define
\[
\Delta_1(\sigma,s):= \kappa_1/(\sigma s)+\kappa_2\sigma+\kappa_3 \sigma s,\qquad \Delta_2 (\sigma,s,\mu):=\Delta_1(\sigma,s)+8\beta_0^2\sigma s \log \left(  1/\mu\right)/\varepsilon_0,
\]
where $\kappa_0,\kappa_1,\kappa_2$ and $\kappa_3$ are constants given by
\[
\kappa_1= R^2/\varepsilon_0,\,\,
\kappa_2= (\kappa_f+2 \nu_g^2+2\kappa_g^2R^2)/\varepsilon_0-\beta_0,\,\,
\kappa_3=2\beta_0 + \varepsilon_0/2+ 8\beta_0^2\log  (32\beta_0^2/\varepsilon_0^2)/\varepsilon_0.
\]

We can also observe that
 $\Delta_1(\sigma,s)$ and $\Delta_2(\sigma,s,\mu)$ are exactly the same as the right-hand sides of  (\ref{eq:aux1}) and (\ref{eq:aux2}) in Lemma  \ref{lem:Yu1}, respectively. Therefore, in view of  Lemma \ref{lem:aux5}, the following lemma is a direct consequence of Lemma  \ref{lem:Yu1}.
\begin{lemma}\label{lem:lambda}
Let the conditions of Lemma \ref{lem:aux5} be satisfied. It follows that
\begin{equation}\label{eq:lambda-Exp}
\mathbb E[\|y^k\|]\leq \Delta_1(\sigma,s).
\end{equation}
Moreover, for any constant $0<\mu<1$, we have
$
{\rm Prob}\,[\|y^k\|\geq \Delta_2(\sigma,s,\mu)]\leq\mu.
$
\end{lemma}

The next lemma is a technical result, which will be used in estimating the difference $l_{\sigma}(x^k,y)-l_{\sigma}(x,y^k)$.
\begin{lemma}\label{lem-4.1}
Let $(x^k,y^k)$ be generated by Algorithm  \ref{ASalm}. Then for $k=1,2,\ldots$ and $x \in X$,
\begin{equation}\label{eq:3.1aa}
\begin{array}{l}
\|\Pi_{{\cal K}^{\circ}}(y^k+\sigma G(x,\xi_k))\|^2-\|\Pi_{{\cal K}^{\circ}}(y^k+\sigma l^k_g(x))\|^2\geq 0.
\end{array}
\end{equation}
\end{lemma}
{\bf Proof}. Let $\eta_1=\|\Pi_{{\cal K}^{\circ}}(y^k+\sigma G(x,\xi_k))\|^2$ and $\eta_2=\|\Pi_{{\cal K}^{\circ}}(y^k+\sigma l^k_g(x))\|^2$, which implies that
\begin{equation}\label{eq:aa2}
\eta_1=\min\{ \|u'\|^2: u'+y^k+\sigma G(x,\xi_k) \in {\cal K}\}\ \mbox{ and }\ \eta_2=\min\{\|u'\|^2: u'+y^k+\sigma l^k_g(x) \in {\cal K}\}.
\end{equation}
In view of the graph-convexity of the set-valued mapping  $x \rightarrow G(x,\xi)-{\cal K}$, we have for $x\in X$,
$$
l_g^k(x)-G(x,\xi_k)={\rm D} G(x^k,\xi_k)(x-x^k)) \in {\cal K}.
$$
Thus, if $u'+y^k+\sigma G(x,\xi_k)\in {\cal K}$, we have
$$
u'+y^k+\sigma l^k_g(x) =u'+y^k+\sigma G(x,\xi_k)+\sigma(l_g^k(x)-G(x,\xi_k)) \in {\cal K},
$$
which implies that the feasible region for $\eta_1$ is a subset of that for  $\eta_2$. Therefore, we have $\eta_1 \geq \eta_2$. The proof is completed.\hfill $\Box$\\

The result in the following proposition gives the bound of the argument Lagrange function at $(x^k,y^k)$  generated by
Algorithm \eqref{ASalm}.
\begin{proposition}\label{prop-4.1}
Let $(x^k,y^k)$ be generated by Algorithm  \ref{ASalm}. Then for $(x,y)\in X \times {\cal Y}$ and $k=1,2,\ldots$,
\begin{equation}\label{eq:3.1a}
\begin{array}{l}
l_{\sigma}(x^k,y)-l_{\sigma}(x,y^k) \\[6pt]
\ \  \leq (f(x^k)-F(x^k,\xi_k))-(f(x)-F(x,\xi_k))+\langle \Pi_{{\cal K}^{\circ}}(y^k+\sigma g(x)),G(x,\xi_k)-g(x)\rangle\\[6pt]
\ \ \quad +\langle \Pi_{{\cal K}^{\circ}}(y+\sigma g(x^k)),g(x^k)-G(x^k,\xi_k)-{\rm D}G(x^k,\xi_k)(x^{k+1}-x^k)\rangle+(\beta_0^2/2+\kappa_f/2++2\nu_g^2)\sigma\\[6pt]
\ \ \quad +\displaystyle \frac{1}{2\sigma}(\|x-x^k\|^2-\|x-x^{k+1}\|^2)+\displaystyle \frac{1}{2\sigma}\left[\|y-y^k\|^2-\|y-y^{k+1}\|^2\right].
\end{array}
\end{equation}
\end{proposition}
{\bf Proof}.
It follows from (\ref{eq:opt-x-1}) and the convexity of $F(\cdot,\xi)$ that for $x \in X$,
\begin{equation}\label{eq:a1}
\begin{array}{l}
F(x^k,\xi_k)-F(x,\xi_k)-\displaystyle \frac{1}{2\sigma}\left[ \|\Pi_{{\cal K}^{\circ}}(y^k+\sigma l_g^k(x))\|^2 \right]\\[6pt]
\ \ \leq -\displaystyle \frac{1}{2\sigma}\|y^{k+1}\|^2-\displaystyle \frac{1}{2\sigma}\|x^{k+1}-x^k\|^2-\langle \nabla_xF(x^k,\xi_k),x^{k+1}-x^k\rangle+\displaystyle \frac{1}{2\sigma}(\|x-x^k\|^2-\|x-x^{k+1}\|^2)\\[6pt]
\ \ \leq -\displaystyle \frac{1}{2\sigma}\|y^{k+1}\|^2+\displaystyle \frac{\sigma}{2}\|\nabla_xF(x^k,\xi_k)\|^2+\displaystyle \frac{1}{2\sigma}(\|x-x^k\|^2-\|x-x^{k+1}\|^2)\\[6pt]
\ \ \leq
-\displaystyle \frac{1}{2\sigma}\|y^{k+1}\|^2+\displaystyle \frac{\sigma}{2}\kappa_f^2+\displaystyle \frac{1}{2\sigma}(\|x-x^k\|^2-\|x-x^{k+1}\|^2).
\end{array}
\end{equation}
Noting from the definition of $y^{k+1}$ that
$
y^{k+1}={\rm argmax}_y \left\{{\rm D}_y l_{\sigma}^k(x^{k+1},y^k)(y-y^k)-\displaystyle \frac{1}{2\sigma}\|y-y^k\|^2\right\},
$
we obtain that for any $y \in {\cal Y}$,
$$
l_{\sigma}^k(x^{k+1},y)\leq l_{\sigma}^k(x^{k+1},y^k)+{\rm D}_y l_{\sigma}^k(x^{k+1},y^k)(y^{k+1}-y^k)-\displaystyle \frac{1}{2\sigma}\left[\|y^{k+1}-y^k\|^2-\|y-y^k\|^2+\|y-y^{k+1}\|^2\right].
$$
From the definition of $l^k_{\sigma}$, we obtain from the above inequality that
$$
\begin{array}{l}
 \|\Pi_{{\cal K}^{\circ}}(y+\sigma l_g^k(x^{k+1}))\|^2-\|y\|^2\\[6pt]
\ \ \leq  \|\Pi_{{\cal K}^{\circ}}(y^k+\sigma l_g^k(x^{k+1}))\|^2-\|y^k\|^2 -\|y^{k+1}-y^k\|^2\\[6pt]
\ \ \quad +2\left[ \langle \Pi_{{\cal K}^{\circ}}(y^k+\sigma l_g^k(x^{k+1}))-y^k, y^{k+1}-y^k\rangle \right]+\|y-y^k\|^2-\|y-y^{k+1}\|^2\\[6pt]
\ \ \leq \|y^{k+1}\|^2-\|y^k\|^2+\|y^{k+1}-y^k\|^2+\|y-y^k\|^2-\|y-y^{k+1}\|^2.
\end{array}
$$
Thus, in view of (\ref{eq:6}), we get
\begin{equation}\label{eq:a3}
\begin{array}{l}
\displaystyle \frac{1}{2\sigma}\left[ \|\Pi_{{\cal K}^{\circ}}(y+\sigma l_g^k(x^{k+1}))\|^2-\|y\|^2 \right]\\[6pt]\ \ \leq
 \displaystyle \frac{1}{2\sigma}[\|y^{k+1}\|^2-\|y^k\|^2]
 +\displaystyle \frac{1}{2\sigma}\left[\|y-y^k\|^2-\|y-y^{k+1}\|^2\right]+\displaystyle \frac{1}{2}\beta_0^2\sigma.
 \end{array}
\end{equation}
 For $(x,y)\in X \times {\cal Y}$, we can express $l_{\sigma}(x^k,y)-l_{\sigma}(x,y^k)$ as follows
 $$
 \begin{array}{rcl}
 l_{\sigma}(x^k,y)-l_{\sigma}(x,y^k) &= &
 F(x^k,\xi_k)-F(x,\xi_k)-\displaystyle \frac{1}{2\sigma}\left[ \|\Pi_{{\cal K}^{\circ}}(y^k+\sigma l_g^k(x))\|^2 \right]\\[6pt]
 && +(f(x^k)-F(x^k,\xi_k))-(f(x)-F(x,\xi_k))\\[6pt]
 && -\displaystyle \frac{1}{2\sigma}\left[\|\Pi_{{\cal K}^{\circ}}(y^k+\sigma G(x,\xi_k))\|^2-\|\Pi_{{\cal K}^{\circ}}(y^k+\sigma l^k_g(x))\|^2\right]\\[12pt]
 && +\displaystyle \frac{1}{2\sigma}\left[\|\Pi_{{\cal K}^{\circ}}(y^k+\sigma G(x,\xi_k))\|^2-\|\Pi_{{\cal K}^{\circ}}(y^k+\sigma g(x))\|^2\right]\\[12pt]
 && +\displaystyle \frac{1}{2\sigma}\left[ \Pi_{{\cal K}^{\circ}}(y+\sigma g(x^k))\|^2-\|y\|^2 \right]+\displaystyle
 \frac{1}{2\sigma} \|y^k\|^2.
  \end{array}
 $$
 Noting that ${\rm D} (1/2 \|\Pi_{{\cal K}^{\circ}}(y)\|^2)=\Pi_{{\cal K}^{\circ}}(y)$ and $\Pi_{{\cal K}^{\circ}}(y)$ is Lipschitz  continuous with modulus $1$, it follows from \cite[Lemma 5.7]{Beck2017} that
 $$
 \begin{array}{l}
\displaystyle \frac{1}{2\sigma}\left[\|\Pi_{{\cal K}^{\circ}}(y^k+\sigma G(x,\xi_k))\|^2-\|\Pi_{{\cal K}^{\circ}}(y^k+\sigma g(x))\|^2\right]\\[8pt]
 \ \ \leq  \displaystyle  \frac{1}{\sigma}\langle \Pi_{{\cal K}^{\circ}}(y^k+\sigma g(x)),\sigma [G(x,\xi_k)-g(x)] \rangle+\displaystyle  \frac{\sigma}{2}\|G(x,\xi_k)-g(x)\|^2.
 \end{array}
 $$
 Thus we obtain from (\ref{eq:a1}), $\|G(x,\xi_k)\|\leq \nu_g$ and $\|g(x)\|\leq \nu_g$ for any $x \in X$ that
 $$
 \begin{array}{rcl}
 l_{\sigma}(x^k,y)-l_{\sigma}(x,y^k)
 &\leq &\displaystyle
 \frac{1}{2\sigma} [\|y^k\|^2-\|y^{k+1}\|^2]+(f(x^k)-F(x^k,\xi_k))-(f(x)-F(x,\xi_k))\\[6pt]
 &&+ \langle \Pi_{{\cal K}^{\circ}}(y^k+\sigma g(x)),G(x,\xi_k)-g(x) \rangle+(2
 \nu_g^2+\kappa_f/2)\sigma\\[12pt]
 && +\displaystyle \frac{1}{2\sigma}(\|x-x^k\|^2-\|x-x^{k+1}\|^2)+\displaystyle \frac{1}{2\sigma}\left[ \|\Pi_{{\cal K}^{\circ}}(y+\sigma g(x^k))\|^2-\|y\|^2 \right].
  \end{array}
 $$
Combining with (\ref{eq:a3}), we obtain
$$
 \begin{array}{rcl}
 l_{\sigma}(x^k,y)-l_{\sigma}(x,y^k) &\leq& (f(x^k)-F(x^k,\xi_k))-(f(x)-F(x,\xi_k))\\[6pt]
 && +  \langle \Pi_{{\cal K}^{\circ}}(y^k+\sigma g(x)), [G(x,\xi_k)-g(x)] \rangle+(\beta_0^2/2+\kappa_f/2+2
 \nu_g^2)\sigma\\[6pt]
 && +\displaystyle \frac{1}{2\sigma}(\|x-x^k\|^2-\|x-x^{k+1}\|^2)+\displaystyle \frac{1}{2\sigma}\left[\|y-y^k\|^2-\|y-y^{k+1}\|^2\right]\\[12pt]
 && +\displaystyle \frac{1}{2\sigma}\left[ \|\Pi_{{\cal K}^{\circ}}(y+\sigma g(x^k))\|^2- \|\Pi_{{\cal K}^{\circ}}(y+\sigma l^k_g(x^{k+1}))\|^2\right]\\[12pt]
 &\leq& (f(x^k)-F(x^k,\xi_k))-(f(x)-F(x,\xi_k))\\[6pt]
 && +  \langle \Pi_{{\cal K}^{\circ}}(y^k+\sigma g(x)), G(x,\xi_k)-g(x) \rangle+(\beta_0^2/2+\kappa_f/2+2
 \nu_g^2)\sigma\\[12pt]
 &&+ \langle \Pi_{{\cal K}^{\circ}}(y+\sigma g(x^k)), g(x^k)- l^k_g(x^{k+1})\rangle\\[6pt]
 && +\displaystyle \frac{1}{2\sigma}(\|x-x^k\|^2-\|x-x^{k+1}\|^2)+\displaystyle \frac{1}{2\sigma}\left[\|y-y^k\|^2-\|y-y^{k+1}\|^2\right],
\end{array}
$$
which implies (\ref{eq:3.1a}). The proof is completed. \hfill $\Box$

Let us analyze convergence properties of the update (\ref{Itformula}), where the error function is defined as
\begin{equation}\label{errorM}
\epsilon_{l_{\sigma}}(z)=[l_{\sigma} (x,y^*)-l_{\sigma}(x^*,y^*)]+[l_{\sigma} (x^*,y^*)-l_{\sigma}(x^*,y)]=l_{\sigma} (x,y^*)-l_{\sigma}(x^*,y),
\end{equation}
where  $z^*=(x^*,y^*)$ is a saddle point of $l_{\sigma}$  defined in (\ref{ACD}).

In the following lemma, we derive a bound of $\epsilon_{l_{\sigma}}(\hat{z}^N)$ in terms of the averaged iterate $\hat{z}^N=(\hat{x}^N,\hat{y}^N)$, where
\begin{equation}\label{hz}
\hat{x}^N=\frac{1}{N}\sum_{k=1}^{N}x^k,\quad \hat{y}^N=\frac{1}{N}\sum_{k=1}^{N}y^k,
\end{equation}
and $N$ is a fixed iteration number.
\begin{proposition}\label{lem:cons}
Let  assumptions {\bf (A1), (C1), (D1)-(D4)} be satisfied. Then,  we have for positive integer $N$,
\begin{equation}\label{exp-gap}
\begin{array}{rcl}
\mathbb E[\epsilon_{l_{\sigma}}(\hat{z}^N)]&
\leq & (\beta_0^2/2+\kappa_f/2+2
 \nu_g^2)\sigma + \kappa_g[\|y^*\|+\nu_g \sigma][\kappa_f+\nu_g\kappa_g\sigma]\sigma \\[6pt] &&+\displaystyle \kappa_g^2[\|y^*\|+\nu_g \sigma]\sigma /N\displaystyle \sum_{k=0}^{N-1}\mathbb E \|y^k\|+(\|y^*\|^2+R^2)/(2\sigma N).
 \end{array}
\end{equation}
\end{proposition}
{\bf Proof}. It follows from (\ref{eq:3.1a}) and {\bf (C1)} that
 \begin{equation}\label{eq:3.1aaa}
\begin{array}{rl}
&\mathbb E[\phi_{\sigma}(x^k,y^*)-\phi_{\sigma}(x^*,y^k) ]\\[6pt]
&\ \ \leq \mathbb E\|y^*+\sigma g(x^k)\|\|{\rm D}G(x^k,\xi_k)\|\|x^{k+1}-x^k\|+(\beta_0^2/2+\kappa_f/2+2
 \nu_g^2)\sigma\\[6pt]
&\ \  \ \ \ \    +\displaystyle \frac{1}{2\sigma}\mathbb E(\|x^k-x^*\|^2-\|x^{k+1}-x^*\|^2)+\displaystyle \frac{1}{2\sigma}\mathbb E\left[\|y^*-y^k\|^2-\|y^*-y^{k+1}\|^2\right]\\[6pt]
&\ \ \leq [\|y^*\|+\sigma \kappa_g]\kappa_g\mathbb E\|x^{k+1}-x^k\|+(\beta_0^2/2+\kappa_f/2+2
 \nu_g^2)\sigma\\[6pt]
 &\ \ \ \ \ \  +\displaystyle \frac{1}{2\sigma}\mathbb E(\|x^k-x^*\|^2-\|x^{k+1}-x^*\|^2)+\displaystyle \frac{1}{2\sigma}\mathbb E\left[\|y^*-y^k\|^2-\|y^*-y^{k+1}\|^2\right].
\end{array}
\end{equation}
Thus we have from Corollary \ref{coro-xd}
that
 \begin{equation}\label{eq:3.1aaa}
 \begin{array}{rl}
&\mathbb E[\phi_{\sigma}(x^k,y^*)-\phi_{\sigma}(x^*,y^k) ]\\[6pt]
& \ \ \leq [\|y^*\|+\sigma \kappa_g]\kappa_g\mathbb E\sigma \Big\{[\kappa_f+\nu_g\kappa_g\sigma]+\kappa_g\|y^k\|\Big\}
+(\beta_0^2/2+\kappa_f/2+2
 \nu_g^2)\sigma\\[6pt]
 &\ \ \ \ \ \  +\displaystyle \frac{1}{2\sigma}\mathbb E(\|x^k-x^*\|^2-\|x^{k+1}-x^*\|^2)+\displaystyle \frac{1}{2\sigma}\mathbb E\left[\|y^*-y^k\|^2-\|y^*-y^{k+1}\|^2\right]\\[6pt]
&\ \ \leq (\beta_0^2/2+\kappa_f/2+2
 \nu_g^2)\sigma + \kappa_g[\|y^*\|+\nu_g \sigma][\kappa_f+\nu_g\kappa_g\sigma]\sigma + \kappa_g^2[\|y^*\|+\nu_g \sigma]\sigma \mathbb E \|y^k\|\\[6pt]
 &\ \ \ \ \ \   +\displaystyle \frac{1}{2\sigma}\mathbb E(\|x^k-x^*\|^2-\|x^{k+1}-x^*\|^2)+\displaystyle \frac{1}{2\sigma}\mathbb E\left[\|y^*-y^k\|^2-\|y^*-y^{k+1}\|^2\right].
\end{array}
\end{equation}
Summing over $\{1,\ldots,N\}$ in (\ref{eq:3.1aaa}), we obtain
$$
\begin{array}{rl}
\displaystyle \frac{1}{N} \sum_{k=0}^{N-1}
\mathbb E[\phi_{\sigma}(x^k,y^*)-\phi_{\sigma}(x^*,y^k) ]&\leq  (\beta_0^2/2+\kappa_f/2+2
 \nu_g^2)\sigma + \kappa_g[\|y^*\|+\nu_g \sigma][\kappa_f+\nu_g\kappa_g\sigma]\sigma\\[6pt]
 &\ \  \quad +\displaystyle \frac{ \kappa_g^2[\|y^*\|+\nu_g \sigma]\sigma }{N}\displaystyle \sum_{k=0}^{N-1}\mathbb E \|y^k\|+\displaystyle \frac{ \|y^*-y^0\|^2+\|x^*-x^0\|^2}{2\sigma N}.
\end{array}
$$
By the Jensen inequality, $y^0=0$ and $\|x^*-x^0\|\leq R$, we obtain that
$
\mathbb E[\epsilon_{l_{\sigma}}(\hat{z}^N)]\leq \displaystyle \frac{1}{N}\sum_{k=1}^{N} \mathbb E[l_{\sigma}(x^k,y^*)-l_{\sigma}(x^*,y^k) ],
$
which implies (\ref{exp-gap}). \hfill $\Box$
\begin{theorem}\label{expected-rate}
Let  assumptions {\bf (A1), (C1), (D1)-(D4)} be satisfied. If we take $\sigma=1/\sqrt{N}$  in Algorithm \ref{ASalm}, where $N > \nu_g^2\max\{\kappa_g^2,1\}$ is a fixed iteration number. Then, we have
\begin{equation}\label{exp-gap-rate}
\mathbb E[\epsilon_{l_{\sigma}}(\hat{z}^N)]\leq c_l\displaystyle \frac{1}{\sqrt{N}}+\displaystyle \frac{\kappa_g^2(\|y^*\|+1)\kappa_2}{N},
\end{equation}
where $c_l=(\beta_0^2/2+\kappa_f/2+2
 \nu_g^2)+\kappa_g(\|y^*\|+1)(\kappa_f+1)+
\kappa_g^2(\|y^*\|+1)(\kappa_1+\kappa_3)+\displaystyle (\|y^*\|^2+R^2)/2$.
\end{theorem}
{\bf Proof}. It follows from (\ref{exp-gap}) and $N > \nu_g^2\max\{\kappa_g^2,1\}$ that
\begin{equation}\label{exp-gap-rate1}\mathbb E[\epsilon_{\phi_{\sigma}}(\hat{z}^N)]\leq  \frac{\beta_0^2/2+\kappa_f/2+2
 \nu_g^2}{\sqrt{N}} + \frac{\kappa_g[\|y^*\|+1][\kappa_f+1]}{\sqrt{N} } +\displaystyle \frac{ \kappa_g^2[\|y^*\|+1] }{N^{3/2}}\displaystyle \sum_{k=0}^{N-1}\mathbb E \|y^k\|+\displaystyle \frac{ \|y^*\|^2+R^2}{2\sigma N}.
\end{equation}
From
(\ref{eq:lambda-Exp}) of Lemma \ref{lem:lambda}, by taking $s=\lceil \sqrt{N}\rceil$, we have
$
\mathbb E \|y^k\| \leq \kappa_1+\kappa_3+\kappa_2/\sqrt{N}.
$
Combining this with (\ref{exp-gap-rate1}), we obtain (\ref{exp-gap-rate}). \hfill $\Box$\\

To develop the high probability guarantee of the update
\eqref{hz} with ${\hat{z}^k}_k$ being generated by Algorithm \ref{ASalm}, we need the following condition.

\begin{itemize}
\item[{\bf (D5)}]
There exists a constant $\nu_f>0$ such that for any $x\in X$,
$
\mathbb E[\exp\{|F(x,\xi)-f(x)\|^2/\nu_f^2\}]\leq\exp\{1\}.
$	
\end{itemize}

\begin{theorem}\label{pexpected-rate}
Let  assumptions {\bf (A1), (C1), (D1)-(D5)} be satisfied. If we take $\sigma=1/\sqrt{N}$  in Algorithm \ref{ASalm}, where $N > \nu_g^2\max\{\kappa_g^2,1\}$ is a fixed iteration number. Then, we have
 \begin{equation}\label{Hprob}
{\rm Prob}\left[\epsilon_{l_{\sigma}}(\hat{z}^N)\leq \varrho_1\displaystyle\frac{1}{\sqrt{N}}+
\varrho_2\displaystyle\frac{\log N}{\sqrt{N}}+
\varrho_3\displaystyle \frac{\log^2N}{\sqrt{N}}+\varrho_4
\displaystyle\frac{\log N}{N}+\varrho_5
\displaystyle\frac{1}{N}\right]\geq 1-\displaystyle \frac{5}{\sqrt[3]{N}},
\end{equation}
where
$$
\begin{array}{l}
\,\varrho_1=(\beta_0^2/2+\kappa_f/2+6\nu_g^2)+(\|y^*\|+1)((\kappa_1+\kappa_3)\kappa_g+\kappa_f+1)\kappa_g+\displaystyle (\|y^*\|^2+R^2)/2,\\[2pt]
\,\varrho_2=2\left(\nu_f+\|y^*\|\nu_g+(\kappa_1+\kappa_3)\nu_g+8\beta_0^2(\|y^*\|+1)\kappa_g^2/\varepsilon_0\right),\\[2pt]
\,\varrho_3=32\beta_0^2\nu_g/\varepsilon_0,\ \  \varrho_4=2\kappa_2\nu_g,\ \ \varrho_5=
(\|y^*\|+1)\kappa_2\kappa_g^2.
\end{array}
$$
 \end{theorem}
 {\bf Proof}. If follows from (\ref{eq:xdd}) that
 $$
\|x^{k+1}-x^k\|\leq \sigma \Big\{[\kappa_f+\nu_g\kappa_g\sigma]+\kappa_g\|y^k\|\Big\}.
$$
From the above inequality and (\ref{eq:3.1a}), we obtain
$$
\begin{array}{l}
 \langle \Pi_{{\cal K}^{\circ}}(y^*+\sigma g(x^k)),g(x^k)-G(x^k,\xi_k)-{\rm D}G(x^k,\xi_k)(x^{k+1}-x^k)\rangle\\[4pt]
\ \ \leq[\|y^*\|+\nu \sigma]\kappa_g\|x^{k+1}-x^k\|+
\langle y^*,g(x^k)-G(x^k,\xi_k\rangle+2\nu_g^2\sigma\\[4pt]
\ \ \leq [(4\nu_g^2+\kappa_f/2+\beta_0^2/2)+(\|y^*\|+\nu_g \sigma)\kappa_g(\kappa_f+\nu_g\kappa_g\sigma)]\sigma\\[4pt]
\ \ \ \ \ +[\|y^*\|+\nu_g \sigma]\kappa_g^2\sigma \|y^k\|+\langle y^*, g(x^k)-G(x^k,\xi_k)\rangle,
\end{array}
$$
which implies that
\begin{align}
\epsilon_{l_{\sigma}}(z^k)
&\ \ \leq [(6\nu_g^2+\kappa_f/2+\beta_0^2/2)+(\|y^*\|+1)(\kappa_f+1)\kappa_g]\sigma+[\|y^*\|+1]\kappa_g^2\sigma \|y^k\|+f(x^k)\nonumber\\
&\quad\
 -F(x^k,\xi_k)+F(x^*,\xi_k)-f(x^*)+\langle y^*, g(x^k)-G(x^k,\xi_k)\rangle+\langle y^k, G(x^*,\xi_k)-g(x^*)\rangle\nonumber\\ \label{eq:3.1ab}
 &\quad\
  +\displaystyle \frac{1}{2\sigma}\left[\|x^*-x^k\|^2-\|x^*-x^{k+1}\|^2\right]+\displaystyle \frac{1}{2\sigma}\left[\|y^*-y^k\|^2-\|y^*-y^{k+1}\|^2\right].
\end{align}
From the Jensen inequality,  $y^0=0$ and $\|x^0-x^*\|\leq R$ , we obtain
\begin{equation}\label{eq:3.1ab}
\begin{array}{rl}
\epsilon_{l_{\sigma}}(\hat{z}^N)&
\leq [(6\nu_g^2+\kappa_f/2+\beta_0^2/2)+(\|y^*\|+1)(\kappa_f+1)\kappa_g]\sigma+ [\|y^*\|+1]\kappa_g^2\displaystyle \frac{\sigma }{N}\sum_{k=1}^{N}\|y^k\|\\[6pt]
& \quad+\displaystyle \frac{\|y^*\|^2+R^2}{2\sigma N}
 +\left(\displaystyle \frac{1}{N} \sum_{k=1}^{N}f(x^k)-\displaystyle \frac{1}{N} \sum_{k=1}^{N}F(x^k,\xi_k)\right)+\left(\displaystyle \frac{1}{N} \sum_{k=1}^{N}F(x^*,\xi_k)-f(x^*)\right)\\[6pt]
&\quad +\displaystyle \frac{1}{N} \sum_{k=1}^{N}[\langle y^*, g(x^k)-G(x^k,\xi_k)\rangle]+\displaystyle \frac{1}{N} \sum_{k=1}^{N}[\langle y^k, G(x^*,\xi_k)-g(x^*)\rangle].
\end{array}
 \end{equation}
 In view of {\bf (D5)}, we have from Lemma \ref{Lem4.1Lan} that for any $\rho>0$,
 \begin{equation}\label{eq:A1}
 {\rm Prob}\left\{\displaystyle \frac{1}{N} \sum_{k=1}^{N}f(x^k)-\displaystyle \frac{1}{N} \sum_{k=1}^{N}F(x^k,\xi_k)\leq \displaystyle \frac{\nu_f \rho}{\sqrt{N}}\right\}\geq 1-\exp\{-\rho^2/3\}
 \end{equation}
 and
 \begin{equation}\label{eq:A2}
 {\rm Prob}\left\{\displaystyle \frac{1}{N} \sum_{k=1}^{N}F(x^*,\xi_k)-\displaystyle \frac{1}{N} \sum_{k=1}^{N}f(x^*)\leq \displaystyle \frac{\nu_f \rho}{\sqrt{N}}\right\}\geq 1-\exp\{-\rho^2/3\}.
 \end{equation}
 Define
 $
 \eta_t=\langle y^*, g(x^t)-G(x^t,\xi_t)\rangle \mbox{ and } \sigma_t=2\|y^*\|\nu_g.
 $
Then $\mathbb E|_{\xi_{[t-1]}} [\eta_t]=0$  and
$E|_{\xi_{[t-1]}}[\exp\{\eta_t^2/\sigma_t^2\}] \leq \exp\{1\}$; namely, the conditions in  Lemma \ref{Lem4.1Lan} are satisfied.
It follows from  Lemma \ref{Lem4.1Lan} that for any $\rho>0$,
\begin{equation}\label{eq:A3}
 {\rm Prob}\left\{\displaystyle \frac{1}{N} \sum_{k=1}^{N}[\langle y^*, g(x^k)-G(x^k,\xi_k)\rangle]\leq \displaystyle
 \frac{2\|y^*\|\nu_g\rho}{\sqrt{N}}\right\}\geq 1-\exp\{-\rho^2/3\}.
\end{equation}
Choosing $s=\lceil \sqrt{N}\rceil$ and $\mu=\exp(-\rho)/N$ in Lemma \ref{lem:lambda}, we get
\begin{equation}\label{eq:A4a}
 {\rm Prob}\left\{\|y^k\|\leq \Delta_3(N,\rho)\right\} \geq 1-\exp(-\rho)/N,
 \end{equation}
 where
 $
 \Delta_3(N,\rho)=\kappa_1+\kappa_3+\kappa_2/\sqrt{N}+8\beta_0^2(\rho+\log N)/\varepsilon_0.
 $
 Thus we obtain
 \begin{equation}\label{eq:A4}
 {\rm Prob}\left\{\frac{\sigma }{N}\sum_{k=1}^{N}\|y^k\|\leq \left(\displaystyle \frac{\kappa_1+\kappa_3}{\sqrt{N}}+ \displaystyle \frac{\kappa_2}{N}+\displaystyle
 \frac{8\beta_0^2}{\varepsilon_0\sqrt{N}}(\rho+\log N)\right)\right\}\geq 1-\exp(-\rho)/N.
 \end{equation}
 Let
 $
 \zeta_t=\langle y^k,G(x^*,\xi_k)-g(x^*)\rangle\ \mbox{ and }\ \chi_t=2\Delta_3(N,\rho)\nu_g.
 $
Then $\mathbb E|_{\xi_{[t-1]}} [\zeta_t]=0$  and
$E|_{\xi_{[t-1]}}[\exp\{\eta_t^2/\chi_t^2\}] \leq \exp\{1\}$ with probability $1-\exp(-\rho)/N$. It follows from  Lemma \ref{Lem4.1Lan} that for any $\rho>0$,
$$
 {\rm Prob}\left\{\displaystyle \frac{1}{N} \sum_{k=1}^{N}\zeta_t\leq \displaystyle \frac{\rho}{N}\displaystyle
  \sqrt{\displaystyle \sum_{t=1}^{N}\chi_t^2}\right\}\geq (1-\exp\{-\rho^2/3\})(1-\exp(-\rho)/N),
$$
which implies that
 \begin{equation}\label{eq:A5}
 {\rm Prob}\left\{\displaystyle \frac{1}{N} \sum_{k=1}^{N}[\langle y^k, g(x^*)-G(x^*,\xi_k)\rangle]\leq \displaystyle
 \displaystyle \frac{\rho}{N} \sqrt{\displaystyle \sum_{t=1}^{N}\chi_t^2}\right\}\geq (1-\exp\{-\rho^2/3\})(1-\exp(-\rho)/N)
\end{equation}
with
\begin{equation}\label{eq:A5aa}
 \begin{array}{l}
 \displaystyle \frac{\rho}{N} \sqrt{\displaystyle \sum_{t=1}^{N}\chi_t^2}=\frac{2\rho\Delta_3(N,\rho)\nu_g}{\sqrt{N}
 }= \displaystyle \frac{2\rho}{\sqrt{N}}\left(\kappa_1+\kappa_3+\displaystyle \frac{\kappa_2}{\sqrt{N}}+
 \displaystyle \frac{8\beta_0^2}{\varepsilon_0}(\log N+\rho)  \right)\nu_g.
 \end{array}
\end{equation}

Combining (\ref{eq:3.1ab}) with (\ref{eq:A1})-(\ref{eq:A5aa}), by choosing $\rho=\log N$, we obtain
\begin{equation}\label{Hproba}
\begin{array}{l}{\rm Prob}\left[\epsilon_{l_{\sigma}}(\hat{z}^N)\leq \varrho_1\displaystyle\frac{1}{\sqrt{N}}+
\varrho_2\displaystyle\frac{\log N}{\sqrt{N}}+
\varrho_3\displaystyle \frac{\log^2N}{\sqrt{N}}+\varrho_4
\displaystyle\frac{\log N}{N}\right]\\[10pt]
\ \ \geq (1-\exp\{-\log^2 N /3\})^4(1-\exp(-\log N)/N)^2\\[6pt]
\ \ > (1-1/\sqrt[3]{N})^4(1-1/N^2)^2
>1-5/\sqrt[3]{N}.
\end{array}
\end{equation}
\hfill $\Box$

From Theorem \ref{expected-rate} and \ref{pexpected-rate},  Algorithm \ref{ASalm} exhibits the same expected sublinear convergence rate  as Algorithm \ref{ASPGalg} and
${\rm O}(\log^2(N)N^{-1/2})$  minimax optimality measure bound with high probability.
 \section{Preliminary Numerical Experiments}\label{Sec:5}
 We  applied the proposed algorithms to solve three classes of problems. Firstly, we tested the SAPS method for stochastic strongly convex-concave minimax problems with three  regularization functions including 1-norm, 2-norm and maximum function. And the impact of the step sizes $\gamma_k$ are analyzed. Secondly, we demonstrated the SAPS method for general stochastic convex-concave minimax  problems and observed the performance  of Algorithm \ref{ASPGalg} with different dimensions and constraints. Thirdly, we used the LSAAL method in Algorithm \ref{ASalm} to solve the multi-class Neyman-Pearson classification, which is equivalent to a finite-sum convex-concave minimax problem. We selected three multi-class classification LIBSVM data sets to report the performance of LSAAL. In this section, all numerical experiments were implemented by MATLAB R2019a on a laptop with Intel(R) Core(TM) i5-6200U 2.30GHz and 8GB memory.

\begin{table}[h]
\caption{Datasets used in the experiments}\label{table1}%
\begin{tabular}{|c|c|c|c|c|}
 \hline
Data Set  & No. of Classes & No. of Data Points: $m$ & No. of Variables: $n$ & Reference                     \\ \midrule
connect-4 & 3              & 67557                   & 126                   & \cite{C2011}\\
covtype   & 7              & 581012                  & 54                    &            \cite{C2011}                   \\
news20    & 20             & 615935                  & 62061                 &     \cite{C2011}                          \\
 \hline
\end{tabular}
\end{table}
\subsection{Stochastic strongly convex-concave minimax
problems}
We tested the SAPS method in Algorithm \ref{ASPGalg} for solving the following stochastic strongly convex-concave minimax
problem
\begin{equation}\label{SSccm1}
  \min_{x \in \Re^n} \max_{y \in \Re^n}\Big\{\mu\vartheta (x)+\mathbb E_{\xi}[\xi^Tx\xi^Ty]-\mu\vartheta(y)\Big\},
\end{equation}
where $\xi$ is a real random vector uniformly distributed in the hypercube $[0, 1]^n$ and $\vartheta:\Re^n \rightarrow \overline \Re$ is a regular term or projection. In our experiments, let $\mu=1$, $n=3$, and $\vartheta$ be chosen as 1-norm, 2-norm and maximum function, respectively.

\begin{figure}[h]
  \centering
  \subfigure[$\vartheta(\cdot)=\|\cdot\|_1$]{
    \label{fig:subfig:10} 
    \includegraphics[width=1.9in]{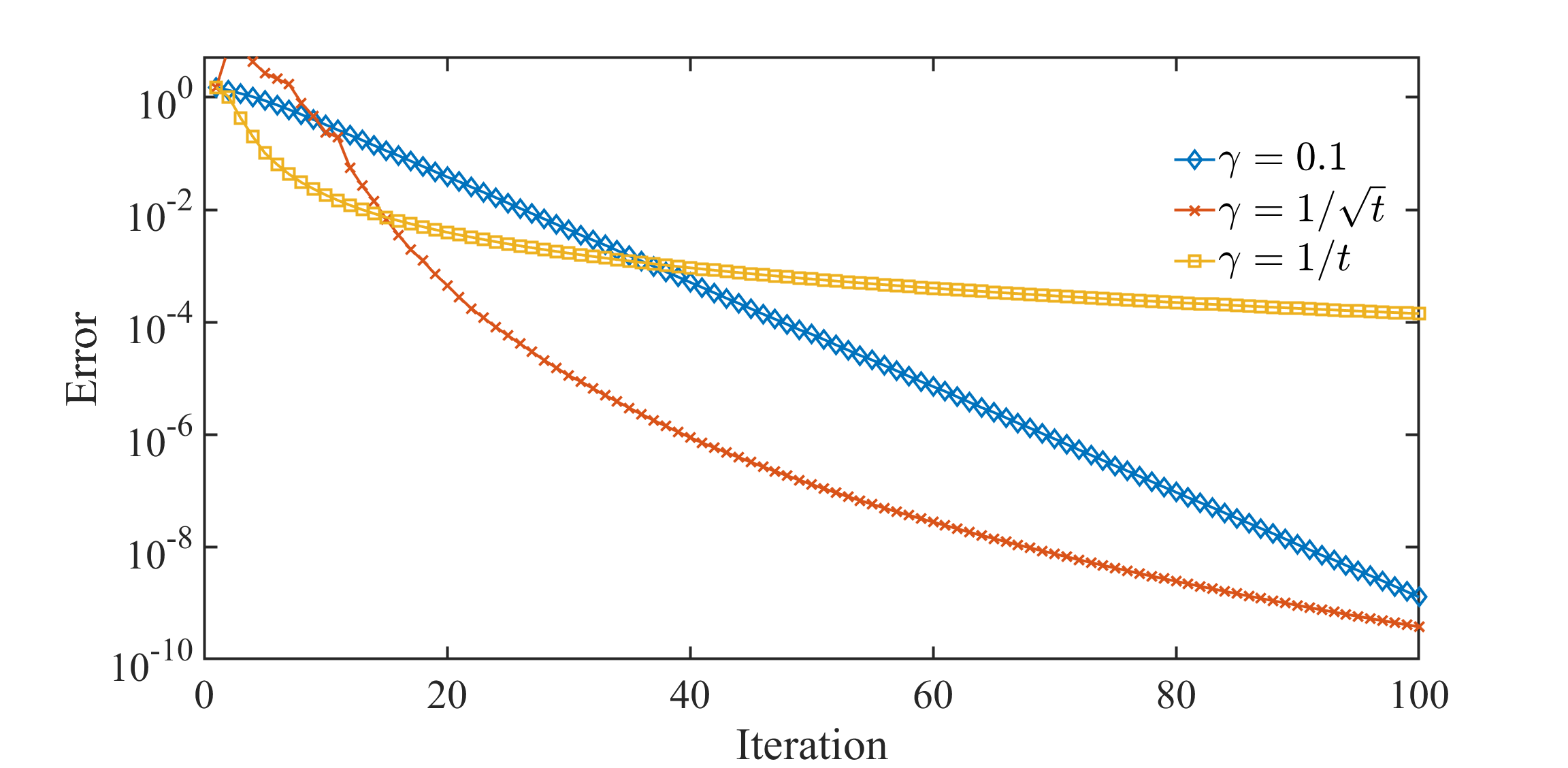}}
  \subfigure[$\vartheta(\cdot)=\|\cdot\|_2$]{
    \label{fig:subfig:100} 
    \includegraphics[width=1.9in]{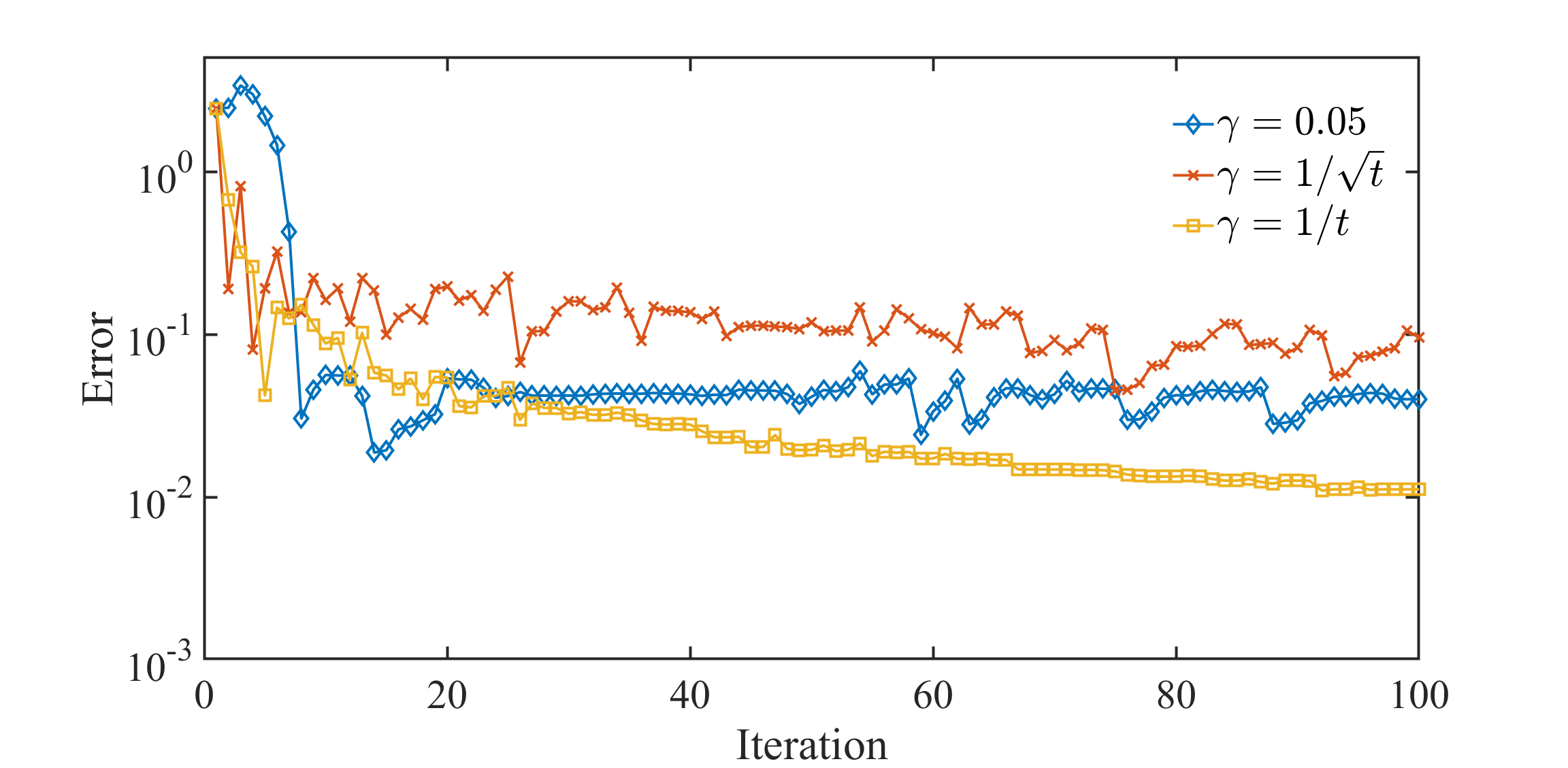}}
    \subfigure[$\vartheta(\cdot)=\max\{\cdot,0\}$]{
    \label{fig:subfig:1000} 
    \includegraphics[width=1.9in]{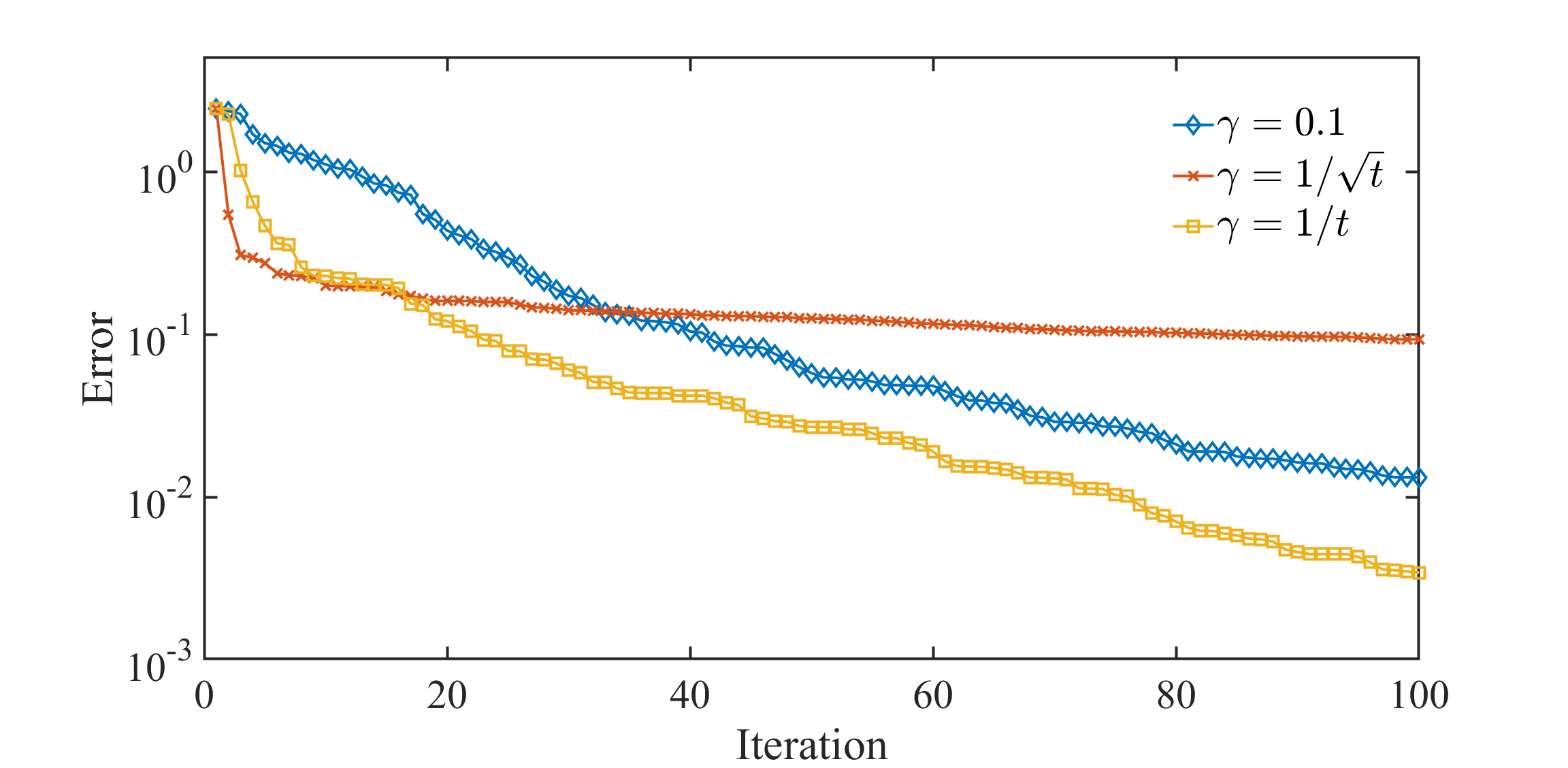}}
  \caption{The trend of the error for \eqref{SSccm1} with respect to iteration.}
  \label{fig:subfig3} 
\end{figure}

For solving \eqref{SSccm1},
we randomly draw a sample $\xi_k$ from the hypercube $[0, 1]^n$ at the $k$-th iteration in Algorithm \ref{ASPGalg} and compute $z^k$ in \eqref{eq:compact-z} as the output of Algorithm \ref{ASPGalg}. In the following experiments, the initial point $z^0:=(x^0,y^0)$ is randomly selected and  the error is computed by
${\rm Error}:=\|z^k-z^*\|_2,$
where $z^*$ is defined as \eqref{eq:compact-z}. We detect the impact of the step sizes $\gamma_k$ on the performance of Algorithm \ref{ASPGalg} with $\gamma_k=1/t,\ 1/\sqrt{t}$ and constant step size under different regular terms and projection. From Figure \ref{fig:subfig3}, we can observe that for minimax problems with 2-norm and maximum function, compared to other step sizes, $\gamma_k=1/t$ allows Algorithm \ref{ASPGalg} to obtain a faster convergence rate. Moreover, as can be seen from the three examples, $\gamma_k=1/t$ makes the convergence more robust, which implies that the variance of the error generated by Algorithm \ref{ASPGalg} is smaller. Therefore, $\gamma_k=1/t$ is a better choice for solving stochastic problems, which further confirms the conclusions numerically.

\subsection{General stochastic convex-concave minimax  problems}
We considered the SAPS method in Algorithm \ref{ASPGalg} for solving a general stochastic convex-concave minimax  problem as follows
\begin{equation}\label{SSccm2}
  \min_{x \in \Re^n} \max_{y \in \Re^n}\Big\{\mu\vartheta (x)+\mathbb E_{u_1,u_2,v_1,v_2}[1-\tanh(v_1\cdot\langle x,u_1\rangle)\tanh(v_2\cdot\langle y,u_2\rangle)]-\mu\vartheta(y)\Big\},
\end{equation}
where the feature vector $u_1,\ u_2$ are  real random vectors uniformly distributed in the hypercube $[0,1]^n$ and $v_1,\ v_2\in\{-1,1\}$ represent the labels of the feature vector $u$ in the binary classification problem. Here we choose $v_1 = {\rm sign}\langle \bar{x},u_1\rangle $ for some $\bar{x}\in \Re^n$ and $v_2 = {\rm sign}\langle \bar{y},u_2\rangle $ for some $\bar{y}\in \Re^n$. Let $\mu=1$ and $\vartheta:\Re^n \rightarrow \overline \Re$ be chosen as 1-norm, 2-norm and maximum function, respectively.

For Algorithm \ref{ASPGalg}, we select the sub-samples $u_1^k,\ u_2^k$
uniformly randomly from $u_i^k\in[0,1]^n, \ i=1,2$ at $k$th iteration and calculate $$v_1^k = {\rm sign}\langle \bar{x},u_1^k\rangle,\ v_2^k = {\rm sign}\langle \bar{y},u_2^k\rangle.$$ Moreover, we compute $\widetilde z^k$ in \eqref{s2.16} as the output of Algorithm \ref{ASPGalg}. In order to verify the performance of our stochastic algorithm, it is necessary to find an approximate optimal solution to \eqref{SSccm2}, since it is difficult to obtain ``true'' optimal solutions for \eqref{SSccm2} in high dimensions. Therefore, \eqref{SSccm2} is approximated by
$$\min_{x \in \Re^n} \max_{y \in \Re^n}\Big\{\mu\vartheta (x)+\sum^m_{i=1}[1-\tanh(v^i_1\cdot\langle x,u^i_1\rangle)\tanh(v^i_2\cdot\langle y,u^i_2\rangle)]-\mu\vartheta(y)\Big\},$$
where $m$ is a sufficiently large sample size. Let $\hat{z}^*=(\hat{x}^*,\hat{y}^*)$ be the optimal solution of the above problem. In the following experiments, set the sample size $m=500$, the initial points $z^0:=(x^0,y^0)$ are randomly selected and  the relative average error is computed by
${\rm Error}:=\|\widetilde z^k-\hat{z}^*\|_2/\|z^0-\hat{z}^*\|_2.$
The step size $\gamma_k$ in Algorithm \ref{ASPGalg} is selected as $\gamma_k=1/\sqrt{t}$.
\begin{figure}[h]
 \centering

 \subfigure[ $n=50$]{
 \begin{minipage}[t]{0.5\linewidth}
 \centering
 \includegraphics[width=0.8\textwidth]{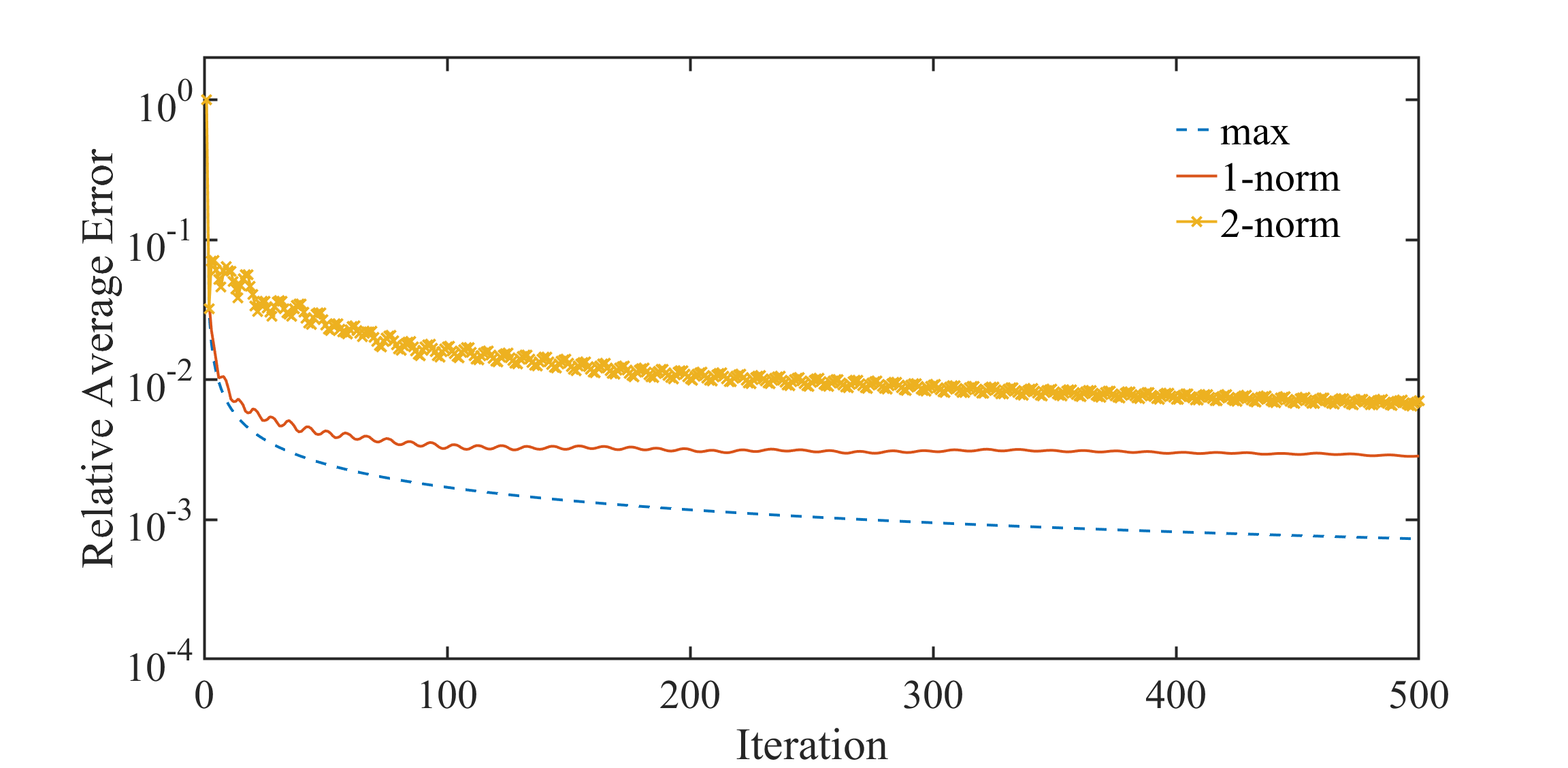}
 \end{minipage}%
 }%
 \subfigure[$n=300$]{
 \begin{minipage}[t]{0.5\linewidth}
 \centering
 \includegraphics[width=0.8\textwidth]{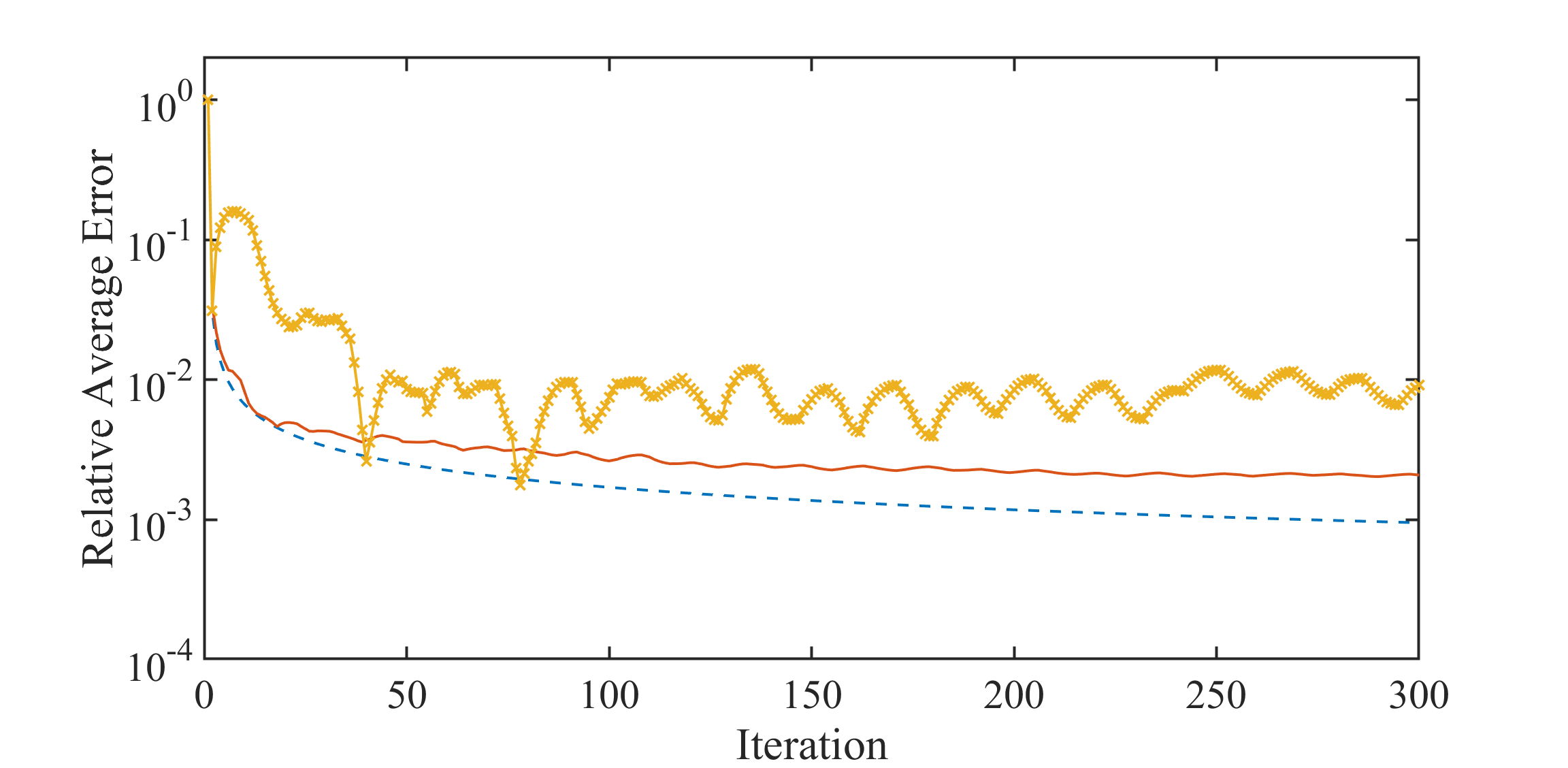}
 \end{minipage}%
 }%

 \subfigure[$n=1000$]{
 \begin{minipage}[t]{0.5\linewidth}
 \centering
 \includegraphics[width=0.8\textwidth]{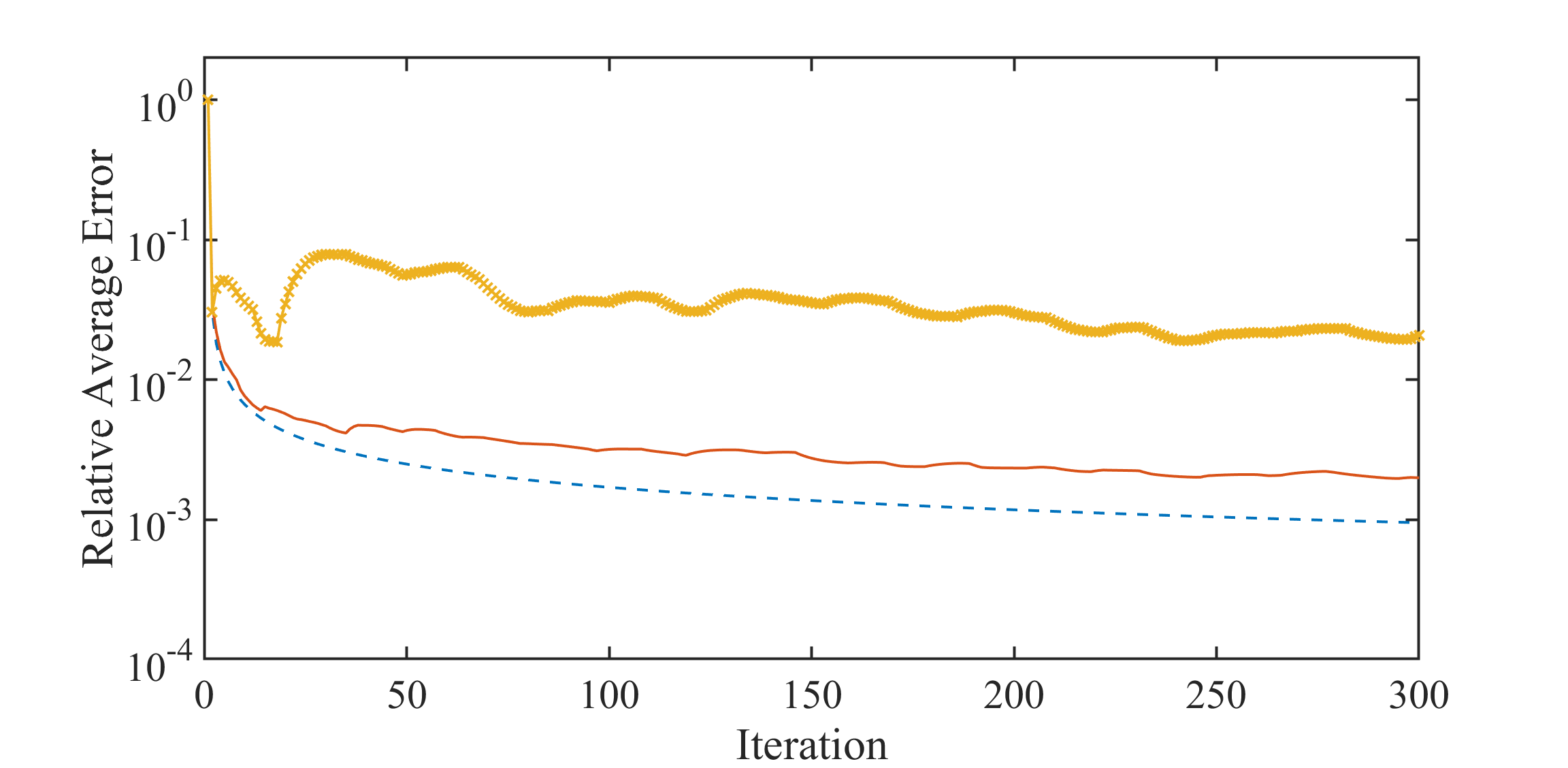}
 \end{minipage}%
 }%
 \subfigure[$n=10000$]{
 \begin{minipage}[t]{0.5\linewidth}
 \centering
 \includegraphics[width=0.8\textwidth]{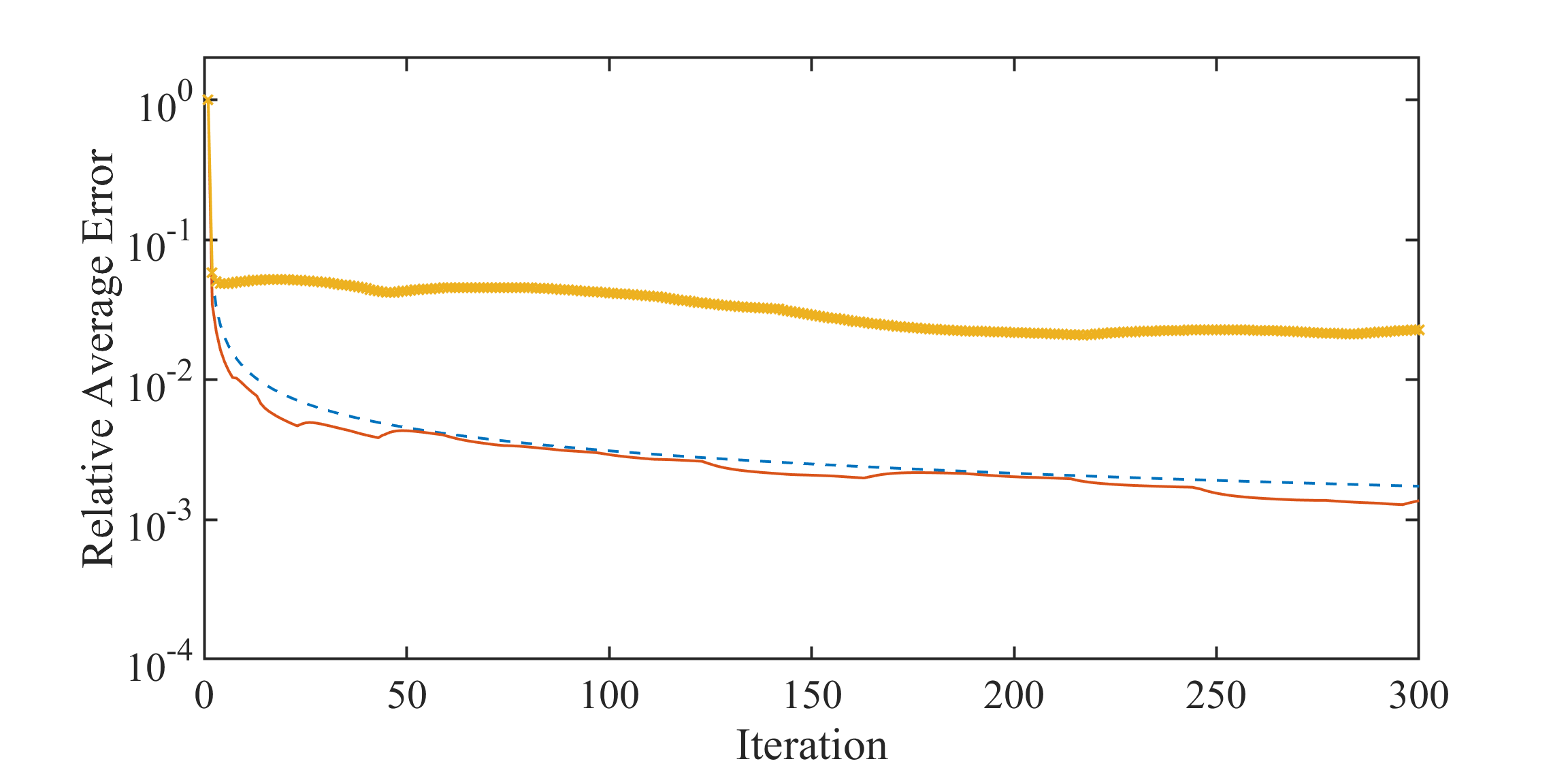}
 \end{minipage}%
 }%

 \centering
 \caption{ The trend of the relative error for solving (\ref{SSccm2}) with respect to iteration}\label{fig13}
 \end{figure}

Figure \ref{fig13} shows the performance of Algorithm \ref{ASPGalg} for (\ref{SSccm2}) with different dimensions and nonsmooth terms. As seen from Figure \ref{fig13}, Algorithm \ref{ASPGalg} has a rapid convergence in early iterations, and then maintains a linear convergence rate.
\subsection{Multi-class Neyman-Pearson classification}
We tested the LSAAL method in Algorithm \ref{ASalm} for the well-known multi-class Neyman-Pearson classification (see \cite{Lin2020} for more details) as follow
\begin{equation}\label{MCNPC1}
\begin{array}{cl}
\displaystyle \min_{\footnotesize{\begin{array}{c}x_i \in \Re^n_i,
\|x_i\|_2\leq\lambda,\\
\forall i=1,2,\ldots,m\end{array}}} & \sum\limits_{l\neq1}\mathbb E[\phi(x_1^T\psi_1-x_l^T\psi_1)]\\[10pt]
\mbox{s. t.} &\sum\limits_{l\neq i}\mathbb E[\phi(x_i^T\psi_i-x_l^T\psi_i)]\leq r_i,\ i=2,3,\ldots ,m,
\end{array}
\end{equation}
where  $\phi:\mathbb{R}\rightarrow\mathbb{R}$ is a non-increasing
convex loss function defined as
$
\phi(x)=\log(1+\exp(-x)).
$
For $m$ classes of data, $\psi_i$
, $i = 1, 2, \ldots, m,$ denotes a random
variable defined using the distribution of data points associated with the $i$-th class and the value of $r_i$
is chosen to capture the misclassification cost of class $i$. Here $\lambda$ is a
regularization parameter. In numerical tests, we compared Algorithm \ref{ASalm} with the following three algorithms on data sets  in Table \ref{table1}.

 \begin{description}
     \item[LAAM]  The deterministic  linearized
approximation augmented Lagrange method.
     \item[LSLM] The  proposed in \cite{YMNeely2017} for stochastic convex programming.
      \item[SALM]   The stochastic augmented
Lagrangian method proposed in \cite{ZZXW22} for stochastic convex programming.
 \end{description}

Let $\psi_i$ follow the empirical distribution over the data set of class $i$ for $i = 1, 2,\ldots, m$, which
implies that all the expectations in \eqref{MCNPC1} become finite-sample averages over data classes. We set  parameters $\lambda = 5$ and $r_i = m-1$ for $i = 2, 3,\ldots, m$.

In order to test the performance of the algorithms more comprehensively,
we use $z^k$ in \eqref{eq:compact-z} and $\hat{z}^N$ in \eqref{hz} as the output of the algorithms to verify the convergence, respectively.
In the following experiments,  the initial point $z^0:=(x^0,y^0)$ is selected as the unit vector.
Due to the difficulty in solving the optimal solution of the stochastic convex programming \eqref{MCNPC1},
we use the KKT condition to verify whether the output of the algorithms satisfies the optimality condition of \eqref{MCNPC1}.
Let $\nabla_zl(z^k)=[\nabla_xl(x^k,y^k),\nabla_yl(x^k,y^k)]$, where $l(x,y)$ is defined in \eqref{eq:Lxi}.
The relative error and the relative average error are computed by
$${\rm RError}:=\frac{\min_{i=1,2,\cdots,k}\nabla_zl(z^i)}{\nabla_zl(z^0)}\
\ \mbox{and}\ \
{\rm RAError}:=\frac{1/k\sum^k_{i=1}\nabla_zl(\hat{z}^i)}{\nabla_zl(z^0)}.$$
In the first numerical experiment with the relative error,
the step sizes are selected as $\gamma_k=1,1,0.5$ in connect-4, covtype,
news20, respectively. In the second one with the relative average error,
the step sizes are selected as $\gamma_k=1,0.5,0.01$.

\begin{figure}[h]
 \centering

\subfigure[connect-4]{ 
    \includegraphics[width=1.9in]{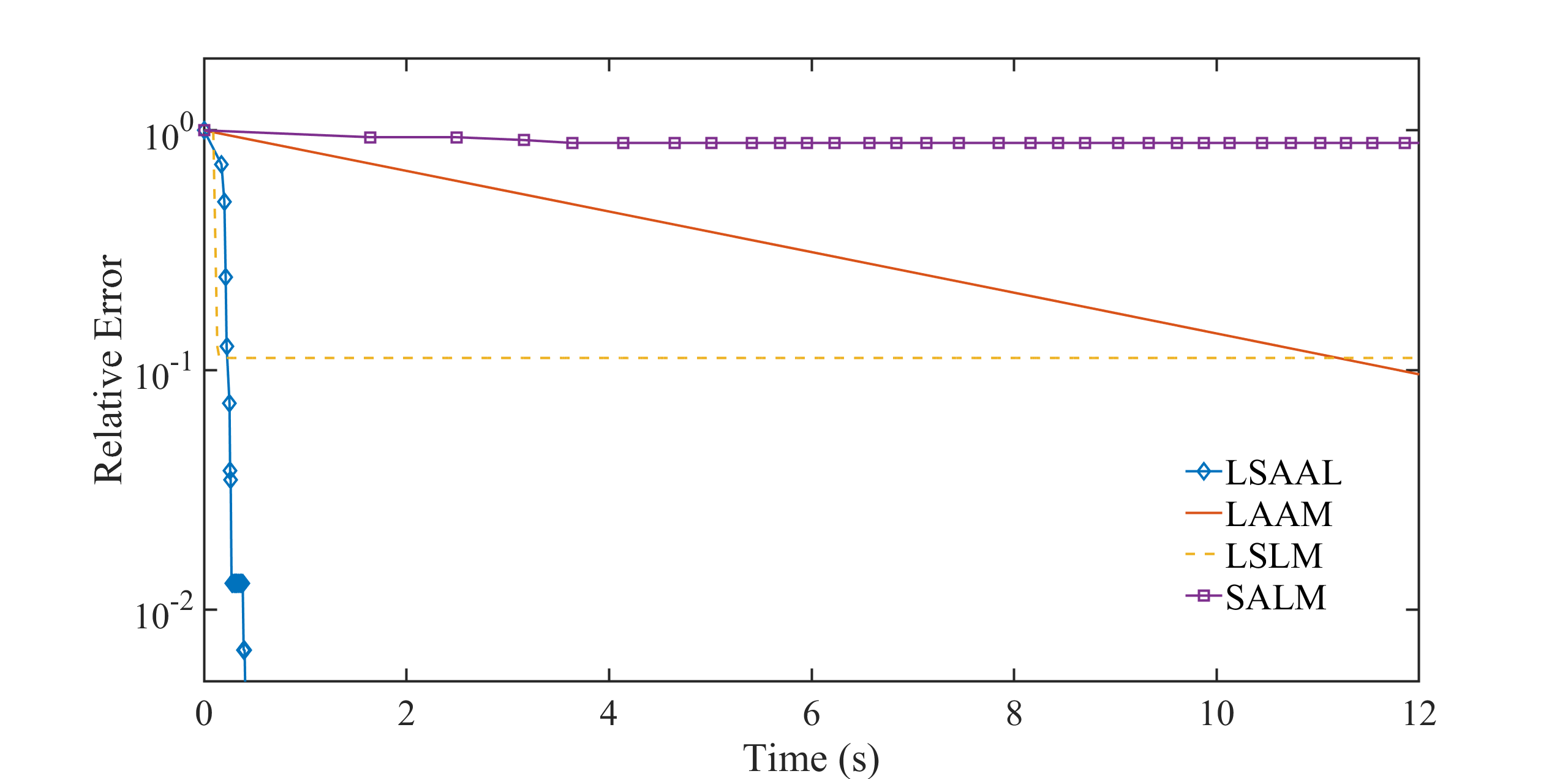}}
  \subfigure[covtype]{
    \includegraphics[width=1.9in]{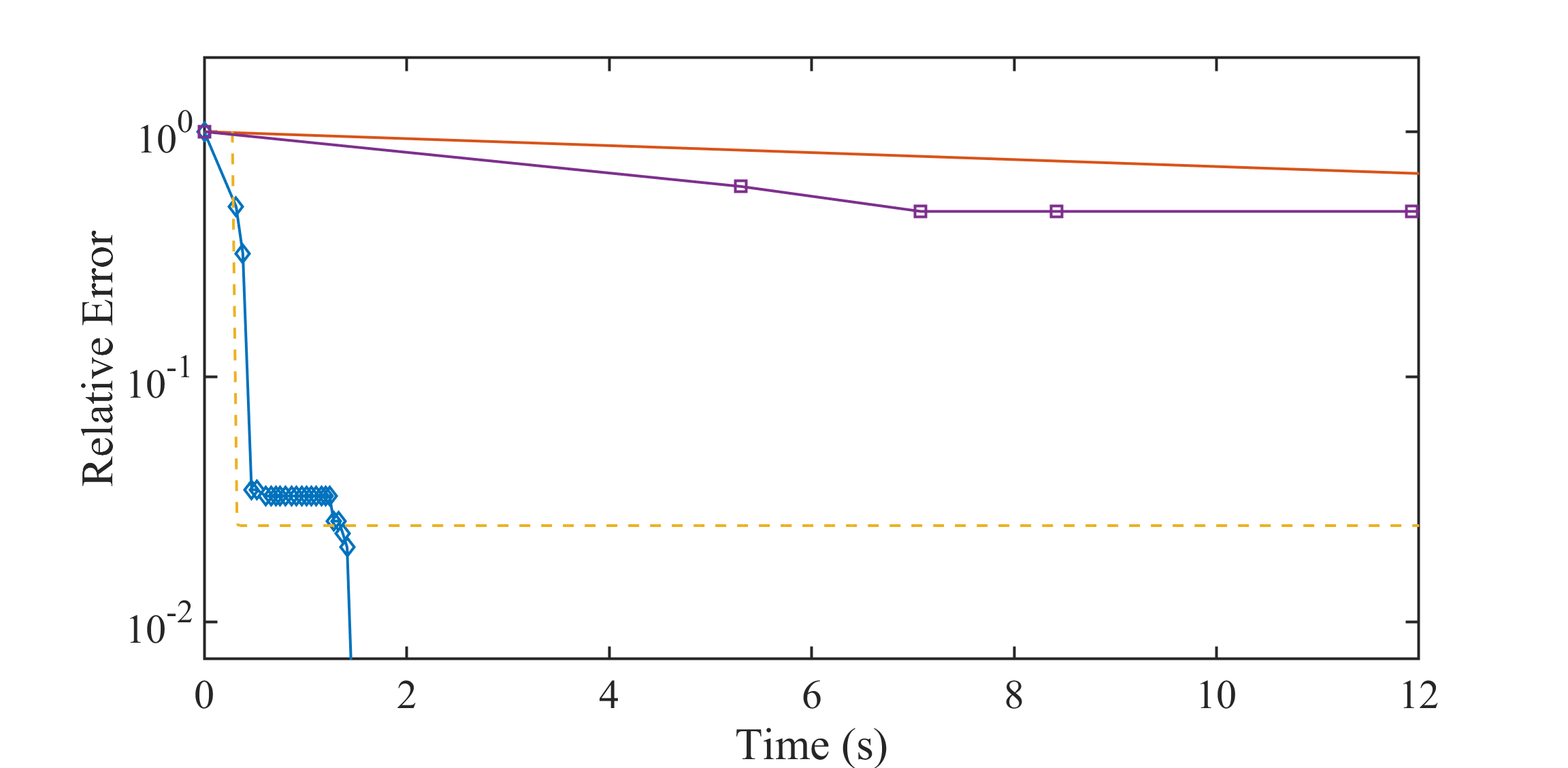}}
    \subfigure[news20]{
    \label{fig:subfig:1000} 
    \includegraphics[width=1.9in]{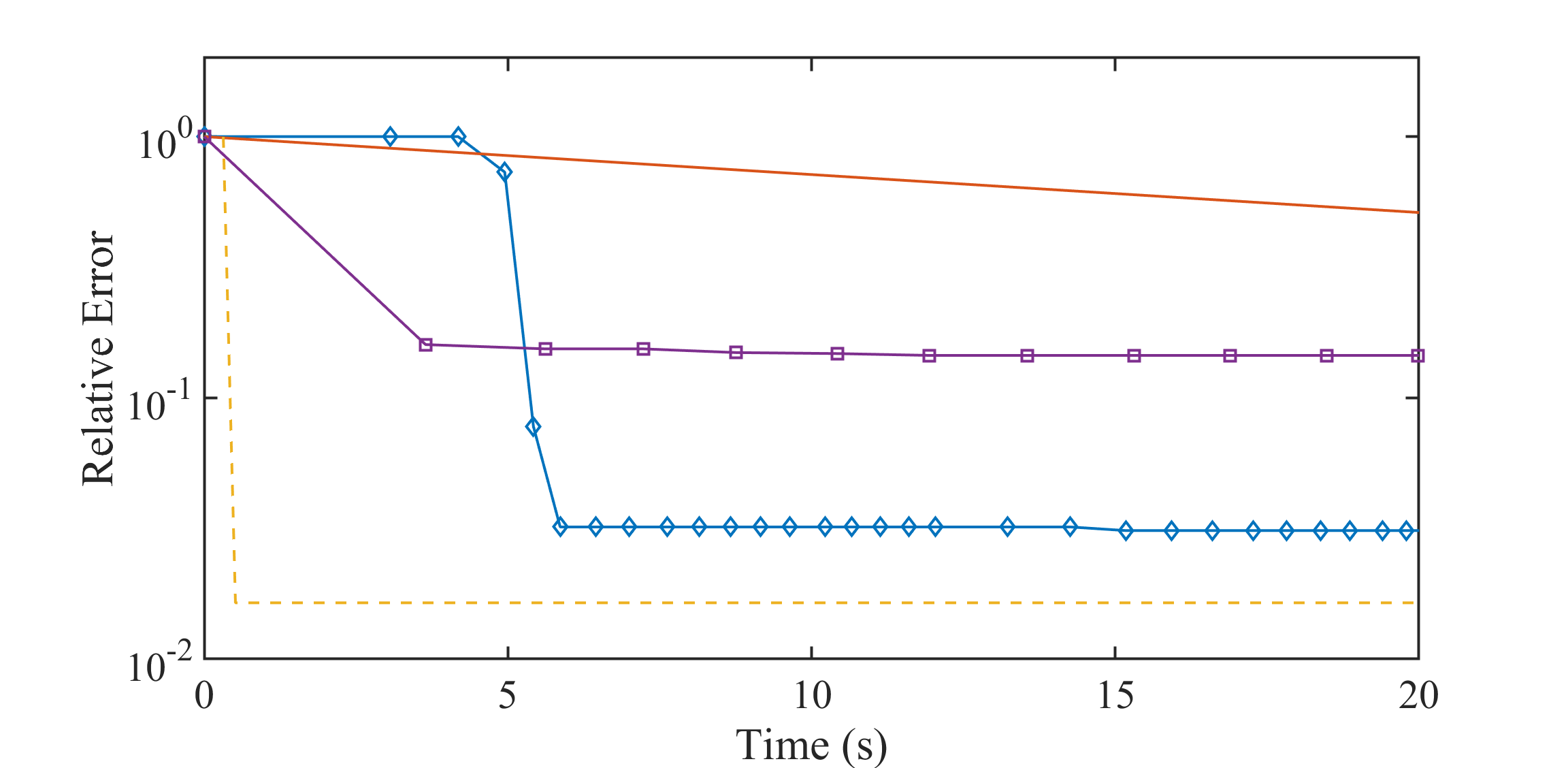}}

\subfigure[connect-4]{ 
    \includegraphics[width=1.9in]{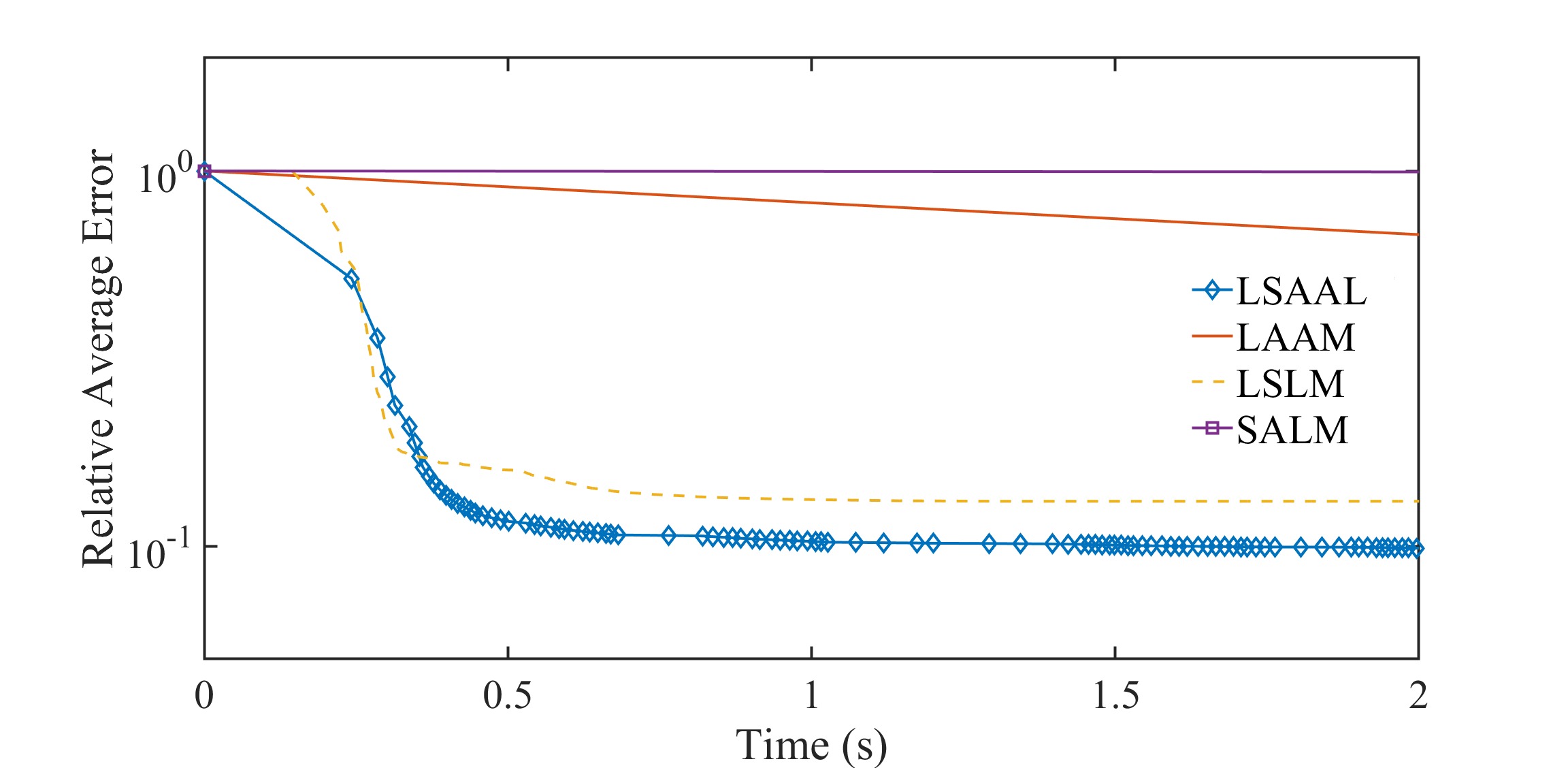}}
  \subfigure[covtype]{
    \includegraphics[width=1.9in]{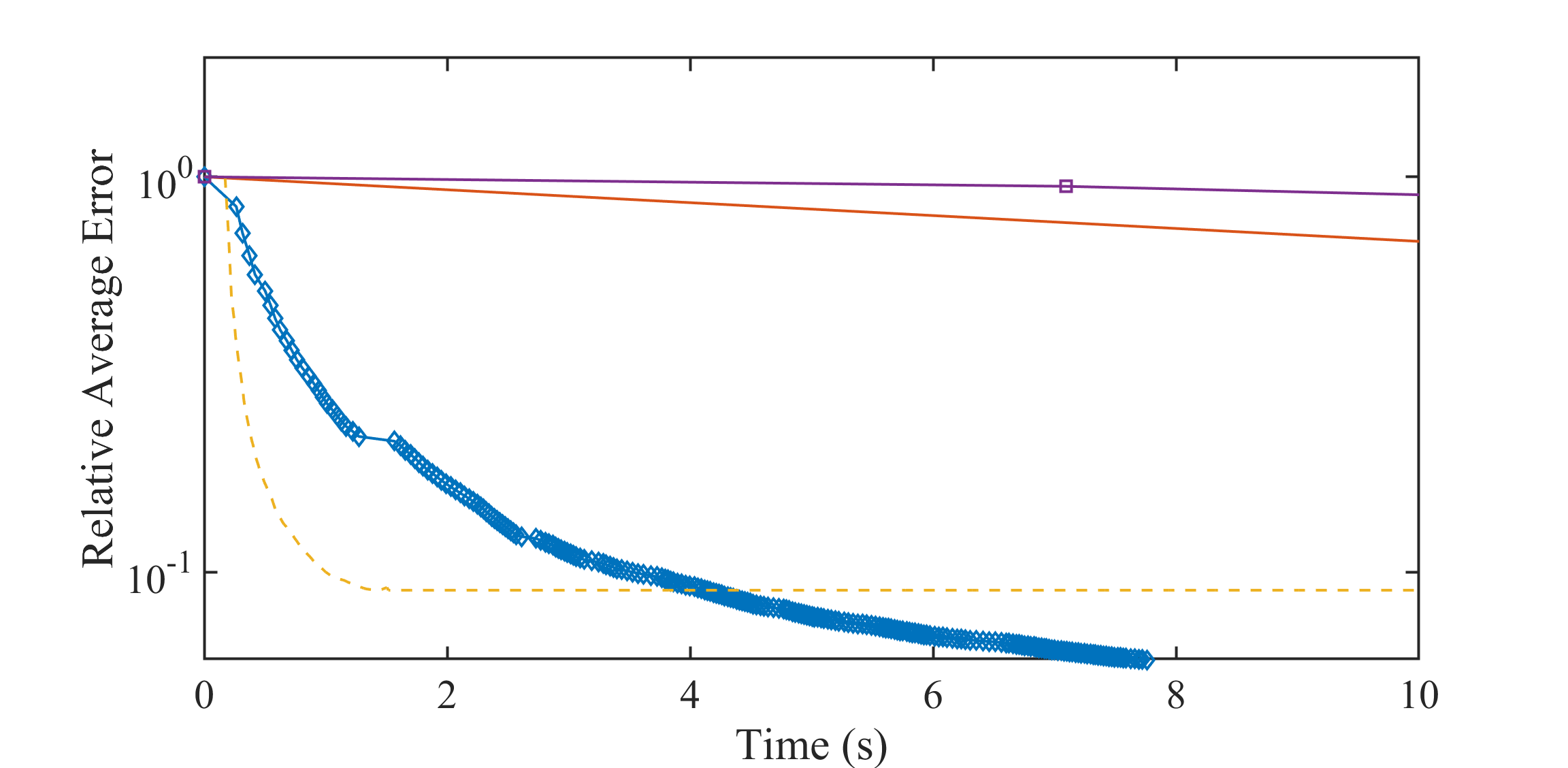}}
    \subfigure[news20]{
    \label{fig:subfig:1000} 
    \includegraphics[width=1.9in]{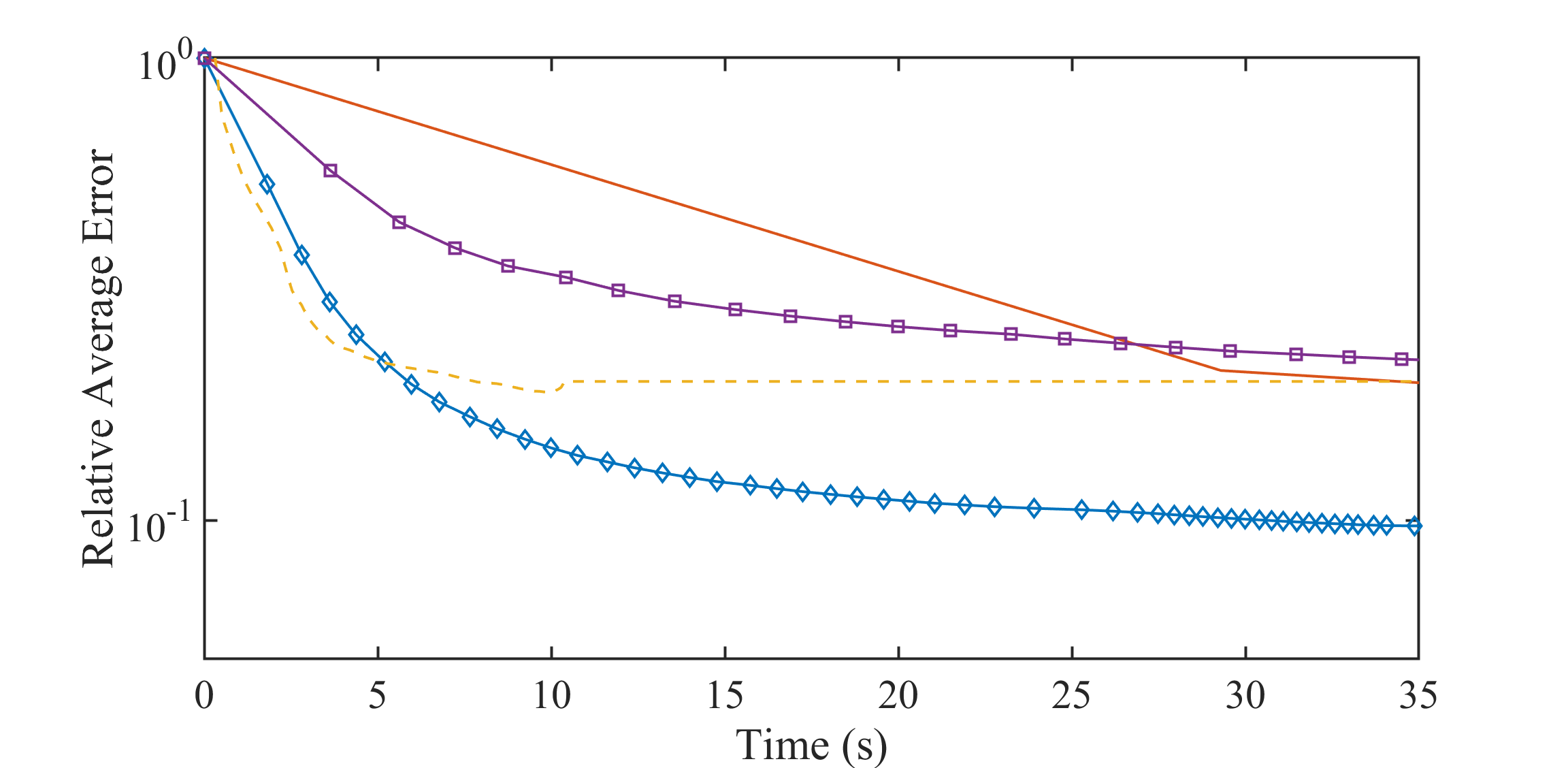}}

 \centering
 \caption{ The trend of the relative error for solving (\ref{MCNPC1}) with respect to CPU time}\label{fig23}
 \end{figure}

 The comparison between four algorithms is presented in Figure \ref{fig23}.
Here, ``Time (s)'' denotes the CPU time in seconds. Generally, LSAAL
performs faster than the other three algorithms.
Besides, the error generated in LSAAL decreases more
rapidly compared with other three algorithms in the early stage of iteration, especially on data set connect-4. This reveals that the combination of
linearization techniques and augmented Lagrange methods indeed does
benefit to the numerical performance of the algorithm.
 \section{Some dicussion}\label{Sec:6}
For the deterministic  minimax optimization, there are many numerical methods in the literature, especially there are many mature numerical algorithms for the convex-concave minimax optimization.
However, for the stochastic  minimax optimization, much attention should be paid to  stochastic approximation  algorithms. In this paper, we have studied the stochastic approximation proximal gradient method for stochastic convex-concave optimization of the form (\ref{minimax}) and the linearized stochastic approximation augmented Lagrangian  method for solving the minimax optimization arising from the stochastic convex conic optimization problems. We have shown that in the expectation sense  the SAPS and LSAAL methods
exhibit  sublinear  convergence rate in terms of the minimax optimality measure if the parameters in the algorithm are properly chosen.  Moreover, the large-deviation properties of the SAPS and LSAAL methods have been established under standard light-tail assumptions. The preliminary numerical experiments have demonstrated that the proposed algorithms are effective for solving the stochastic convex-concave minimax optimization. To the best of our knowledge, SAPS and LSAAL are the first stochastic algorithms for
solving the stochastic nonsmooth convex-concave minimax optimization and stochastic convex conic optimization. We believe that the SAPS and LSAAL methods extends the existing stochastic methods and are used to solve saddle point problems in
 both modern machine learning and tradition research
areas such as saddle point problems, numerical partial differential equations and various types of bi-level optimization problems.

There are some recent  works for the smooth nonconvex-concave  stochastic optimization, see for instance \cite{Bot2020} and \cite{Luo20}. For further work,  an interesting topic is how to use the techniques in this paper to investigate  the stochastic approximation proximal subgradient method to find optimal solutions of stochastic nonsmooth nonconvex-concave  optimization (\ref{minimax}) defined in \cite{Jin2019}.

\section*{Appendix}
\setcounter{equation}{0}
\setcounter{subsection}{0}
\setcounter{lemma}{0}

\renewcommand{\theequation}{A.\arabic{equation}}
\renewcommand{\thesubsection}{A.\arabic{subsection}}

\renewcommand{\thelemma}{A.\arabic{lemma}}

In the proposed algorithms, we have to minimize a strongly convex function at each iteration, the following lemma plays an important role.
\begin{lemma}\label{prox-ineq}
Let ${\cal Z}$ be a finite-dimensional Hilbert space and $\varphi:{\cal Z} \rightarrow \overline \Re$ be a proper lower semicontinuous convex function. Let $z^c \in {\cal Z}$ be given and ${\cal T}: {\cal Z} \rightarrow {\cal Z}$ be a positive definite self-adjoint operator. Then the problem
$$
\min \,\, \varphi (z)+\displaystyle \frac{1}{2}\|z-z^c\|_{\cal T}^2
$$
has a unique solution, denoted by $z^+$. For any $z \in {\cal Z}$,
\begin{equation}\label{Ineq-p}
\varphi (z)+\displaystyle \frac{1}{2}\|z-z^c\|_{\cal T}^2-\displaystyle \frac{1}{2}\|z-z^+\|_{\cal T}^2
\geq \varphi (z^+)+\displaystyle \frac{1}{2}\|z^+-z^c\|_{\cal T}^2.
\end{equation}
\end{lemma}

 By $\xi_{[k]}=(\xi_1,\ldots,\xi_k)$, we denote the history of the process
$\xi_1,\xi_2,\ldots,$  up to time $k$. Unless stated otherwise, all relations between random variables
are supposed to hold almost surely.

The following two lemmas are directly from \cite[Lemma 5]{YMNeely2017} and \cite[Lemma 4.1]{Lan2020}, which are important in the convergence analysis.

\begin{lemma}\label{lem:Yu1}
Let $\{Z_t, t \geq 0\}$ be a discrete time stochastic process adapted to a filtration $\{{\cal F}_t, t\geq
0\}$ with $Z_0 = 0$ and ${\cal F}_0 = \{\emptyset, \Omega\}$. Suppose there exist an integer $t_0 >0$, real constants $\theta>0$, $\delta_{\max}>0$ and $ 0 <\zeta \leq \delta_{\max}$ such that
\[
\begin{array}{rl}
|Z_{t+1}-Z_t| & \leq \delta_{\max},\\[12pt]
\mathbb{E}[Z_{t+t_0}-Z_t\,|\, {\cal F}_t] & \leq \left
\{
\begin{array}{ll}
t_0 \delta_{\max}, & \mbox{if } Z_t < \theta,\\[6pt]
-t_0\zeta, & \mbox{if } Z_t \geq \theta,
\end{array}
\right.
\end{array}
\]
hold for all $t \in \{1,2,\ldots\}.$ Then the following properties are satisfied.
\begin{itemize}
\item[(i)] The  following inequality holds,
\begin{equation}\label{eq:aux1}
\mathbb{E}[Z_t] \leq \theta +t_0 \delta_{\max}+t_0  \frac{4 \delta_{\max}^2}{\zeta}\log \left[  \frac{8 \delta_{\max}^2}{\zeta^2} \right],\ \forall t \in \{1,2,\ldots\}.
\end{equation}
\item[(ii)] For any constant $0 < \mu <1$, we have
    $$
    \Pr\left[Z_t\geq z\right] \leq \mu,\ \forall t \in \{1,2,\ldots\},
    $$
    where
\begin{equation}\label{eq:aux2}
    z=\theta +t_0 \delta_{\max}+t_0  \frac{4 \delta_{\max}^2}{\zeta}\log \left[  \frac{8 \delta_{\max}^2}{\zeta^2}\right]+t_0  \frac{4 \delta_{\max}^2}{\zeta}\log\left( \frac{1}{\mu} \right).
\end{equation}
\end{itemize}
\end{lemma}

\begin{lemma}\label{Lem4.1Lan}
 Let $\xi_{[t]}=\{\xi_1,\ldots, \xi_t\}$ be a sequence of i.i.d. random variables, and
$\zeta_t=\zeta_t(\xi_{[t]})$ be deterministic Borel functions of $\xi_{[t]}$ such that $\mathbb E|_{\xi_{[t-1]}} [\zeta_t]=0$ a.s. and
$E|_{\xi_{[t-1]}}[\exp\{\zeta_t^2/\sigma_t^2\}] \leq \exp\{1\}$ a.s., where $\sigma_t>0$ are deterministic. Then
$$
\forall \lambda \geq 0: {\rm Pr} \left\{\sum_{t=1}^N\zeta_t>\lambda \sqrt{\sum_{t=1}^N \sigma_t^2}\right\}\leq \exp\{-\lambda^2/3\}.
$$
\end{lemma}


\begin{thebibliography}{99}
\bibitem{Akhavan2021}A. Akhavan, M. Pontil and A. Tsybakov, Distributed zero-order optimization under adversarial noise, Advances in Neural Information Processing Systems, 34, 2021, pp. 10209-10220.
\bibitem{Beck2017}
A. Beck, First-Order Methods
in Optimization, Society for Industrial and Applied Mathematics,
Philadelphia, 2017.
\bibitem{Ben2009} A. Ben-Tal, L. E. Ghaoui and A. Nemirovski, Robust Optimization, Princeton University Press, 2009.
    \bibitem{Ben-Nun2019}T. Ben-Nun and T. Hoefler, Demystifying parallel and distributed deep learning: An in-depth  concurrency analysis, ACM Computing Surveys (CSUR), 52:4, 2019, pp. 1-43.
\bibitem{BS2000} J. F. Bonnans and A. Shapiro, {Perturbation Analysis of Optimization Problems},
New York, Springer, 2000.

\bibitem{Bot2020}
R. I. Bo\c{t} and A. B{\"o}hm, Alternating proximal-gradient steps for (stochastic)
nonconvex-concave minimax problems, arXiv preprint arXiv:2007.13605, 2020.
\bibitem{Bottou2018}L. Bottou, F. E. Curtis and J. Nocedal,  Optimization methods for large-scale machine learning, SIAM review, 60:2, 2018, pp. 223-311.

\bibitem{Carmon2020}Y. Carmon, Y. Jin, A. Sidford and K. Tian,  Coordinate methods for matrix games, In 2020 IEEE 61st Annual Symposium on Foundations of Computer Science (FOCS), 2020, pp. 283-293.
\bibitem{C2011}C. C. Chang and C. J. Lin, LIBSVM: A library for support vector machines.
ACM Transactions on Intelligent Systems and Technology, 2:27, 2011, pp. 1-27. Software
available at http://www.csie.ntu.edu.tw/~cjlin/libsvm.
\bibitem{chi2019}Y. Chi, Y. M. Lu and Y. Chen, Nonconvex optimization meets low-rank matrix factorization: An
  overview, IEEE Transactions on Signal Processing, 67:20, 2019, pp. 5239-5269.

\bibitem{Davis2019}D. Davis and D. Drusvyatskiy,
Stochastic model-based minimization of weakly convex functions, SIAM Journal on Optimization, 29:1, 2019, pp. 207-239.

\bibitem{Dietterich1995}T. Dietterich,  Overfitting and undercomputing in machine learning, ACM Computing Surveys (CSUR), 27:3, 1995, pp. 326-327.
\bibitem{Domingos2000}   P. Domingos,  Bayesian averaging of classifiers and the overfitting problem, In International Conference on Machine Learning,  747,  2000, pp. 223-230.
\bibitem{Du2017}    S. S. Du, J. Chen, L. Li, L. Xiao
and D. Zhou, Stochastic variance reduction
methods for policy evaluation, In Proceedings of the 34 th International Conference on Machine
Learning,  70, 2017, pp. 1049-1058.
 \bibitem{Duchi2015}  J. C.  Duchi, M. I. Jordan, M. J. Wainwright and A. Wibisono,  Optimal rates for zero-order convex optimization: The power of two function evaluations, IEEE Transactions on Information Theory, 61:5, 2015, pp. 2788-2806.
 \bibitem{Dvinskikh2022}  D.  Dvinskikh, V. Tominin, Y. Tominin and A. Gasnikov,  Gradient-free optimization for non-smooth minimax problems with maximum value of adversarial noise, arXiv preprint arXiv:2202.06114, 2022.
\bibitem{Farnia2020}F. Farnia and A. Ozdaglar, Train simultaneously, generalize better: Stability of gradient-based minimax learners, International Conference on Machine Learning, 2021, pp. 3174-3185.
\bibitem{Good14}
I. Goodfellow, J. Pouget-Abadie, M. Mirza, B. Xu, D. Warde-Farley, S.
Ozair, A. Courville and Y. Bengio, Generative adversarial nets, In Advances in
Neural Information Processing Systems, 2014, pp. 2672-2680.
\bibitem{Gower19}R.M. Gower, N. Loizou, X. Qian, A. Sailanbayev, E. Shulgin and P. Richt{\'a}rik, Sgd: General analysis and improved rates, International Conference on Machine Learning, 2019, pp. 5200-5209.
\bibitem{Huang2022}
 F. H. Huang, S. Q. Gao, J. Pei and  H. Huang,
Accelerated zeroth-order and first-order momentum
methods from mini to minimax optimization, Journal of Machine Learning Research, 23, 2022, pp. 1-70.
\bibitem{HuangZ2022}K. Huang and S. Zhang,  New first-order algorithms for stochastic variational inequalities, SIAM Journal on Optimization, 32:4, 2022, pp. 2745-2772.
\bibitem{Jin2019}
C. Jin, P. Netrapalli  and M. I. Jordan, What is local optimality in nonconvex-nonconcave minimax
optimization? In International Conference on Machine Learning, 2020, pp. 4880-4889.
%




%

\bibitem{Kolluri2020}J. Kolluri, V. K. Kotte, M. S. B. Phridviraj and S. Razia, Reducing overfitting problem in machine learning using novel $L1/4$ regularization method, In 2020 4th International Conference on Trends in Electronics and Informatics, 48184, 2020, pp. 934-938.
\bibitem{Lan2020}
     { G. Lan,} {First-order and Stochastic Optimization Methods for Machine Learning},
     {Springer Series in the Data Sciences}, {Springer, Cham}, {2020}.
     \bibitem{Lei21}Y. Lei, Z. Yang, T. Yang and Y. Ying, Stability and generalization of stochastic gradient methods for minimax problems, International Conference on Machine Learning, 2021, pp. 6175-6186.
\bibitem{Levy20}D. Levy, Y. Carmon, J. C. Duchi and A. Sidford, Large-scale methods for distributionally robust optimization, Advances in Neural Information Processing Systems, 33, 2020, pp. 8847-8860.
    \bibitem{Li19}X. Li and F. Orabona, On the convergence of stochastic gradient descent with adaptive stepsizes, The 22nd International Conference on Artificial Intelligence and Statistics, 2019, pp. 983-992.
\bibitem{Lin20}  Q. Lin, S. Nadarajah, N. Soheili and T. Yang, A data efficient and feasible level set method for stochastic convex optimization with expectation constraints, Journal of Machine Learning Research, 21:143, 2020, pp. 1-45.
 \bibitem{LinT20}  T. Lin, C. Jin and M. I. Jordan,  Near-optimal algorithms for minimax optimization, In Conference on Learning Theory, 2020, pp. 2738-2779.

    \bibitem{Liu19a}  M. Liu, Y. Mroueh, J. Ross, W. Zhang, X. Cui, P. Das and T. Yang, Towards better understanding of adaptive gradient algorithms in generative
adversarial nets, arXiv preprint arXiv:1912.11940, 2019.

  \bibitem{Liu19b}
S. Liu, S. Lu, X. Chen, Y. Feng, K. Xu, A. Al-Dujaili, M.
Hong and U. Obelilly, Min-max optimization without gradients: convergence and applications to black-box evasion and poisoning attacks, In International conference on machine learning, 2020, pp. 6282-6293.
\bibitem{Luo21}L. Luo, G. Xie, T. Zhang and Z. Zhang, Near optimal stochastic algorithms for finite-sum unbalanced convex-concave minimax optimization, arXiv preprint arXiv:2106.01761, 2021.
\bibitem{Luo20} L. Luo, H. Ye and T. Zhang, Stochastic recursive gradient descent ascent for stochastic nonconvex-strongly-concave minimax problems, arXiv preprint
arXiv:2001.03724, 2020.






%
%
\bibitem{Lan2009}
A. Nemirovski, A. Juditsky, G. Lan and A. Shapiro, Robust stochastic approximation approach to stochastic programming, SIAM Journal on Optimization, 19, 2009, pp. 1574-1609.
\bibitem{Mahdavi(2011)} M. Mahdavi, R. Jin and T. Yang,  Trading regret for efficiency: online convex optimization with long term constraints, {Journal of Machine Learning Research}, {13}:3, 2011, pp. 2503-2528.
\bibitem{Mehrdad(2013)}M. Mahdavi, T. Yang and R. Jin, Stochastic convex optimization with multiple objectives,  in Proceedings of Neural Information Processing Systems, 2013, pp. 1115-1123.

\bibitem{Polyak1990}
   { B. T. Polyak,}
      {New stochastic approximation type procedures},
   {Automat. i Telemekh.},
    {7}, {1990}, pp. {98-107}.
     \bibitem{PJ1992}
   { B. T. Polyak and A. B. Juditsky},
     {Acceleration of stochastic approximation by averaging},
    {SIAM Journal on Control and Optimization},
   {30}:4, {1992}, pp. {838-855}.
\bibitem{Pu21}S. Pu and A. Nedi{\'c}, Distributed stochastic gradient tracking methods, Mathematical Programming, 187, 2021, pp. 409-457.

%
%
%
%
%

\bibitem{Rockafellar76a}
    R. T. Rockafellar, { Monotone operators and the proximal point algorithm},
         SIAM Journal on Control and Optimization, 14, 1976, pp. 877-898.

\bibitem{Schmidt2017}M. Schmidt, N. Le~Roux and F. Bach, Minimizing finite sums with the stochastic average gradient, Mathematical Programming, 162, 2017, pp. 83-112.
    \bibitem{Sery2020}T. Sery and M. Cohen, On analog gradient descent learning over multiple access fading
  channels, IEEE Transactions on Signal Processing, 68, 2020, pp. 2897-2911.
      \bibitem{Shapiro2009}A. Shapiro, D. Dentcheva and A. Ruszczynski, Lectures on Stochastic Programming: Modeling and Theory, MOS-SIAM Series on Optimization, 2021.
\bibitem{Q20}Q. Tran-Dinh, D. Liu and  L. M. Nguyen, Hybrid variance-reduced sgd algorithms for nonconvex-concave minimax problems, arXiv preprint arXiv:2006.15266, 2020.

%
\bibitem{Wai18}
H. Wai, Z. Yang, Z. Wang and M. Hong, Multi-agent reinforcement
learning via double averaging primal-dual optimization, In Advances in Neural Information Processing Systems,  2018, pp. 9649-9660.
\bibitem{Wai19}
H. Wai, M. Hong, Z. Yang, Z. Wang and K. Tang, Variance
reduced policy evaluation with smooth function approximation, In Advances in Neural
Information Processing Systems,  2019, pp. 5784-5795.

\bibitem{Wang2020}
Z. Wang, K. Balasubramanian, S. Ma and M. Razaviyayn,
Zeroth-order algorithms for nonconvex minimax problems with improved complexities,
arXiv preprint arXiv:2001.07819, 2020.
\bibitem{Xu2021}Z. Xu, J. Shen, Z. Wang and Y. Dai,  Zeroth-order alternating randomized gradient projection algorithms for general nonconvex-concave minimax problems, arXiv preprint arXiv:2108.00473, 2021.
\bibitem{Xu2023} Z. Xu, H. Zhang, Y. Xu and G. Lan,   A unified single-loop alternating gradient projection algorithm for nonconvex-concave and convex-nonconcave minimax problems,  Mathematical Programming, 2023, pp. 1-72.
\bibitem{Xu2020}
T. Xu, Z. Wang, Y. Liang and H. V. Poor, Enhanced first and zeroth order
variance reduced algorithms for min-max optimization, arXiv preprint arXiv:2006.09361,
2020
\bibitem{Xian21}W. Xian, F. Huang, Y. Zhang and H. Huang, A faster decentralized algorithm for nonconvex minimax problems, Advances in Neural Information Processing Systems, 34, 2021, pp. 25865-25877.
\bibitem{YY20}Y. Yan, Y. Xu, Q. Lin, W. Liu and T. Yang, Optimal epoch stochastic gradient descent ascent methods for min-max  optimization, Advances in Neural Information Processing Systems, 33, 2020, pp. 5789-5800.
\bibitem{Yang2020}J. Yang, N. Kiyavash and N. He, Global convergence and variance-reduced optimization for a class of nonconvex-nonconcave minimax problems, arXiv preprint arXiv:2002.09621,  2020.
\bibitem{Ying2019}X. Ying, An overview of overfitting and its solutions, In Journal of Physics: Conference series, 1168:022022, 2019.
 \bibitem{YMNeely2017}
H. Yu, M. Neely and X. Wei,  Online convex optimization with stochastic constraints,
Advances in Neural Information Processing Systems, 30, 2017.
\bibitem{ZZXW22}L. Zhang, Y. Zhang, X. Xiao and J. Wu, Stochastic approximation proximal method of multipliers for convex stochastic programming, Mathematics of Operations Research, 48:1, 2023, pp. 177-193.
    \bibitem{ZAG21}X. Zhang, N. S. Aybat and M.  G{\"u}rb{\"u}zbalaban, Robust accelerated primal-dual methods for computing saddle points,  arXiv preprint arXiv:2111.12743,  2021.
\end{thebibliography}
\end{document}